\documentclass[a4paper]{article}

\usepackage{latexsym,amssymb,enumerate}
\usepackage{psboxit}
\usepackage[centertags,leqno]{amsmath}
\usepackage{amsthm}
\usepackage{natbib}
\usepackage{epsfig}

\newtheorem{theorem}{Theorem}

\newtheorem{lemma}{Lemma}[section]

\newtheorem{assumption}{Assumption}[section]
\newtheorem{remark}[lemma]{Remark}

\newcommand{\bs}{\bar {s}}%
\newcommand{\bS}{\bar{S}}%
\newcommand{\QY}{\bar{Q}_Y}
\newcommand{\barQY}{\bar{Q}_Y}
\newcommand{\ohneQY}{{Q}_Y}
\newcommand{\bm}{\bar{m}}

\newcommand{\Dpos}{\D}

\newcommand{\Isl}{\ensuremath{\MCI}}
\newcommand{\VIM}{\ensuremath{\MCV}}
\newcommand{\mycomment}[1]{}
\usepackage{bbm}
\usepackage{amsmath}
\usepackage{ifthen}
\usepackage{calc}

 \providecommand{\C}{{\ensuremath{\mathbf{C}}}}

 \providecommand{\D}{{\ensuremath{\mathbf{D}}}}
 
 \providecommand{\E}{{\ensuremath{\mathbf{E}}}}

 \providecommand{\I}{{\ensuremath{\mathbf{I}}}}

 \providecommand{\N}{{\ensuremath{\mathbbm{N}}}}

 \renewcommand{\P}{{\ensuremath{\mathbf{P}}}}
 \providecommand{\P}{{\ensuremath{\mathbf{P}}}}
 
 \providecommand{\R}{{\ensuremath{\mathbbm{R}}}}

 \providecommand{\Z}{{\ensuremath{\mathbbm{Z}}}}
 \providecommand{\1}{{\ensuremath{\mathbbm{1}}}}

 \providecommand{\MCF}{{\ensuremath{\mathcal F}}}

 \providecommand{\MCI}{{\ensuremath{\mathcal I}}}

 \providecommand{\MCV}{{\ensuremath{\mathcal V}}}

 \providecommand{\Yh} {{\ensuremath{\hat{Y}}}}

 \providecommand{\Yt} {{\ensuremath{\tilde{Y}}}}

 \providecommand{\ct} {{\ensuremath{\tilde{c}}}}

 \providecommand{\gt} {{\ensuremath{\tilde{g}}}}

 \providecommand{\xt} {{\ensuremath{\tilde{x}}}}
 \providecommand{\yt} {{\ensuremath{\tilde{y}}}}

 \providecommand{\Qb} {{\ensuremath{\bar{Q}}}}
 
 \providecommand{\Sb} {{\ensuremath{\bar{S}}}}

 \providecommand{\cb} {{\ensuremath{\bar{c}}}}

 \providecommand{\mb} {{\ensuremath{\bar{m}}}}

 \renewcommand{\sb} {{\ensuremath{\bar{s}}}}

 \providecommand{\wb} {{\ensuremath{\bar{w}}}}

\providecommand{\ldt}  {{\ensuremath{\tilde{\ld}}}}

\providecommand{\Psit}{\ensuremath{\tilde{\Psi}}}%

\providecommand{\qed}{\hfill\mbox{$\Box $}}

\providecommand{\ru}[1]    {{{(#1)}}}
\providecommand{\rub}[1]   {{{\bigl(#1\bigr)}}}
\providecommand{\rubb}[1]  {{{\biggl(#1\biggr)}}}
\providecommand{\ruB}[1]   {{{\Bigl(#1\Bigr)}}}

\providecommand{\eckb}[1]  {{{\bigl[#1\bigr]}}}
\providecommand{\eckbb}[1] {{{\biggl[#1\biggr]}}}
\providecommand{\eckB}[1]  {{{\Bigl[#1\Bigr]}}}

\providecommand{\curlb}[1]  {{{\bigl\{#1\bigr\}}}}

\providecommand{\curlB}[1]  {{{\Bigl\{#1\Bigr\}}}}

\providecommand{\abs}[1]  {{\ensuremath{|#1|}}}
\providecommand{\absb}[1] {{\ensuremath{\bigl|#1\bigr|}}}
\providecommand{\absbb}[1]{{\ensuremath{\Bigl|#1\Bigr|}}}
\providecommand{\absB}[1] {{\ensuremath{\Bigl|#1\Bigr|}}}

\providecommand{\al}      {{\ensuremath{\alpha}}}
\providecommand{\ga}      {{\ensuremath{\gamma}}}

\providecommand{\ld}      {{\ensuremath{\lambda}}}

\providecommand{\eps}     {{\ensuremath{\varepsilon}}}
\providecommand{\dl}      {{\ensuremath{\delta}}}

\providecommand{\om}      {{\ensuremath{\omega}}}

\providecommand{\limeps}{{\ensuremath{{\displaystyle \lim_{\eps \ra 0}}}}}

\providecommand{\limn}  {{\ensuremath{{\displaystyle \lim_{n \ra \infty}}}}}

\providecommand{\limt}  {{\ensuremath{{\displaystyle \lim_{t \ra \infty}}}}}

\providecommand{\limT}  {{\ensuremath{{\displaystyle \lim_{T \ra \infty}}}}}

\providecommand{\limyO} {{\ensuremath{{\displaystyle \lim_{y \ra 0}}}}}
\providecommand{\limz}  {{\ensuremath{{\displaystyle \lim_{z \ra \infty}}}}}

\providecommand{\limld} {{\ensuremath{{\displaystyle \lim_{\lambda \ra \infty}}}}}
\providecommand{\limldO}{{\ensuremath{{\displaystyle \lim_{\lambda \ra 0}}}}}
\providecommand{\limb}{{\ensuremath{{\displaystyle \lim_{b \ra \infty}}}}}

\providecommand{\liminft}  {{\ensuremath{{\displaystyle \liminf_{t \ra \infty}}}}}

\providecommand{\limsupn}  {{\ensuremath{{\displaystyle \limsup_{n \ra \infty}}}}}

\providecommand{\tlimn}{\lim_{n \ra \infty}}

\providecommand{\tliminft}  {\liminf_{t \ra \infty}}

\providecommand{\tlimsupt}  {\limsup_{t \ra \infty}}

\providecommand{\qasn}  {\ensuremath{\quad\text{as }n\to\infty}}

\providecommand{\qqast}  {\ensuremath{\qquad\text{as }t\to\infty}}

\providecommand{\qqasldO}{\ensuremath{\qquad\text{as }\ld\to0}}

\providecommand{\qqastk}{\ensuremath{\qquad(\text{as }t\to\infty})}

\providecommand{\ra}{\rightarrow}

\providecommand{\upa}{\uparrow}

\providecommand{\downa}{\downarrow}

\providecommand{\lra}{\longrightarrow}

\providecommand{\lran} {\xrightarrow{n  \ra\infty}}
\providecommand{\lrat} {\xrightarrow{t  \ra\infty}}

\providecommand{\lram} {\xrightarrow{m  \ra\infty}}

\providecommand{\lraldO}{\xrightarrow{\ld\ra 0}}

\providecommand{\lrayO}{\xrightarrow{y\ra0}}

\providecommand{\wlim}{\xrightarrow{\,w\,}}

\providecommand{\wlimeps}%
      {{\ensuremath{\stackrel{\eps \rightarrow \infty}%
                            {\Longrightarrow}}}}

\providecommand{\Var}{\operatorname{Var}}

\providecommand{\supp}{\operatorname{supp}}

\providecommand{\mal}{\ensuremath{{\displaystyle \cdot}}}

\providecommand{\fa}{\ensuremath{\;\;\forall\;}}

\providecommand{\Law}[2][]{{\ensuremath{\mathcal L^{#1}\left(#2\right)}}}
\providecommand{\Lawb}[2][]{{\ensuremath{\mathcal L^{#1}\bigl(#2\bigr)}}}

\providecommand{\clearemptydoublepage}%
    {\newpage{\pagestyle{empty}\cleardoublepage}}

\providecommand{\eqd}{{\ensuremath{\;\overset{d}{=}\;}}}

\newcommand{\cadlag}{c\`adl\`ag}

\begin{document}

\title{\bfseries The Virgin Island Model}

\author{Martin Hutzenthaler\footnote{Research supported by the DFG in the Dutch German Bilateral
Research Group "Ma\-the\-ma\-tics of Random Spatial Models from
Physics and Biology" (FOR 498)} \footnote{Research supported by EPSRC Grant no GR/T19537/01}\\
Johann Wolfgang Goethe-University Frankfurt\\
Rober-Mayer-Stra\ss e 10\\
60325 Frankfurt, Germany\\
Email: hutzenth@math.uni-frankfurt.de}

\date{}

\maketitle
\def\thefootnote{\fnsymbol{footnote}}
\def\@makefnmark{\hbox to\z@{$\m@th^{\@thefnmark}$\hss}}
\footnotesize\rm\noindent
\footnote[0]{{\it AMS\/\ {\rm 2000} subject
classifications} 60K35, 92D25}
\normalsize\rm

\begin{abstract}
  A continuous mass population model with local competition
  is constructed
  where every emigrant colonizes an unpopulated island.
  The population founded by an emigrant is
  modeled as excursion from zero of an one-dimensional
  diffusion.
  With this excursion measure, we construct a process which we call
  Virgin Island Model.
  A necessary and sufficient condition for extinction of the
  total population is obtained for finite initial total mass.
\end {abstract}

\noindent
\section{Introduction}%
\label{sec:introduction}
This paper is motivated by  an open question on
a system of interacting locally regulated diffusions.
In~\cite{HW07}, a sufficient condition
for local extinction is established for such a system.
In general, however, there is no criterion available for
global extinction, that is, convergence of the total mass process
to zero when started in finite total mass.

The method of proof for the local extinction result in~\cite{HW07}
is a comparison with a mean field model $(M_t)_{t\geq0}$ which solves
\begin{equation}  \label{eq:def:M}
  dM_t =\kappa(\E M_t-M_t)dt+h(M_t)dt+\sqrt{2g(M_t)}dB_t
\end{equation}
where $(B_t)_{t\geq0}$ is a standard Brownian motion and
where $h,g\colon[0,\infty)\to\R$ are suitable functions satisfying
$h(0)=0=g(0)$.
This mean field model arises as the limit
as $N\to\infty$ (see Theorem 1.4 in~\cite{Szn89} for the case $h\equiv 0$)
of the following system of interacting locally regulated diffusions
on $N$ islands with uniform migration
\begin{equation}  \begin{split}\label{eq:def:N_islands_model}
  dX_t^N(i)=&\kappa\eckb{\frac{1}{N}
          \sum_{j=0}^{N-1}X_t^N(j)-X_t^N(i)}\,dt\\
         & +h\rub{X_t^N(i)}dt
         +\sqrt{2g\rub{X_t^N(i)}}\,dB_t(i)\qquad i=0,\ldots,N-1.
\end  {split}     \end  {equation}
For this convergence, $X_0^N(0),\dots,X_0^N(N-1)$ may be assumed to be
independent and identically distributed with the law of
$X_0^N(0)$ being independent of $N$.
The intuition behind the comparison with the mean field model is that
if there is competition (modeled through the functions $h$ and $g$
in~\eqref{eq:def:N_islands_model})
among individuals and resources are everywhere
the same, then the best strategy for survival of the population
is to spread out in space as quickly as possible.

The results of~\cite{HW07} cover translation invariant initial measures
and local extinction. For general $h$ and $g$, not much is known about extinction
of the total mass process.
Let the solution $(X_t^N)_{t\geq0}$ of~\eqref{eq:def:N_islands_model} be
started in $X_0^N(i)=x\1_{i=0}$, $x\geq0$.
We prove in a forthcoming paper under suitable conditions on the parameters
that the total mass $\abs{X_t^N}:=\sum_{i=1}^N X_t^N(i)$
converges as $N\to\infty$.
In addition, we show in that paper that the limiting
process dominates the total mass process of the corresponding
system of interacting locally regulated diffusions started in finite
total mass.
Consequently, a global extinction result for the limiting process
would imply a global extinction result for systems of locally regulated diffusions.

In this paper we introduce and study a model which we call
\emph{Virgin Island Model} and which is the limiting
process of $(X_t^N)_{t\geq0}$ as $N\to\infty$.
Note that in the process $(X_t^N)_{t\geq0}$ an emigrant moves to a given
island with probability $\tfrac{1}{N}$.
This leads to the characteristic property of the Virgin Island Model
namely every emigrant moves to an unpopulated island. Our main
result is a necessary and
sufficient condition (see~\eqref{eq:bed_fuer_extinction} below)
for global extinction for the Virgin Island Model.
Moreover, this condition is fairly explicit in terms of the
parameters of the model.

Now we define the model.
On the $0$-th island evolves a diffusion $Y=(Y_t)_{t\geq0}$ with
state space $\R_{\geq0}$ given by the strong
solution of the stochastic
differential equation
\begin{equation}  \label{eq:def:Y}
  dY_t=-a(Y_t)\,dt+h(Y_t)\,dt+\sqrt{2g(Y_t)}dB_t,\qquad Y_0=y\geq0,
\end  {equation}
where $(B_t)_{t\geq0}$ is a standard Brownian motion.
This diffusion models the total mass of a population
and is the diffusion limit of
near-critical branching particle processes where both the offspring mean
and the offspring variance are regulated by the total population.
Later, we will specify conditions on $a, h$ and $g$ so that
$Y$ is well-defined. For now, we restrict our attention to the
prototype example of a Feller branching diffusion with logistic growth
in which $a(y)=\kappa y$, $h(y)=\gamma y(K-y)$ and
$g(y)=\beta y$ with $\kappa,\ga,K,\beta>0$.
Note that zero is a trap for $Y$, that is, $Y_t=0$ implies
$Y_{t+s}=0$ for all $s\geq0$.

Mass emigrates from the $0$-th island at rate $a(Y_t)\,dt$
and colonizes unpopulated islands.
A new population should evolve as the process $(Y_t)_{t\geq0}$.
Thus, we need the law of excursions of $Y$ from the trap zero.
For this, define the
set of excursions from zero by 
\begin{equation}  \label{eq:def:U}
  U:=\bigl\{\chi\in\mathbf{C}\rub{(-\infty,\infty),[0,\infty)}\colon
     T_0\in(0,\infty],
     \,\chi_t=0\fa t\in(-\infty,0]\cup[T_0,\infty)\bigr\}
\end  {equation}
where $T_y=T_y(\chi):=\inf\{t>0\colon\chi_t=y\}$ is the first hitting time
of $y\in[0,\infty)$.
The set $U$ is furnished with locally uniform convergence.
Throughout the paper, $\C\ru{S_1,S_2}$ and $\D\ru{S_1,S_2}$ denote
the set of continuous functions and the set of \cadlag\ functions,
respectively, between two intervals $S_1,S_2\subset\R$.
Furthermore, define
\begin{equation}  \label{eq:abuse_of_notation}
  \Dpos:=\curlb{\chi\in\D\rub{(-\infty,\infty),[0,\infty)}\colon\chi_t=0\ \fa t<0}.
\end  {equation}
The \emph{excursion measure} $\ohneQY$ is a
$\sigma$-finite measure on $U$.
It has been constructed by Pitman and Yor~\cite{PY82} as follows:
Under $\ohneQY$, the trajectories come from zero according to an
entrance law and then move according to the law of $Y$.
Further characterizations of $Q_Y$ are given in~\cite{PY82}, too.
For a discussion on the excursion theory of one-dimensional diffusions,
see~\cite{SVY07}.
We will give a definition of $\ohneQY$ later.

\begin{figure}[ht]
  \epsfig{file=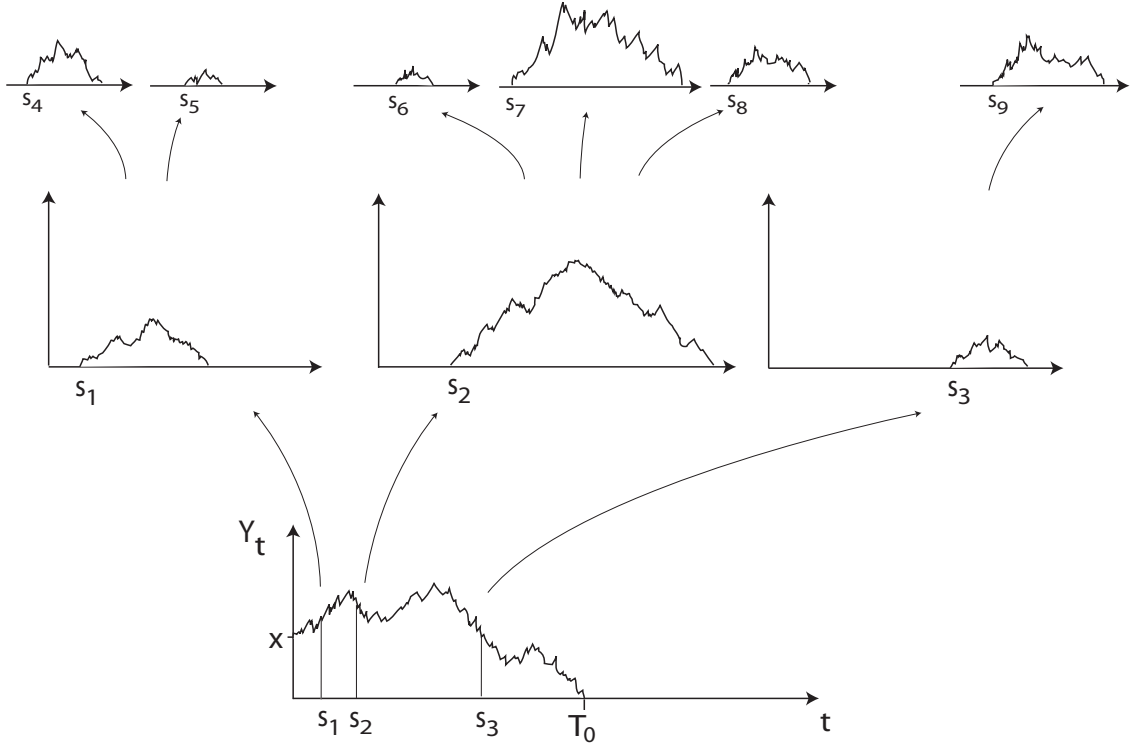, clip=,
      width=\linewidth}
  \caption{\footnotesize Subtree of the Virgin Island Model.
      Only offspring islands with a certain excursion height are drawn.
      Note that infinitely many islands are colonized e.g.\ between
      times $s_1$ and $s_2$.}
   \label{f:excursion_tree}
\end{figure}

Next we construct all islands which are colonized from the $0$-th
island and call these islands the first generation.
Then we construct
the second generation which is the collection of all islands which 
have been colonized from islands of the first generation, and so on.
Figure~\ref{f:excursion_tree}
illustrates the resulting tree of
excursions.
For the generation-wise construction, we use a method to index
islands which keeps track of which island has been colonized from
which island.
An island is identified with a triple which indicates  its mother island,
the time of its colonization and the population size on the island
as a function of time.
For $\chi\in\Dpos$, let
\begin{equation}
  \Isl_0^\chi:=\curlb{\rub{\emptyset,0,\chi}}
\end  {equation}
be a possible $0$-th island.
For each $n\geq1$ and $\chi\in\Dpos$, define
\begin{equation}
  \Isl_{n}^\chi:=\curlb{\rub{\iota_{n-1},s,\psi}\colon
    \iota_{n-1}\in\Isl_{n-1}^\chi,(s,\psi)\in [0,\infty)\times \Dpos}
\end  {equation}
which we will refer to as the set of all possible islands of the $n$-th
generation with fixed $0$-th island $(\emptyset,0,\chi)$.
This notation should be read as follows.
The island $\iota_n=(\iota_{n-1},s,\psi)\in\Isl_n^\chi$ has been colonized
from island $\iota_{n-1}\in \Isl_{n-1}^\chi$ at time $s$ and carries total
mass $\psi(t-s)$ at time $t\geq0$.
Notice that
there is no mass on an island before the time of its colonization.
The island space is defined by
\begin{equation}
  \Isl:=\{\emptyset\}\cup\bigcup_{\chi\in\Dpos}\Isl^\chi\ \text{ where }\ 
  \Isl^\chi:=\bigcup_{n\geq0}\Isl_{n}^\chi.
\end  {equation}
Denote by $\sigma_\iota:=s$ the colonization time of island $\iota$ if
$\iota=(\iota^{'},s,\psi)$ for some $\iota{'}\in\Isl$.
Furthermore, let
$\{\Pi^{\iota}\colon \iota\in\Isl\setminus\{\emptyset\}\}$
be a set of
Poisson point processes on $[0,\infty)\times\Dpos$
with intensity measure
\begin{equation}  \label{eq:intensity_measure_Pi}
  \E\eckb{ \Pi^{(\iota,s,\chi)}(dt\otimes d\psi)}
  =a\rub{\chi(t-s)}\,dt\otimes \ohneQY(d\psi)\quad \iota\in\Isl.
\end  {equation}
For later use, let $\Pi^\chi:=\Pi^{(\emptyset,0,\chi)}$.
We assume that the family $\{\Pi^\iota\colon\iota\in\Isl^\chi\}$ is independent for
every $\chi\in\Dpos$.

The Virgin Island Model is defined recursively generation by generation.
The $0$-th generation only consists of the $0$-th island
\begin{equation}  \label{eq:first_island}
  \VIM^{(0)}:=\curlb{\rub{\emptyset,0,Y}}.
\end  {equation}
The $(n+1)$-st generation, $n\geq0$,
is the (random) set of all islands
which have been colonized from islands of the $n$-th generation
\begin{equation} \label{eq:nth_island}
  \VIM^{(n+1)}:= \curlb{\rub{\iota_n,s,\psi}\in\Isl\colon
    \iota_n\in\VIM^{(n)}, \Pi^{\iota_n}\rub{\{(s,\psi)\}}>0}.
\end  {equation}
The set of all islands is defined by
\begin{equation} \label{eq:all_island}
  \VIM:=\bigcup_{n\geq0}\VIM^{(n)}.
\end  {equation}
The total mass process of the Virgin Island Model is defined by
\begin{equation}  \label{eq:def:V}
  V_t:=\sum_{\rub{\iota,s,\psi}\in\VIM}\psi(t-s),\quad t\geq0.
\end  {equation}
Our main interest concerns the behaviour of the
law $\Law{V_t}$ of $V_t$ as $t\to\infty$.

The following observation is crucial for understanding the behavior
of $(V_t)_{t\geq0}$ as $t\to\infty$.
There is an inherent
branching structure in the Virgin Island Model.
Consider as new ``time coordinate'' the number of island
generations.
One offspring island together with all its offspring islands
is again a Virgin Island Model 
but with the path $(Y_t)_{t\geq0}$ on the $0$-th island
replaced by an excursion path.
Because of this branching structure,
the Virgin Island Model is a multi-type
Crump-Mode-Jagers branching process
(see~\cite{Ja75} under ``general branching process'')
if we consider islands as individuals and $[0,\infty)\times\Dpos$
as type space.
We recall that a single-type Crump-Mode-Jagers process is a particle process
where every particle $i$ gives birth to particles at the time
points of a point process $\xi_i$ until its death at time $\ld_i$,
and $(\ld_i,\xi_i)_i$ are independent and identically distributed.
The literature on Crump-Mode-Jagers processes assumes
that the number of offspring per individual is finite in every finite
time interval.
In the Virgin Island Model, however,
every island has infinitely many offspring islands
in a finite time interval because $Q_Y$ is an infinite measure.

The most interesting question about the Virgin Island Model is
whether or not the process survives with positive probability as $t\to\infty$.
Generally speaking, branching particle processes survive
if and only if the expected number of offspring per particle is strictly
greater than one, e.g.\ the Crump-Mode-Jagers process survives
if and only if $\E\xi_i[0,\ld_i]>1$.
For the Virgin Island Model, the offspring of an island $(\iota,s,\chi)$
depends on the emigration intensities $a\rub{\chi(t-s)}dt$.
It is therefore not surprising that
the decisive parameter for survival
is the expected ``sum'' over those emigration intensities
\begin{equation}  \label{eq:man_hours}
  \int\int_0^\infty a\rub{\chi_t}\,dt\, Q_Y(d\chi).
\end  {equation}
We denote the expression in~\eqref{eq:man_hours} as
``expected total emigration intensity'' of the Virgin Island Model.
The observation that~\eqref{eq:man_hours} is the decisive parameter
plus an explicit formula for~\eqref{eq:man_hours} leads
to the following main result.
In Theorem~\ref{thm:extinction_VIM}, we will prove
that the Virgin Island Model survives
with strictly positive probability if and only if
\begin{equation} \label{eq:bed_fuer_extinction_intro}
  \int_0^\infty \frac{a(y)}{g(y)}
    \exp\ruB{\int_0^y\frac{-a(u)+h(u)}{g(u)}\,du} \,dy>1.
\end  {equation}
Note that the left-hand side of~\eqref{eq:bed_fuer_extinction_intro}
is equal to $\int_0^\infty a(y)m(dy)$ where $m(dy)$ is the
speed measure of the one-dimensional diffusion~\eqref{eq:def:Y}.
The method of proof for the extinction result is to study an integral
equation (see Lemma~\ref{l:integral_laplace_trafo})
which the Laplace transform of the total mass $V$ solves.
Furthermore, we will show in  Lemma~\ref{l:Q_explicit_man_hours}
that the expression
in~\eqref{eq:man_hours} is equal to the left-hand side
of~\eqref{eq:bed_fuer_extinction_intro}.

Condition~\eqref{eq:bed_fuer_extinction_intro} already appeared
in~\cite{HW07} as necessary and sufficient condition for existence
of a nontrivial invariant measure for the mean field model,
see Theorem 1 and Lemma 5.1 in~\cite{HW07}.
Thus, the total mass  process of the 
Virgin Island Model dies out if and only if
the mean field model~\eqref{eq:def:M} dies out.
The following duality indicates why the same condition appears
in two situations which seem to be fairly different at first view.
If $a(x)=\kappa x$, $h(x)=\gamma x(K-x)$ and $g(x)=\beta x$ with
$\kappa, \gamma, \beta>0$, that is, in the case of Feller branching
diffusions with logistic growth,
then model~\eqref{eq:def:N_islands_model}
is dual to itself, see Theorem 3 in~\cite{HW07}.
If $(X_t^N)_{t\geq0}$ indeed approximates the Virgin Island Model
as $N\to\infty$, then -- for this choice of parameters --
the total mass process $(V_t)_{t\geq0}$
is dual to the mean field model.
This duality would directly imply that -- in the case of Feller
branching diffusions with logistic growth --
global extinction of the Virgin Island
Model is equivalent to local extinction of the mean field model.

An interesting quantity of the Virgin Island process is
the area under the path of $V$.
In Theorem~\ref{thm:expected_man_hours}, we prove
that the expectation of
this quantity is finite exactly in the subcritical situation
in which case we give an expression
in terms of $a$, $h$ and $g$.
In addition, in the critical case and in the supercritical case,
we obtain the asymptotic behaviour of the expected 
area under the path of $V$
up to time $t$
\begin{equation}   \label{eq:expected_area_under_V_time_t}
  \int_0^t\E^x V_s\,ds
\end  {equation}
as $t\to\infty$ for all $x\geq0$.
More precisely, the order of~\eqref{eq:expected_area_under_V_time_t}
is $O(t)$ in the critical case. For the supercritical case, let
$\al>0$ be the Malthusian parameter defined by
\begin{equation}   \label{eq:Malthusian_intro}
  \int_0^\infty\ruB{ e^{-\al u}\int a\rub{\chi_u} Q_Y(d\chi)}\,du=1.
\end  {equation}
It turns out that the expression
in~\eqref{eq:expected_area_under_V_time_t}
grows exponentially with rate $\al$ as $t\to\infty$.

The result of Theorem~\ref{thm:expected_man_hours} in the supercritical
case suggests that the event that $(V_t)_{t\geq0}$ grows exponentially with rate $\al$
as $t\to\infty$ has positive probability.
However, this is not always true.
Theorem~\ref{thm:general:xlogx} proves that $e^{-\al t}V_t$
converges in distribution to a random variable $W\geq0$. Furthermore,
this variable is not identically zero if and only if
\begin{equation}  \label{eq:xlogx_intro}
  \int \rubb{\int_0^\infty a\rub{\chi_s} e^{-\al s}\,ds}
       \log^+\rubb{\int_0^\infty a\rub{\chi_s} e^{-\al s}\,ds}Q_Y(d\chi)<\infty
\end  {equation}
where $\log^+(x):= \max\{0,\log(x)\}$.
This $(x \log x)$-criterion is similar to the Kesten-Stigum
Theorem (see~\cite{KS66}) for
multidimensional Galton-Watson processes.
Our proof follows Doney~\cite{Don72} 
who establishes an $(x \log x)$-criterion for
Crump-Mode-Jagers processes.

Our construction introduces as new ``time coordinate'' the number of island
generations.
Readers being interested in a construction of the Virgin Island Model in the
original time coordinate -- for example in a relation between $V_t$ and
$(V_s)_{s<t}$ -- are referred to Dawson and Li (2003)~\cite{DaLi03}.
In that paper, a superprocess with dependent spatial motion and interactive
immigration is constructed as the pathwise unique solution of a stochastic
integral equation driven by a Poisson point process whose intensity measure
has as one component the excursion measure of the Feller branching diffusion.
In a special case
(see equation (1.6) in~\cite{DaLi03} with $x(s,a,t)=a$,
$q(Y_s,a)=\kappa Y_s(\R)$ and $m(da)=\1_{[0,1]}(a)\,da$),
this is just the Virgin Island Model with~\eqref{eq:def:Y} replaced by
a Feller branching diffusion, i.e.\ $a(y)=\kappa y$, $h(y)=0$, $g(y)=\beta y$.
It would be interesting to know whether existence and uniqueness of
such stochastic integral equations still hold if the excursion measure
of the Feller branching diffusion is replaced by $Q_Y$.

Models with competition have been studied by various authors.
Mueller and Tribe (1994)~\cite{MT94} and
Horridge and Tribe (2004)~\cite{HT04} investigate an one-dimensional SPDE analog
of interacting Feller branching diffusions with logistic growth which can
also be viewed as KPP equation with branching noise.
Bolker and Pacala (1997)~\cite{BP97} propose a branching random walk
in which the individual mortality rate is increased by a weighted sum
of the entire population.
Etheridge (2004)~\cite{Eth04} studies two diffusion limits hereof.
The ``stepping stone version of the Bolker-Pacala model'' is a system
of interacting Feller branching diffusions with non-local logistic growth.
The ``superprocess version of the Bolker-Pacala model'' is an analog of this
in continuous space.
Hutzenthaler and Wakolbinger~\cite{HW07}, motivated by~\cite{Eth04},
investigated interacting diffusions with local competition which is an analog of
the Virgin Island Model but with mass migrating on $\Z^d$ instead of
migration to unpopulated islands. 

%
%
\noindent
\section{Main results}%
\label{sec:main_results}
The following assumption guarantees existence and uniqueness of
a strong $[0,\infty)$-valued solution of
equation~\eqref{eq:def:Y}, see e.g.\ Theorem IV.3.1
in~\cite{IkWa}.
Assumption~\ref{a:A1} additionally requires that $a(\cdot)$ is
essentially linear.
\begin{assumption}\label{a:A1}
 The three functions $a\colon[0,\infty)\to[0,\infty)$,
 $h\colon[0,\infty)\to\R$ and $g\colon[0,\infty)\to[0,\infty)$
 are  locally Lipschitz continuous in $[0,\infty)$ and satisfy
 $a(0)=h(0)=g(0) = 0$. 
 The function $g$ is strictly positive on $(0,\infty)$.
 Furthermore, $h$ and $\sqrt{g}$ satisfy
 the linear growth condition
 \begin{equation}  \label{eq:growth_condition}
   \limsup_{x\to\infty}\frac{0\vee h(x)+\sqrt{g(x)}}{x}<\infty
 \end  {equation}
 where $x\vee y$ denotes the maximum of $x$ and $y$.
 In addition, $c_1\mal x\leq a(x)\leq c_2 \mal x$ holds for
 all $x\geq0$ and for some constants $c_1,c_2\in(0,\infty)$.
\end  {assumption}
The key ingredient in the construction of the Virgin Island Model
is the law of excursions of $(Y_t)_{t\geq0}$ from the
boundary zero.
Note that under Assumption~\ref{a:A1}, zero is an absorbing boundary
for~\eqref{eq:def:Y},
i.e.\ $Y_t=0$ implies $Y_{t+s}=0$ for all $s\geq0$.
As zero is not a regular point, it is not possible to apply
the well-established It\^o excursion theory.
Instead we follow Pitman and Yor~\cite{PY82} and obtain a
$\sigma$-finite measure $\barQY$
-- to be called \emph{excursion measure} --
on $U$
(defined in \eqref{eq:def:U}).
For this, we
additionally assume that $(Y_t)_{t\geq0}$ hits zero
in finite time with positive probability.
The following assumption formulates
a necessary and sufficient condition for this 
(see Lemma 15.6.2 in~\cite{KT2}).
To formulate the assumption, we define
\begin{equation} \label{eq:def:bar_s}
  \bs(z):=\exp\ruB{-\int_1^z\frac{-a(x)+h(x)}{g(x)}\,dx},\quad
  \bS(y):=\int_0^y\bar{s}{(z)}\,dz,\quad
  z,y>0.
\end  {equation}
Note that $\bS$ is a scale function,
that is,
\begin{equation}  \label{eq:def:scale_function}
  \P^y\rub{T_b(Y)<T_c(Y)}=\frac{\bS (y)-\bS (c)}{\bS (b)-\bS (c)}
\end  {equation}
holds for all $0\leq c< y< b<\infty$, see Section 15.6 in~\cite{KT2}.
\begin{assumption} \label{a:A2}
  The functions $a$, $g$ and $h$ satisfy
  \begin{equation} \label{eq:V_A2}
    \int_0^x\bS(y)\frac{1}{g(y)\bs(y)} \,dy<\infty
  \end  {equation}
  for some $x>0$.
\end  {assumption}
\noindent
Note that if Assumption~\ref{a:A2} is satisfied, then~\eqref{eq:V_A2} holds
for all $x>0$.

Pitman and Yor~\cite{PY82} construct the excursion measure $\Qb_Y$
in three different ways one being as follows.
The set of excursions reaching level $\dl>0$ has
$\Qb_Y$-measure $1/\Sb(\dl)$.
Conditioned on this event an excursion follows the diffusion
$(Y_t)_{t\geq0}$ conditioned to converge to infinity until this process
reaches level $\dl$. From this time on the excursion follows an independent
unconditioned process. We carry out this construction in detail
in Section~\ref{sec:excursions_from_a_trap_of_one_dimensional_diffusions}.
In addition
Pitman and Yor~\cite{PY82} describe the excursion measure
``in a preliminary way as''
\begin{equation}  \label{eq:limit_with_S}
  \limyO\frac{1}{\bS(y)}\Law[y]{Y}
\end  {equation}
where the limit indicates weak convergence of finite
measures on $\mathbf{C}\rub{[0,\infty),[0,\infty)}$ away from neighbourhoods of the zero-trajectory.
However, they do not give a proof.
Having $\Qb_Y$ identified as the limit in~\eqref{eq:limit_with_S} will
enable us to transfer explicit formulas for $\Law{Y}$ to explicit formulas
for $\Qb_Y$.
We establish the existence of the limit in~\eqref{eq:limit_with_S}
in Theorem~\ref{thm:existence_excursion_measure} below.
For this, let the topology on $\mathbf{C}\rub{[0,\infty),[0,\infty)}$ 
be given by locally uniform convergence.
Furthermore, recall $Y$ from~\eqref{eq:def:Y},
the definition of $U$ from~\eqref{eq:def:U}
and the definition of $\bS$ from~\eqref{eq:def:bar_s}.
%
%
\begin{theorem}  \label{thm:existence_excursion_measure}
  Assume~\ref{a:A1} and~\ref{a:A2}.
  Then there exists a $\sigma$-finite measure $\barQY$
  on $U$ such that
  \begin{equation}  \label{eq:def:wlim}
    \lim_{y\to0}\frac{1}{\bS(y)}\E^yF\ru{Y}=\int F(\chi)\barQY(d\chi)
  \end  {equation}
  for all bounded continuous
  $F\colon\mathbf{C}\rub{[0,\infty),[0,\infty)}\to\R$
  for which there exists an $\eps>0$ such that $F(\chi)=0$ whenever
  $\sup_{t\geq0}\chi_t<\eps$.
\end  {theorem}

For our proof of
the global
extinction result for the Virgin Island Model,
we need the scaling function $\bS$
in~\eqref{eq:def:wlim} to behave essentially linearly in a neighbourhood
of zero.
More precisely, we assume $\bS{'}(0)$ to exist
in $(0,\infty)$.
From definition~\eqref{eq:def:bar_s} of $\bS$ it is clear that
a sufficient condition for this is given by the following assumption.
\begin{assumption}  \label{a:S_bar}
  The integral
  $\int_\eps^1\tfrac{-a(y)+h(y)}{g(y)}\,dy$
  has a limit in $(-\infty,\infty)$ as $\eps\!\to0$.
\end  {assumption}
\noindent
It follows
from
dominated convergence
and from
the local Lipschitz continuity of $a$ and $h$
that
Assumption~\ref{a:S_bar} holds if $\int_0^1\tfrac{y}{g(y)}\,dy$
is finite.

In addition,
we assume that the expected total emigration intensity
of the Virgin Island  Model is finite.
Lemma~\ref{l:finite_man_hours} shows that,
under Assumptions~\ref{a:A1} and \ref{a:A2},
an equivalent condition for this is given in
Assumption~\ref{a:finite_man_hours}.
\begin{assumption}  \label{a:finite_man_hours}
  The functions $a$, $g$ and $h$ satisfy
  \begin{equation}
    \int_x^\infty \frac{a(y)}{g(y)\bs(y)}\,dy<\infty
  \end  {equation}
  for some and then for all $x>0$.
\end  {assumption}
\noindent%
We mention that
if Assumptions~\ref{a:A1}, \ref{a:A2}
and~\ref{a:finite_man_hours} hold, then the process $Y$ hits zero in finite
time almost surely
(see Lemma~\ref{l:finite_time_extinction} and
Lemma~\ref{l:finite_man_hours}).
Furthermore, we give a generic example for $a$, $h$ and $g$ namely
$a(y)=c_1 y$,
$h(y)=c_2 y^{\kappa_1}-c_3y^{\kappa_2}$, $g(y)=c_4y^{\kappa_3}$
with $c_1,c_2,c_3,c_4>0$.
The Assumptions~\ref{a:A1},
\ref{a:A2}, \ref{a:S_bar} and \ref{a:finite_man_hours}
are all satisfied if $\kappa_2>\kappa_1\geq1$ and if
$\kappa_3\in[1,2)$.
Assumption~\ref{a:A2} is not met by
$a(y)=\kappa y$, $\kappa>0$, $h(y)=y$ and $g(y)=y^2$
because then $\bs(y)=y^{\kappa-1}$,
$\bS(y)=y^{\kappa}/\kappa$ and
condition~\eqref{eq:V_A2} fails to hold.

Next we formulate the main result of this paper.
Theorem~\ref{thm:extinction_VIM} proves a nontrivial transition from
extinction to survival.
For the formulation of this result, we define
\begin{equation} \label{eq:def:s}
  s(z):=\exp\ruB{-\int_0^z\frac{-a(x)+h(x)}{g(x)}\,dx},\quad
  S(y):=\int_0^y s(z)\,dz,\quad
  z,y>0,
\end  {equation}
which is well-defined under Assumption~\ref{a:S_bar}.
Note that $\bS(y)=S(y)\bS^{'}(0)$.
Define the excursion measure
\begin{equation} \label{eq:excursion_measure}
  Q_Y:=\bS^{'}(0)\barQY
\end{equation}
and recall the total mass process $(V_t)_{t\geq0}$ from~\eqref{eq:def:V}.

%
%
\begin{theorem}\label{thm:extinction_VIM}
  Assume~\ref{a:A1}, \ref{a:A2}, \ref{a:S_bar} and \ref{a:finite_man_hours}.
  Then the total mass process $(V_t)_{t\geq0}$ started in $x>0$ dies out
  (i.e., converges in probability to zero as $t\to\infty$) if and only if
  \begin{equation} \label{eq:bed_fuer_extinction}
    \int_0^\infty
       \frac{a(y)}{g(y)s(y)} \,dy
    \leq 1.
  \end  {equation}
  If~\eqref{eq:bed_fuer_extinction} fails to hold, then
  $V_t$ converges in distribution as $t\to\infty$
  to a random variable $V_\infty$
  satisfying
  \begin{equation}
    \P^x(V_\infty=0)=1-\P^x( V_\infty=\infty)
     =\E^x\exp\ruB{-q\int_0^\infty a(Y_s)\,ds}
  \end  {equation}
  for all $x\geq0$ and some $q>0$.
\end  {theorem}
\begin{remark}
  The constant $q>0$ is the unique strictly positive fixed-point of a function 
  defined in Lemma~\ref{l:fixed_point}.
\end  {remark}

In the critical case, that is, equality in~\eqref{eq:bed_fuer_extinction},
$V_t$ converges to zero in distribution as $t\to\infty$.
However, it turns out that the expected area under the graph of $V$ is
infinite.
In addition, we obtain in Theorem~\ref{thm:expected_man_hours}
the asymptotic behaviour of the expected area under the graph  of $V$
up to time $t$ as $t\to\infty$.
For this, define
\begin{equation}
  w_a(x):=\int_0^\infty S(x\wedge z)\frac{a(z)}{g(z)s(z)}\,dz,\ x\geq0,
\end  {equation}
and similarly $w_{id}:=w_a$ with $a(z)=z$.
If Assumptions~\ref{a:A1}, \ref{a:A2}, \ref{a:S_bar}
and~\ref{a:finite_man_hours} hold, then $w_a(x)+w_{id}(x)$ is finite for
fixed $x<\infty$; see Lemma~\ref{l:finite_man_hours}.
Furthermore, under 
Assumptions~\ref{a:A1}, \ref{a:A2}, \ref{a:S_bar}
and~\ref{a:finite_man_hours},
\begin{equation}
  w_a^{'}(0)=\int_0^\infty \frac{a(z)}{g(z)s(z)}\,dz<\infty
\end  {equation}
by the dominated convergence theorem.

%
%
\begin{theorem} \label{thm:expected_man_hours}
  Assume~\ref{a:A1}, \ref{a:A2}, \ref{a:S_bar} and
  \ref{a:finite_man_hours}.
  If the left-hand side of~\eqref{eq:bed_fuer_extinction}
  is strictly smaller than one, then the expected
  area under the path of $V$ is equal to
  \begin{equation} \label{eq:man_hours_subcritical}
    \E^x\int_0^\infty V_s\,ds
    =w_{id}(x)+\frac{w_{id}^{'}(0)\;w_a(x)}
          {1-w_a^{'}(0)}\in(0,\infty)
  \end  {equation}
  for all $x\geq0$.
  Otherwise, the left-hand side of~\eqref{eq:man_hours_subcritical}
  is infinite.
  In the critical case, that is,
  equality in~\eqref{eq:bed_fuer_extinction},
  \begin{equation}  \begin{split} \label{eq:man_hours_critical}
    \frac{1}{t}\int_0^t \E^x V_s\,ds\to
     \frac{w_{id}^{'}(0)\;w_a(x)}
      {\int_0^\infty \frac{w_a(y)}{g(y)s(y)}dy}
    \in[0,\infty)\qqast
  \end  {split}     \end  {equation}
  where the right-hand side is interpreted as zero
  if the denominator is equal to infinity.
  In the supercritical case, i.e., if~\eqref{eq:bed_fuer_extinction} fails
  to be true, let $\al>0$ be such that
  \begin{equation}   \label{eq:Malthusian}
    \int_0^\infty e^{-\al s}\int a\rub{\chi_s} Q_Y(d\chi)\,ds=1.
  \end  {equation}
  Then the order of growth of the expected area under the
  path of $(V_s)_{s\geq0}$ up to time $t$ as $t\to\infty$ can be read off from
  \begin{equation} \label{eq:man_hours_supercritical}
    e^{-\al t}\int_0^t \E^x V_s\,ds\to
    \frac{\int_0^\infty      {e^{-\al s}\int\chi_s\ohneQY(d\chi)}\,ds\,
          \mal \int_0^\infty {e^{-\al s}\E^x a(Y_s)}\,ds}
         {\int_0^\infty \ruB{\al s e^{-\al s}\int a\ru{\chi_s}\ohneQY(d\chi)}ds}
    \in(0,\infty)
  \end  {equation}
  for all $x\geq0$.
\end  {theorem}

The following result is an analog of the Kesten-Stigum Theorem, see~\cite{KS66}.
In the supercritical case, $e^{-\al t}V_t$ converges to a random variable $W$
as $t\to\infty$.
In addition, $W$ is not identically zero if and only if the
$(x \log x)$-condition~\eqref{eq:xlogx_intro} holds.
We will prove a more general version hereof in Theorem~\ref{thm:general:xlogx}
below. Unfortunately, we do not know of an explicit formula in terms of
$a$, $h$ and $g$ for the left-hand side of~\eqref{eq:xlogx_intro}.
Aiming at a condition which is easy to verify, we assume instead
of~\eqref{eq:xlogx_intro} that the second moment
$\int(\int_0^\infty a(\chi_s)\,ds)^2 Q(d\chi)$ is finite.
In Assumption~\ref{a:finite_man_hours_squared}, we formulate a condition
which is slightly stronger than that,
see Lemma~\ref{l:Q_explicit_man_hours} below.
\begin{assumption}  \label{a:finite_man_hours_squared}
  The functions $a$, $g$ and $h$ satisfy
  \begin{equation}
    \int_x^\infty a(y)\frac{y+ w_a(y)}{g(y)\bs(y)}\,dy<\infty
  \end  {equation}
  for some and then for all $x>0$.
\end  {assumption}
%
%
\begin{theorem}  \label{thm:xlogx}
  Assume~\ref{a:A1}, \ref{a:A2},
  \ref{a:S_bar} and
  \ref{a:finite_man_hours_squared}.
  Suppose that~\eqref{eq:bed_fuer_extinction} fails to be true
  (supercritical case) and let $\al>0$ be the
  unique solution of~\eqref{eq:Malthusian}.
  Then
  \begin{equation}  \label{eq:weak_convergence_rescaled_VIM}
    \frac{V_t}{e^{\al t}}\wlim W\qqast
  \end  {equation}
  in the weak topology
  and $\P\{W>0\}=\P\{V_\infty>0\}$.
\end  {theorem}

%
%
\noindent
\section{Outline}%
\label{sec:outline}
Theorem~\ref{thm:existence_excursion_measure} will be established in
Section~\ref{sec:excursions_from_a_trap_of_one_dimensional_diffusions}.
Note that
Section~\ref{sec:excursions_from_a_trap_of_one_dimensional_diffusions}
does not depend on  the
sections~\ref{sec:random_characteristics}-\ref{sec:convergence_supercritical}.
We will prove the survival and extinction result of
Theorem~\ref{thm:extinction_VIM} in two steps.
In the first step, we obtain a criterion for survival and
extinction in terms of $Q_Y$. More precisely, we prove that
the process dies out if and only if the expression in~\eqref{eq:man_hours}
is smaller than or equal to one.
In this step, we do not exploit that $Q_Y$ is the
excursion measure of $Y$.
In fact, we will prove an analog of Theorem~\ref{thm:extinction_VIM}
in a more general setting where $Q_Y$ is replaced by some $\sigma$-finite
measure $Q$ and where the islands are counted with random characteristics.
See Section~\ref{sec:random_characteristics} below for the definitions.
The analog of Theorem~\ref{thm:extinction_VIM} is stated in
Theorem~\ref{thm:general:extinction},
see Section~\ref{sec:random_characteristics},
and will be proven in
Section~\ref{sec:extinction_VIM}.
The key equation for its proof is contained in
Lemma~\ref{l:tree_structure} which
formulates the branching structure in the Virgin Island Model.
In the second step, we calculate an expression for~\eqref{eq:man_hours}
in terms of $a,h$ and $g$.
This will be done in Lemma~\ref{l:Q_explicit_man_hours}.
Theorem~\ref{thm:extinction_VIM} is then a corollary
of Theorem~\ref{thm:general:extinction} and of
Lemma~\ref{l:Q_explicit_man_hours}, see Section~\ref{sec:proof_main_theorems}.
Similarly, a more general version of Theorem~\ref{thm:expected_man_hours}
is stated in Theorem~\ref{thm:general:expected_man_hours},
see Section~\ref{sec:random_characteristics} below.
The proofs of Theorem~\ref{thm:expected_man_hours} and of
Theorem~\ref{thm:general:expected_man_hours} are contained in
Section~\ref{sec:proof_main_theorems}
and
Section~\ref{sec:proof_of_theorem_expected_man_hours}, respectively.
As mentioned in Section~\ref{sec:introduction}, a rescaled version of
$(V_t)_{t\geq0}$ converges in the supercritical case. This convergence
is stated in a more general formulation in Theorem~\ref{thm:general:xlogx},
see Section~\ref{sec:random_characteristics} below.
The proofs of Theorem~\ref{thm:xlogx} and of
Theorem~\ref{thm:general:xlogx} are contained in
Section~\ref{sec:proof_main_theorems} and in
Section~\ref{sec:convergence_supercritical}, respectively.

%
%
\noindent
\section{Virgin Island Model counted with random characteristics}%
\label{sec:random_characteristics}
In the proof of the extinction result of Theorem~\ref{thm:extinction_VIM},
we exploit that one offspring island
together with all its offspring islands
is again a Virgin Island Model 
but with a typical excursion instead of $Y$ on the $0$-th island.
For the formulation of this branching property, we need a version
of the Virgin Island Model where the population on the $0$-th island is
governed by $Q_Y$. More generally, we replace the law $\Law{Y}$ of the first
island by some measure $\nu$ and we replace the excursion measure $Q_Y$
by some measure $Q$.
Given two $\sigma$-finite measures $\nu$ and $Q$ 
on the Borel-$\sigma$-algebra of
$\Dpos$, we define the
\emph{Virgin Island Model with initial island measure $\nu$
and excursion measure $Q$} as follows.
Define the random sets of islands $\VIM^{(n),\nu,Q}$, $n\geq0$,
and $\VIM^{\nu,Q}$ through the definitions
\eqref{eq:intensity_measure_Pi},
\eqref{eq:first_island},
\eqref{eq:nth_island} and
\eqref{eq:all_island}
with $\Law{Y}$ and $Q_Y$ replaced by $\nu$ and $Q$,
respectively.
A simple example for $\nu$ and $Q$ is
$\nu(d\chi)=Q(d\chi)=\E\dl_{t\mapsto\1_{t<L}}(d\chi)$
where $L\geq0$ is a random variable and $\dl_{\psi}$ is the Dirac measure
on the path $\psi$.
Then the Virgin Island Model coincides with a Crump-Mode-Jagers
process in which a particle has offspring according to a
rate $a(1)$ Poisson process until its death at time $L$.

Furthermore, our results do not only hold for the total mass process
\eqref{eq:def:V} but more generally when the islands are counted
with random characteristics. This concept is well-known for
Crump-Mode-Jagers processes, see Section 6.9 in~\cite{Ja75}.
Assume that $\phi_\iota=\rub{\phi_\iota(t)}_{t\in\R}$, $\iota\in\Isl$,
are separable
and nonnegative processes with the following properties.
It vanishes on the negative half-axis, i.e.\ $\phi_\iota(t)=0$ for $t<0$.
Informally speaking our main assumption on $\phi_\iota$ is that it
does not depend on the history.
Formally we assume that
\begin{equation}  \label{eq:main_assumption_phi}
  \ruB{\phi_{\rub{\iota,s,\chi}}(t)}_{t\in\R}
  \eqd
  \ruB{\phi_{\rub{\emptyset,0,\chi}}(t-s)}_{t\in\R}
  \ \ \fa\chi\in \Dpos,\iota\in\Isl, s\geq0.
\end  {equation}
Furthermore, we assume that
the family $\{\phi_{\iota},\Pi^{\iota}\colon\iota\in\Isl^\chi\}$ is independent
for each $\chi\in\Dpos$ and $(\om,t,\chi)\mapsto \phi_{(\emptyset,0,\chi)}(t)(\om)$
is measurable.
As a short notation, define
$\phi_\chi(t):=\phi(t,\chi):=\phi_{\ru{\emptyset,0,\chi}}(t)$ for
$\chi\in\Dpos$.
With this, we define
\begin{equation}  \label{eq:def:Vtphi}
  V_t^{\phi,\nu,Q}:=\sum_{\iota\in\VIM^{\nu,Q}}\phi_\iota(t-\sigma_\iota),
   \quad t\geq0,
\end  {equation}
and say that $\rub{V_t^{\phi,\nu,Q}}_{t\geq0}$ is a \emph{Virgin Island process
counted with random characteristics} $\phi$.
Instead of $V_t^{\phi,\dl_\chi,Q}$, we write $V_t^{\phi,\chi,Q}$
for a path $\chi\in\Dpos$ and note that
$(\om,t,\chi)\mapsto V_t^{\phi,\chi,Q}(\om)$ is measurable.
A prominent example for $\phi_\chi$ is the de\-ter\-mi\-nis\-tic random variable
$\phi_{\chi}(t)\equiv\chi(t)$.
In this case, $V_t^{\nu,Q}:=V_t^{\phi,\nu,Q}$ is
the total mass of all islands at time $t$.
Notice that $(V_t)_{t\geq0}$ defined in~\eqref{eq:def:V}
is a special case hereof, namely $V_t=V_t^{\Law{Y},\ohneQY}$.
Another example for $\phi_\chi$ is $\phi(t,\chi)=\chi(t)\1_{t\leq t_0}$.
Then $\VIM_t^{\phi,\chi,Q}$ is the total mass at time $t$ of all islands which
have been colonized in the last $t_0$ time units.
If $\phi(t,\chi)=\int_t^\infty \chi_s\,ds$,
then $V_t^{\phi,\chi,Q}=\int_t^\infty V_s^{\chi,Q}ds$.

As in Section~\ref{sec:main_results}, we need an assumption
which guarantees finiteness of $V_t^{\phi,\nu,Q}$.
\begin{assumption} \label{a:general:finite_moments}
  The function $a\colon[0,\infty)\to[0,\infty)$ is continuous and
  there exist $c_1,c_2\in(0,\infty)$ such that
  $c_1 x\leq a(x)\leq c_2 x$ for all $x\geq0$.
  Furthermore,
  \begin{equation}  \begin{split}
    \sup_{t\leq T}\int \ruB{a\rub{\chi_t}+\E\phi\rub{t,\chi}} \nu(d\chi)&+
    \sup_{t\leq T}\int \ruB{a\rub{\chi_t}+\E\phi\rub{t,\chi}} Q(d\chi)<\infty
  \end  {split}     \end  {equation}
  for every $T<\infty$
\end  {assumption}
\noindent
The analog of Assumption~\ref{a:finite_man_hours} in the
general setting is the following assumption.
\begin{assumption} \label{a:general:finite_man_hours}
  Both the expected emigration intensity of the $0$-th island and of subsequent
  islands are finite:
  \begin{equation} 
    \int \ruB{\int_0^\infty a\rub{\chi_u}\,du}\nu(d\chi)
   +\int \ruB{\int_0^\infty a\rub{\chi_u}\,du}Q(d\chi)<\infty.
  \end  {equation}
\end  {assumption}
\noindent
In Section~\ref{sec:main_results}, we assumed that $(Y_t)_{t\geq0}$
hits zero in finite time with positive probability.
See Assumption~\ref{a:A2} for an equivalent condition.
Together with~\ref{a:finite_man_hours}, this assumption implied
almost sure convergence of $(Y_t)_{t\geq0}$ to zero as $t\to\infty$.
In the general setting, we need a similar but somewhat weaker
assumption.
More precisely, we assume that $\phi(t)$ converges to zero
''in distribution`` both with respect to $\nu$ and with respect to $Q$.
\begin{assumption} \label{a:convergence_to_zero}
  The random processes $\curlb{\rub{\phi_\chi(t)}_{t\geq0}\colon\chi\in
  \Dpos}$ and the measures $Q$ and $\nu$
  satisfy
  \begin{equation}
    \int \ruB{1-\E e^{-\ld \phi(t,\chi)}}\rub{\nu+Q}(d\chi)\to 0\qqast
  \end  {equation}
  for all $\ld\geq 0$.
\end  {assumption}
%
%
Having introduced the necessary assumptions, we now formulate
the extinction and survival result of
Theorem~\ref{thm:extinction_VIM} in the general setting.
\begin{theorem}  \label{thm:general:extinction}
  Let $\nu$ be a probability measure on $\Dpos$
  and let $Q$ be a measure on $\Dpos$.
  Assume~\ref{a:general:finite_moments},
  \ref{a:general:finite_man_hours}
  and~\ref{a:convergence_to_zero}.
  Then the Virgin Island process $(V_t^{\phi,\nu,Q})_{t\geq0}$ counted with
  random characteristics $\phi$ with $0$-th island distribution
  $\nu$ and with excursion measure $Q$ dies out
  (i.e., converges to zero in probability)
  if and only if
  \begin{equation}   \label{eq:general:bed_fuer_extinction}
    \bar{a}:=\int \ruB{\int_0^\infty a\rub{\chi_u}\,du}Q(d\chi)\leq 1.
  \end  {equation}
  In case of survival,
  the process converges weakly as $t\to\infty$
  to a probability measure $\Law{V_\infty^{\phi,\nu,Q}}$ 
  with support in $\{0,\infty\}$ which puts mass
  \begin{equation}  \label{eq:Ws_V_gleich_unendlich}
     \int 1-\exp\ruB{-q\int_0^\infty a\rub{\chi_s}\,ds}\nu(d\chi)
  \end  {equation}
  on the point $\infty$
  where $q>0$ is the unique strictly positive fixed-point
  of
  \begin{equation} \label{eq:fixed_point}
    z\mapsto \int 1-\exp\ruB{-z\int_0^\infty a\rub{\chi_s}\,ds}Q(d\chi),
    \quad z\geq0.
  \end  {equation}
\end  {theorem}
\begin{remark}
  The assumption on $\nu$ to be a probability measure is convenient for
  the formulation in terms of convergence in probability. For a formulation
  in the case of a $\sigma$-finite measure $\nu$, see the proof of the theorem
  in Section~\ref{sec:extinction_VIM}.
\end  {remark}

Next we state Theorem~\ref{thm:expected_man_hours} in the general
setting.
For its formulation, define
\begin{equation}  \label{eq:def:f_nu}
  f^\nu(t):=\int \E\phi(t,\chi)\nu(d\chi),\quad t\geq0,
\end  {equation}
and similarly $f^Q$ with $\nu$ replaced by $Q$.

%
%
\begin{theorem} \label{thm:general:expected_man_hours}
  Assume~\ref{a:general:finite_moments} and
  \ref{a:general:finite_man_hours}.
  If the left-hand side of~\eqref{eq:general:bed_fuer_extinction}
  is strictly smaller than one and if both $f^\nu$ and $f^Q$
  are integrable, then
  \begin{equation} \label{eq:general:man_hours_subcritical}
    \int\E\eckB{ \int_0^\infty V_s^{\phi,\chi,Q}\,ds}\nu(d\chi)
    =\int_0^\infty \!\!\!f^\nu(s)ds
    +\frac{\int_0^\infty f^Q(s)\,ds\,\int
           \int_0^\infty a\rub{\chi_s}\,ds\nu(d\chi)}
          {1-\int \ruB{\int_0^\infty a\rub{\chi_s}\,ds}Q(d\chi)}
  \end  {equation}
  which is finite and strictly positive.
  Otherwise, the left-hand side of~\eqref{eq:general:man_hours_subcritical}
  is infinite.
  If the left-hand side of~\eqref{eq:general:bed_fuer_extinction}
  is equal to one and if both $f^\nu$ and $f^Q$
  are integrable,
  \begin{equation} \label{eq:general:man_hours_critical}
    \limt\frac{1}{t}\int\E\eckB{ \int_0^t V_s^{\phi,\chi,Q}\,ds}\nu(d\chi)
    =\frac{\int_0^\infty f^Q(s)\,ds\,\mal\int
           \int_0^\infty a\rub{\chi_s}\,ds\,\nu(d\chi)}
          {\int_0^\infty s\int  a\rub{\chi_s}Q(d\chi)\,ds}
    <\infty
  \end  {equation}
  where the right-hand side is interpreted as zero
  if the denominator is equal to infinity.
  In the supercritical case, i.e.,
  if~\eqref{eq:general:bed_fuer_extinction} fails
  to be true, let $\al>\!0$ be such that
  \begin{equation}   \label{eq:general:Malthusian}
    \int_0^\infty \ruB{e^{-\al s}\int a\rub{\chi_s} \,Q(d\chi)}\,ds=1.
  \end  {equation}
  Additionally assume that $f^Q$ is continuous a.e.\ with respect
  to the Lebesgue measure, 
  \begin{equation}  \label{eq:supercritical_renewal_theory}
    \sum_{k=0}^\infty \sup_{k\leq t<k+1}\abs{e^{-\al t}f^Q(t)}<\infty
  \end  {equation}
  and that $e^{-\al t}f^\nu(t)\to0$ as $t\to\infty$.
  Then the order of convergence of the expected total intensity  up to
  time $t$ can be read off from
  \begin{equation} \label{eq:general:man_hours_supercritical}
    \limt e^{-\al t}\int \E\eckB{\int_0^t V_s^{\phi,\chi,Q}\,ds}\nu(d\chi)
    =\frac{1}{\al}\,\limt e^{-\al t}\int \E \eckb{V_t^{\phi,\chi,Q}}\nu(d\chi)
  \end  {equation}
  and from
  \begin{equation} \label{eq:general:rescaled_VIM_supercritical}
    \limt e^{-\al t}\int \E \eckb{V_t^{\phi,\chi,Q}}\nu(d\chi)
    =\frac{\int_0^\infty e^{-\al s}f^Q(s)\,ds
           \,\mal\int_0^\infty e^{-\al s}\int a\rub{\chi_s}\nu(d\chi)\,ds}
          {\int_0^\infty s e^{-\al s}\int a\rub{\chi_s}Q(d\chi)\,ds}.
  \end  {equation}
\end  {theorem}

For the formulation of the analog of the
Kesten-Stigum Theorem, denote by
\begin{equation}  \label{eq:def:barm}
  \bm:=\frac{\int_0^\infty e^{-\al s}f^Q(s)\,ds}
                {\int_0^\infty se^{-\al s} \int a\rub{\chi_s}Q(d\chi)\,ds}
                \in(0,\infty)
\end  {equation}
the right-hand side of~\eqref{eq:general:rescaled_VIM_supercritical}
with $\nu$ replaced by $Q$.
Furthermore, define
\begin{equation}  \label{eq:def:A_al}
  A_\al(\chi):=\int_0^\infty a\rub{\chi_s} e^{-\al s}\,ds
\end  {equation}
for every path $\chi\in\Dpos$.
For our proof of Theorem~\ref{thm:general:xlogx}, we additionally
assume the following properties of $Q$.

\begin{assumption}  \label{a:for_SKDist}
  The measure $Q$ satisfies
  \begin{equation}  \label{eq:for_SKDist_1}
    \int\ruB{\int_0^T a(\chi_s)\,ds}^2Q(d\chi)<\infty
  \end  {equation}
  for every $T<\infty$ and
  \begin{equation}
       \sup_{t\geq0}\int\eckbb{\E\phi_\chi(t)
       \int_0^t a\rub{\chi_s} ds}Q(d\chi)<\infty,\quad
       \sup_{t\geq0}\int\E\rub{\phi_\chi^2(t)} Q(d\chi)<\infty.
  \end  {equation}
\end  {assumption}
%
%
\begin{theorem}  \label{thm:general:xlogx}
  Let $\nu$ be a probability measure on $\Dpos$
  and let $Q$ be a measure on $\Dpos$.
  Assume~\ref{a:general:finite_moments},
  \ref{a:general:finite_man_hours},
  \ref{a:convergence_to_zero}
  and~\ref{a:for_SKDist}.
  Suppose that $\bar{a}>1$ (supercritical case) and let $\al>0$ be the
  unique solution of~\eqref{eq:general:Malthusian}.
  Then
  \begin{equation}
    \frac{V_t^{\phi,\nu,Q}}{e^{\al t}\bm}\wlim W\qqast
  \end  {equation}
  in the weak topology
  where $W$ is a nonnegative random variable.
  The variable $W$ is not identically zero if and only if
  \vspace*{-1mm}
  \begin{equation}  \label{eq:xlogx}
    \int A_\al(\chi)\log^+\rub{A_\al(\chi)}Q(d\chi)<\infty
  \end  {equation}
  where $\log^+(x):= \max\{0,\log(x)\}$.
  If~\eqref{eq:xlogx} holds, then
  \begin{equation}  \begin{split} \label{eq:general:EW}
    \E W=\int\eckB{\int_0^\infty e^{-\al s}a\rub{\chi_s}\,ds}\nu(d\chi),
    \P\rub{W=0}=\int\eckB{e^{-q\int_0^\infty a\ru{\chi_s}\,ds}}\nu(d\chi)
  \end{split}     \end{equation}
  where $q>0$ is the unique strictly positive fixed-point
  of~\eqref{eq:fixed_point}.
\end  {theorem}
\begin{remark}
  Comparing~\eqref{eq:general:EW} with~\eqref{eq:Ws_V_gleich_unendlich},
  we see that $\P(W>0)=\P(V_\infty^{\phi,\nu,Q}>0)$.
  Consequently, the Virgin Island process $\rub{V_t^{\phi,\nu,Q}}_{t\geq0}$
  conditioned on not converging to zero grows exponentially fast with
  rate $\al$ as $t\to\infty$.
\end  {remark}

%
%
\noindent
\section{Branching structure}%
\label{sec:branching_structure}
  We mentioned in the introduction that there is an inherent
  branching structure in the Virgin Island Model.
  One offspring island together with all its offspring islands
  is again a Virgin Island Model 
  but with a typical excursion instead of $Y$ on the $0$-th island.
  In Lemma~\ref{l:tree_structure}, we formalize this idea.
  As a corollary thereof, we obtain an integral equation
  for the modified Laplace transform of the Virgin Island Model
  in Lemma~\ref{l:integral_laplace_trafo}
  which is the key equation for our proof of the extinction result
  of Theorem~\ref{thm:extinction_VIM}.
  Recall the notation of Section~\ref{sec:introduction} and of
  Section~\ref{sec:random_characteristics}.

  \newcommand{\Vau}{\operatorname{V}}
  \newcommand{\VIMop}{\mathop{\VIM}}
  \newcommand{\Vspsi}[1]{\sideset{^{(s,\psi)}}{_{#1}^{\phi,\chi,Q}}\Vau}
  \newcommand{\Vspsin}[2]{\sideset{^{(s,\psi)}}{_{#1}^{(#2),\phi,\chi,Q}}\Vau}
  \newcommand{\VIMspsi}[1]{\sideset{^{(s,\psi)}}{^{(#1),\chi}}\VIMop}
  \newcommand{\VIMchi}{\sideset{^{(s,\psi)}}{^{\chi}}\VIMop}
 \begin{lemma}  \label{l:tree_structure}
    Let $\chi\in \Dpos$.
    There exists an independent family 
    \begin{equation}  \label{eq:independent_family}
      \curlB{\ruB{\Vspsi{t}}_{t\geq0}\colon (s,\psi)\in[0,\infty)\times\Dpos}
    \end  {equation}
    of random variables which is independent of $\phi_\chi$ and of
    $\Pi^{\chi}$
    such that
    \begin{equation}  \label{eq:tree_structure}
      V_t^{\phi,\chi,Q}
      =\phi_\chi(t)+
          \!\!\sum_{(s,\psi)\in\Pi^{\chi}}\!\!\!\!\!
              \Vspsi{t}\quad\fa t\geq0
    \end  {equation}
    and such that
    \begin{equation}  \label{eq:Vspsi_eqd_V}
      \ruB{\Vspsi{t}}_{t\geq0}\eqd \ruB{V_{t-s}^{\phi,\psi,Q}}_{t\geq0}
    \end  {equation}
    for all $(s,\psi)\in[0,\infty)\times\Dpos$.
  \end  {lemma}
  \begin{proof}
    Write $\VIM^\chi:=\VIM^{\chi,Q}$
    and $\VIM^{(n),\chi}:=\VIM^{(n),\chi,Q}$. Define
    \begin{equation}  \label{eq:VIMspsi}
      \VIMspsi{1}:=\curlB{\rub{(\emptyset,0,\chi),s,\psi}}\subset\Isl_1^\chi
      \text{ and }
      \VIMchi:=\bigcup\nolimits_{n\geq1}\!\!\!\!\!\!\VIMspsi{n}
    \end  {equation}
    for $(s,\psi)\in[0,\infty)\times\Dpos$ where
    \begin{equation}  \label{eq:VIMspsin}
      \VIMspsi{n+1}:=\curlB{\rub{\iota_n,r,\zeta}\in\Isl_{n+1}^\chi\colon
        \iota_n\in\VIMspsi{n},\Pi^{\iota_n}(r,\zeta)>0}
    \end  {equation}
    for $n\geq1$.
    Comparing~\eqref{eq:VIMspsi} and~\eqref{eq:VIMspsin}
    with~\eqref{eq:nth_island}, we see that
    \begin{equation}  \label{eq:mengenaufteilung}
      \VIM^{(0),\chi}=\curlb{(\emptyset,0,\chi)}\ \text{ and }\ 
               \VIM^{(n),\chi}=\bigcup_{(s,\psi)\in\Pi^\chi}\!\!\!\!\VIMspsi{n}
               \fa n\geq1.
    \end  {equation}
    Define $V_t^{(0),\phi,\chi,Q}=\phi_\chi(t)$ for $t\geq0$ and for $n\geq1$
    \begin{equation}  \begin{split}  \label{eq:def:Vn}
      V_t^{(n),\phi,\chi,Q}:=
        \sum_{(s,\psi)\in\Pi^\chi}\sum_{\iota\in\VIMspsi{n}} \phi_\iota(t-\sigma_\iota)
        =:\sum_{(s,\psi)\in\Pi^\chi}\Vspsin{t}{n}.
    \end{split}     \end{equation}
    Summing over $n\geq0$ we obtain for $t\geq0$
    \begin{equation}
      V_t^{\phi,\chi,Q}=\phi_\chi(t)
        +\sum_{(s,\psi)\in\Pi^\chi}\sum_{n\geq1}\Vspsin{t}{n}
      =:\phi_\chi(t)
        +\sum_{(s,\psi)\in\Pi^\chi}\Vspsi{t}.
    \end  {equation}
    This is equality~\eqref{eq:tree_structure}.
    Independence of the family~\eqref{eq:independent_family}
    follows from independence of $(\Pi^\iota)_{\iota\in\Isl^\chi}$
    and from independence of
    $(\phi_\iota)_{\iota\in\Isl^\chi}$.
    It remains to prove~\eqref{eq:Vspsi_eqd_V}.
    Because of assumption~\eqref{eq:main_assumption_phi} the random
    characteristics $\phi_\iota$ only depends on the last part of $\iota$.
    Therefore
    \begin{equation}  \begin{split}  \label{eq:Vspsin_eqd_V}
      \Vspsin{\cdot}{n}
      &=\sum_{\iota\in\VIMspsi{n}}\phi_\iota\rub{\cdot-\sigma_\iota}\\
      &\eqd \sum_{\tilde{\iota}\in\VIM^{(n-1),\psi,Q}}\phi_{\tilde{\iota}}
         (\cdot-\ru{\sigma_{\tilde{\iota}}+s})
      =V_{\cdot-s}^{(n-1),\psi,Q}.
    \end  {split}     \end  {equation}
    Summing over $n\geq1$ results in~\eqref{eq:Vspsi_eqd_V} and finishes the proof.
  \end  {proof}

  In order to increase readability, we introduce
  the following suggestive symbolic abbreviation
  \begin{equation}   \label{eq:def:I}
    \I\eckB{f\rub{V_t^{\phi,\nu,Q}}}
    :=\int\E f\rub{ V_t^{\phi,\chi,Q}}\nu(d\chi)
    \quad t\geq0, f\in\C\rub{[0,\infty),[0,\infty)}.
  \end  {equation}
  One might want to read this as ``expectation'' with respect to a
  non-probability measure. However,   \eqref{eq:def:I} is not intended
  to define an operator.

  The following lemma proves that
  the Virgin Island Model counted with random characteristics
  as defined in~\eqref{eq:def:Vtphi} is finite.
%
%
%
  \begin{lemma}  \label{l:general:finite_moments}
    Assume~\ref{a:general:finite_moments}.
    Then,
    for every $T<\infty$,
    \begin{equation}    \label{eq:first_moment_V}
      \sup_{t\leq T}\I\eckB{V_t^{\phi,\nu,Q}}
        <\infty.
    \end  {equation}
    Furthermore, if
    \begin{equation}  \label{eq:ass_second_moment}
      \sup_{t\leq T}\int \E\rub{\phi_\chi^2(t)}
         +\ruB{\int_0^T a(\chi_s)\,ds}^2\,Q(d\chi)<\infty,
    \end{equation}
    then there exists a constant $c_T<\infty$ such that
    \begin{equation}  \begin{split}  \label{eq:second_moment_V}
      \lefteqn{\sup_{t\leq T} \I\eckbb{\ruB{V_t^{\phi,\nu,Q}}^2}}\\
      &\leq c_T\ruB{1+\sup_{t\leq T}\int \E\rub{\phi_\chi^2(t)}\, (\nu+Q)(d\chi)
      +\int\ruB{\int_0^T a(\chi_s)ds}^2\nu(d\chi)}
    \end  {split}     \end  {equation}
    for all $\nu$ and the right-hand side of~\eqref{eq:second_moment_V} is finite
    in the special case $\nu=Q$.
  \end  {lemma}
  \begin{proof}
    We exploit the branching property formalized in
    Lemma~\ref{l:tree_structure} and apply Gronwall's inequality.
    Recall $\Vau^{(n),\chi,Q}$ from the proof of
    Lemma~\ref{l:tree_structure}.
    The two equalities~\eqref{eq:def:Vn} and~\eqref{eq:Vspsin_eqd_V}
    imply
    \begin{equation}  \label{eq:recursive_equation0}
      \I\eckB{V_t^{(0),\phi,\nu,Q}}
      =\int \E\phi_\chi(t)\,\nu(d\chi)\leq
      \sup_{s\leq T} \int \E\phi_\chi(s)\,\nu(d\chi)
    \end{equation}
    for $t\leq T$ and for $n\geq1$
    \begin{equation}  \begin{split} \label{eq:recursive_equation}
      \I&\eckB{V_t^{(n),\phi,\nu,Q}}
      =\int \E\eckB{\sum_{(s,\psi)\in\Pi^{\chi}}
        \E\eckb{V_{t-s}^{(n-1),\phi,\psi,Q}}}\nu(d\chi)\\
      &=\int \rubb{\int_0^t\int
        \E\eckb{V_{t-s}^{(n-1),\phi,\psi,Q}}Q(d\psi) a\ru{\chi_s}ds}\,\nu(d\chi)\\
      &\leq \sup_{u\leq T}\int a\ru{\chi_u}\nu(d\chi)\int_0^t
        \I\eckB{V_{s}^{(n-1),\phi,Q,Q}} \,ds.
    \end  {split}     \end  {equation}
    Using Assumption~\ref{a:general:finite_moments}
    induction on $n\geq0$ shows that all expressions
    in~\eqref{eq:recursive_equation0} and in~\eqref{eq:recursive_equation}
    are finite in the case $\nu=Q$.
    Summing~\eqref{eq:recursive_equation} over $n\leq n_0$ we obtain
    \begin{equation}  \begin{split}  \label{eq:recursive_equation2}
      \sum_{n=0}^{n_0}\I\eckB{V_t^{(n),\phi,\nu,Q}}
      \leq \int\E\phi_\chi(u)\nu(d\chi)
       +\int_0^t\sum_{n=0}^{n_0}\I\eckB{V_s^{(n),\phi,Q,Q}}\int a\ru{\chi_{t-s}}\nu(d\chi)\,ds
    \end  {split}     \end  {equation}
    for $t\leq T$. In the special case $\nu=Q$  Gronwall's inequality implies
    \begin{equation}       \label{eq:Gronwall_V}
      \sum_{n=0}^{n_0}\I\eckB{V_t^{(n),\phi,Q,Q}}\leq \sup_{u\leq T}\int\E\phi_\chi(u)Q(d\chi)
      \mal\exp\ruB{ t \sup_{u\leq T}\int a\ru{\chi_u}Q(d\chi)}.
    \end  {equation}
    Summing~\eqref{eq:recursive_equation} over $n\leq n_0$,
    inserting~\eqref{eq:Gronwall_V}
    into~\eqref{eq:recursive_equation}
    and letting $n_0\to\infty$
    we see that~\eqref{eq:first_moment_V}
    follows from Assumption~\ref{a:general:finite_moments}.

    For the proof of~\eqref{eq:second_moment_V},
    note that~\eqref{eq:recursive_equation2} with $\nu=\dl_\chi$
    and~\eqref{eq:first_moment_V}
    imply
    \begin{equation}  \label{eq:sufficient_to_consider_variances}
      \int\ruB{\E V_t^{\phi,\chi,Q}}^2 Q(d\chi)
      \leq\int 2\ruB{\E\phi_\chi(t)}^2+\ct_T\Bigl(\int_0^T a(\chi_s)ds\Bigr)^2 Q(d\chi)<\infty
    \end  {equation}
    for some $\ct_T<\infty$.
    In addition the two equalities~\eqref{eq:def:Vn} and~\eqref{eq:Vspsin_eqd_V}
    together with independence imply
    \begin{equation}
      \int \Var\rub{ V_t^{(0),\phi,\chi,Q}}\nu(d\chi)
         =\int \Var\rub{\phi_\chi(t)}\nu(d\chi)
    \end{equation}
    for $t\geq0$ and for $n\geq1$
    \begin{equation}  \begin{split}   \label{eq:goofy_var}
      \lefteqn{\int \Var\rub{ V_t^{(n),\phi,\chi,Q}}\nu(d\chi)}\\
      &=\int\E\rubb{\sum_{(s,\psi)\in\Pi^\chi}
         \Var \rub{V_{t-s}^{(n-1),\phi,\psi,Q}}}\,\nu(d\chi)\\
      &=\int\int_0^t \ruB{a(\chi_s)
         \int \Var \rub{V_{t-s}^{(n-1),\phi,\psi,Q}} Q(d\psi)}
         ds\,\nu(d\chi)\\
      &\leq\int_0^t
         \int \Var \rub{V_{s}^{(n-1),\phi,\psi,Q}} Q(d\psi) ds
         \mal\sup_{u\leq T}\int a(\chi_u)\,\nu(d\chi).
    \end  {split}     \end  {equation}
    In the special case $\nu=Q$ induction on $n\geq0$ together
    with~\eqref{eq:ass_second_moment} shows that all involved
    expressions are finite.
    A similar estimate as in~\eqref{eq:goofy_var} leads to
    \begin{equation*}  \begin{split}   \label{eq:goofy_var_sum}
      \lefteqn{\int \E\eckB{\Bigl(\sum_{n=0}^{n_0} V_t^{(n),\phi,\chi,Q}\Bigr)^2}\nu(d\chi)
        -\int\Bigl(\E\sum_{n=0}^{n_0} V_t^{(n),\phi,\chi,Q}\Bigr)^2\nu(d\chi)}\\
      &=\int \Var\rub{\phi_\chi(t)}+\E\rubb{\sum_{(s,\psi)\in\Pi^\chi}
         \Var \rub{\sum_{n=1}^{n_0}V_{t-s}^{(n-1),\phi,\psi,Q}}}\,\nu(d\chi)\\
      &=\int \Var\rub{\phi_\chi(t)}+\int_0^t \ruB{a(\chi_s)
         \int \Var \rub{\sum_{n=0}^{n_0-1}V_{t-s}^{(n),\phi,\psi,Q}} Q(d\psi)}
         ds\,\nu(d\chi)\\
      &\leq \int\E\ruB{\phi_\chi^2(t)}\nu(d\chi)+
         \int_0^t
         \int \E\eckB{\Bigl(\sum_{n=0}^{n_0}V_{s}^{(n),\phi,\psi,Q}\Bigr)^2} Q(d\psi) ds
         \mal\sup_{u\leq T}\int a(\chi_u)\,\nu(d\chi).
    \end  {split}     \end  {equation*}
    In the special case $\nu=Q$ 
    Gronwall's inequality together with~\eqref{eq:sufficient_to_consider_variances}
    leads to  
    \begin{equation}  \begin{split}  \label{eq:sup_var_V}
      \lefteqn{\int \E\eckB{\Bigl(\sum_{n=0}^{n_0}V_t^{(n),\phi,\chi,Q}\Bigr)^2}Q(d\chi)}\\
      &\leq\rubb{\int 3\E\ruB{\phi_\chi^2 (t)}
         +\ct_T\Bigl(\int_0^T a(\chi_s)\,ds\Bigr)^2Q(d\chi)}
         \exp\ruB{\sup_{u\leq T} \int a(\chi_u)\,Q(d\chi)T}
    \end{split}     \end{equation}
    which is finite by Assumption~\ref{a:general:finite_moments}
    and assumption~\eqref{eq:ass_second_moment}.
    Inserting~\eqref{eq:sup_var_V} into~\eqref{eq:goofy_var} and 
    letting $n_0\to\infty$ finishes the proof.
  \end  {proof}

  In the following lemma, we establish an integral equation
  for the modified Laplace transform of the Virgin Island Model.
  Recall the definition of $V_t^{\phi,\nu,Q}$ from~\eqref{eq:def:Vtphi}.
  \begin{lemma}  \label{l:integral_laplace_trafo}
    Assume~\ref{a:general:finite_moments}.
    The modified Laplace transform 
    $\I\eckb{1-e^{-\ld V_t^{\phi,\nu,Q}}}$
    of the Virgin Island Model counted with random characteristics $\phi$
    satisfies
    \begin{equation}  \begin{split}  \label{eq:integral_laplace_trafo}
      \lefteqn{\I\eckB{1-e^{-\ld V_t^{\phi,\nu,Q}}}}\\
      &=\int\E\eckB{1-
        \exp\ruB{-\ld\phi_{\chi}(t)-\int_0^\infty
          \I\eckb{1-e^{-\ld V^{\phi,Q,Q}_{t-s}}}
          a\ru{\chi_s}\,ds}}\nu(d\chi)
    \end  {split}     \end  {equation}
    for all $\ld,t\geq0$.
  \end  {lemma}
  \begin{proof}
    Fix $\ld,t\geq0$.
    Applying Lemma~\ref{l:tree_structure},
    \begin{equation*}  \begin{split}
      \lefteqn{\I\eckB{1-e^{-\ld V_t^{\phi,\nu,Q}}}}\\
      &=\int\eckB{1-\E \rub{e^{-\ld\phi_\chi(t)}}\mal
        \E\ruB{\prod_{(s,\psi)\in\Pi^\chi} \E e^{-\ld V^{\phi,\psi,Q}_{t-s}}}}
        \nu(d\chi)\\
      &=\int\eckB{1-\E \rub{e^{-\ld\phi_\chi(t)}}\mal
        \exp\ruB{-\int_0^\infty\int 1-
          \E e^{-\ld V^{\phi,\psi,Q}_{t-s}}Q(d\psi)a\ru{\chi_s}\,ds}}
          \nu(d\chi)\\
      &=\int\E\eckB{1-
        \exp\ruB{-\ld\phi_{\chi}(t)-\int_0^\infty
          \I\eckb{1-e^{-\ld V^{\phi,Q,Q}_{t-s}}}
          a\ru{\chi_s}\,ds}}\nu(d\chi).
    \end  {split}     \end  {equation*}
    This proves the assertion.
  \end  {proof}

%
\noindent
\section{Proof of Theorem~\ref{thm:general:expected_man_hours}}%
\label{sec:proof_of_theorem_expected_man_hours}

  Recall the definition of $(V_t^{\phi,\nu,Q})_{t\geq0}$
  from~\eqref{eq:def:Vtphi},
  $f^\nu$ from~\eqref{eq:def:f_nu}
  and the notation $\I$ from~\eqref{eq:def:I}.
  We begin with the supercritical case and let $\al>0$ be the
  Malthusian parameter which is the unique solution
  of~\eqref{eq:general:Malthusian}.
  Define
  \begin{equation}  \label{eq:def:mu_nu}
    m^\nu(t):=\I\eckB{V_t^{\phi,\nu,Q}}\qquad
    \mu^\nu(ds):=\int a\rub{\chi_s} \,\nu(d\chi)\,ds
  \end  {equation}
  for $t\geq0$. In this notation, equation~\eqref{eq:recursive_equation}
  with $\nu$ replaced by $Q$ reads as
  \begin{equation}  \label{eq:recursive_m_Q}
    e^{-\al t}m^Q(t) =e^{-\al t}f^Q(t)+\int_0^t 
      e^{-\al (t-s)} m^Q(t-s) e^{-\al s}\mu^Q(ds).
  \end  {equation}
  This is a renewal equation for $e^{-\al t}m^Q(t)$.
  By definition of $\al$, $e^{-\al s}\mu^Q(ds)$ is a probability measure.
  From Lemma~\ref{l:general:finite_moments}
  we know that $m^Q$ is bounded on finite intervals.
  By assumption,
  $f^Q$ is continuous Lebesgue-a.e.\ and
  satisfies~\eqref{eq:supercritical_renewal_theory}.
  Hence, we may apply standard renewal theory
  (e.g.\ Theorem 5.2.6 of~\cite{Ja75}) and obtain
  \begin{equation}  \label{eq:supercrit_Q}
    \lim_{t\to\infty}e^{-\al t} m^Q(t)
    =\frac{\int_0^\infty e^{-\al s}f^Q(s)\,ds}
          {\int_0^\infty se^{-\al s}\mu^Q(ds)}<\infty.
  \end  {equation}
  Multiply equation~\eqref{eq:recursive_equation} by $e^{-\al t}$,
  recall $e^{-\al t}f^\nu(t)\to0$ as $t\to\infty$
  and apply the dominated convergence theorem together
  with~\ref{a:general:finite_man_hours} to obtain
  \begin{equation}  \begin{split}  \label{eq:toarriveat}
    \lim_{t\to\infty}e^{-\al t}m^\nu(t)
    =\int_0^\infty e^{-\al s} \lim_{t\to\infty}e^{-\al(t-s)}m^Q(t-s)
      \mu^\nu(ds).
  \end  {split}     \end  {equation}
  Insert~\eqref{eq:supercrit_Q} to obtain
  equation~\eqref{eq:general:rescaled_VIM_supercritical}.
  An immediate consequence of
  the existence of the limit on the left-hand side
  of~\eqref{eq:toarriveat} is
  equation~\eqref{eq:general:man_hours_supercritical}
  \begin{equation}
    e^{-\al t}\int_0^t m^\nu(s)\,ds
    =\int_0^\infty e^{-\al s}\mal e^{-\al(t-s)}m^\nu(t-s)\,ds
    \lrat\frac{1}{\al}\mal\limt e^{-\al t}m^\nu(t)
  \end  {equation}
  where we used the dominated convergence theorem.

  Next we consider the subcritical and the critical case.
  Define
  \begin{equation} 
    \bar{x}^\nu(t):=\int_0^t\I\eckB{V_s^{\phi,\nu,Q}}\,ds,\quad t\geq0.
  \end  {equation}
  In this notation, equation~\eqref{eq:recursive_equation}
  integrated over $[0,t]$ reads as
  \begin{equation}  \label{eq:renewal_equation}
    \bar{x}^\nu(t)
    =\int_0^t f^\nu(s)\,ds+\int_0^t \bar{x}^Q(t-u)\mu^\nu(du),\quad t\geq0.
  \end  {equation}
  In the subcritical case, $f^Q$ and $f^\nu$ are integrable.
  Theorem 5.2.9 in~\cite{Ja75} applied to~\eqref{eq:renewal_equation}
  with $\nu$ replaced by $Q$ implies
  \begin{equation}  \label{eq:limxQ}
    \lim_{t\to\infty}\bar{x}^Q(t)
     =\frac{\int_0^\infty f^Q(s)ds}{1-\mu^Q\rub{[0,\infty)}}.
  \end  {equation}
  Letting $t\to\infty$ in~\eqref{eq:renewal_equation}, dominated convergence
  and $\mu^\nu\rub{[0,\infty)}<\infty$ imply
  \begin{equation}  \begin{split}
    \lim_{t\to\infty}\bar{x}^\nu(t)=\int_0^\infty f^\nu(s)ds+
     \int_0^\infty \lim_{t\to\infty}\bar{x}^Q(t-u)\mu^\nu(du).
  \end  {split}     \end  {equation}
  Inserting~\eqref{eq:limxQ} results
  in~\eqref{eq:general:man_hours_subcritical}.
  In the critical case,
  similar arguments lead to
  \begin{equation}  \begin{split}
    \lefteqn{\lim_{t\to\infty}\frac{1}{t}\bar{x}^\nu(t)}\\
    &=\lim_{t\to\infty}\frac{1}{t}\int_0^t f^\nu(s)ds+
     \int_0^\infty \lim_{t\to\infty}\frac{t-u}{t}
       \lim_{t\to\infty}\frac{1}{t-u}\bar{x}^Q(t-u)\mu^\nu(du)\\
    &=\frac{\int_0^\infty f^Q(s)\,ds}{\int_0^\infty u\mu^Q(du)}
    \mu^\nu\rub{[0,\infty)}.
  \end  {split}     \end  {equation}
  The last equality follows from~\eqref{eq:renewal_equation}
  with $\nu$ replaced by $Q$
  and Corollary 5.2.14 of~\cite{Ja75}
  with $c:=\int_0^\infty f^Q(s)\,ds$, $n:=0$ and
  $\theta:=\int_0^\infty u\mu^Q(du)$.
  Note that the assumption $\theta<\infty$
  of this corollary is not necessary for this conclusion.
\qed

%
\section{Extinction and survival in the Virgin Island Model.
  Proof of Theorem~\ref{thm:general:extinction}}%
\label{sec:extinction_VIM}
  Recall the definition of $(V_t^{\phi,\nu,Q})_{t\geq0}$
  from~\eqref{eq:def:Vtphi}
  and the notation $\I$ from~\eqref{eq:def:I}.
As we pointed out in Section~\ref{sec:main_results},
the expected total emigration intensity of the Virgin Island Model
plays an important role.
The following lemma provides us with some properties of the
modified Laplace transform of the total emigration intensity.
These properties are crucial for
our proof of Theorem~\ref{thm:general:extinction}.
%
%
%
\begin{lemma} \label{l:fixed_point}
  Assume \ref{a:general:finite_man_hours}.
  Then the function
  \begin{equation}
    k(z):=\int 1-\exp\ruB{-z\int_0^\infty a\rub{\chi_s}\,ds}Q(d\chi),
    \quad z\geq0,
  \end  {equation}
  is concave with at most two fixed-points.
  Zero is the only fixed-point if and only if
  \begin{equation} \label{eq:k_prime}
    k^{'}(0)=\int\int_0^\infty a\rub{\chi_s}\,ds\, Q(d\chi)\leq1.
  \end  {equation}
  Denote by $q$ the maximal
  fixed-point. Then we have for all $z\geq0$:
  \begin{eqnarray}
    z\leq k(z) &\implies& z\leq q \label{aa:1}\\ 
    z\geq k(z)\wedge z>0 &\implies & z\geq q.  \label{aa:2}
  \end  {eqnarray}
\end  {lemma}
\begin{proof}
  If $\int_0^\infty a\rub{\chi_s}\,ds=0$ for $Q$-a.a.\ $\chi$, then
  $k\equiv0$ and zero is the only fixed-point.
  For the rest of the proof, we assume w.l.o.g.\ that
  $\int\rub{\int_0^\infty a\ru{\chi_s}\,ds}Q(d\chi)>0$.

  The function $k$ has finite values because of $1-e^{-c}\leq c$, $c\geq0$,
  and Assumption~\ref{a:general:finite_man_hours}.
  Concavity of $k$ is inherited from the concavity of
  $x\mapsto 1-e^{-xc}$, $c\geq0$.
  Using dominated convergence
  together with Assumption~\ref{a:general:finite_man_hours},
  we see that
  \begin{equation} \label{gargamel}
    \frac{k(z)}{z}=
      \int \frac{1-\exp\rub{-z\int_0^\infty a\ru{\chi_s}\,ds}}{z}Q(d\chi)
      \xrightarrow{z\to\infty} 0.
  \end  {equation}
  In addition, dominated convergence
  together with Assumption~\ref{a:general:finite_man_hours}
  implies
  \begin{equation}
    k^{'}(z)=\int \eckB{ \int_0^\infty a\rub{\chi_s}\,ds\,
      \exp\ruB{-z\int_0^\infty a\rub{\chi_s}\,ds}}Q(d\chi)\,\quad z\geq0.
  \end  {equation}
  Hence, $k$ is strictly concave.
  Thus, $k$ has a fixed-point which is not zero
  if and only if $k^{'}(0)>1$.
  The implications~\eqref{aa:1} and~\eqref{aa:2} follow from the
  strict concavity of $k$.
\end  {proof}
The method of proof (cf. Section 6.5 in~\cite{Ja75}) of the extinction result for a
Crump-Mode-Jagers process $(J_t)_{t\geq0}$
is to study an equation for $\ru{\E e^{-\ld J_t}}_{t\geq0,\ld\geq0}$.
The Laplace transform $\ru{\E e^{-\ld J_t}}_{\ld\geq0}$ converges
monotonically to $\P(J_t=0)$ as $\ld\to\infty$, $t\geq0$.
Furthermore, $\P(J_t=0)=\P(\exists s\leq t\colon J_s=0)$
converges monotonically to the extinction probability $\P(\exists s\geq0\colon J_s=0)$
as $t\to\infty$.
Taking monotone limits in the equation for
$\ru{\E e^{-\ld J_t}}_{t\geq0,\ld\geq0}$ results in an equation for the extinction
probability. In our situation, there is an equation for
the modified Laplace transform $(L_t(\ld))_{t>0,\ld>0}$
as defined in~\eqref{eq:def:L_t} below. However, the monotone limit of $L_t(\ld)$
as $\ld\to\infty$ might be infinite.
Thus, it is not clear how to transfer the above method of proof.
The following proof of Theorem~\ref{thm:extinction_VIM} directly 
establishes the convergence of the modified Laplace transform.
\begin{proof}[\textbf{\upshape Proof of Theorem~\ref{thm:general:extinction}}]
  Recall $q$ from
  Lemma~\ref{l:fixed_point}.
  In the first step, we will prove
  \begin{equation}  \label{eq:def:L_t}
    L_t:=L_t(\ld):=\I\eckb{1-e^{-\ld V^{\phi,Q,Q}_{t}}} \to q\qqastk
  \end  {equation}
  for all $\ld>0$. Set $L_t(0):=0$.
  It follows from Lemma~\ref{l:general:finite_moments} that
  $(L_t)_{t\leq T}$ is bounded for every finite $T$.
  Lemma~\ref{l:integral_laplace_trafo} with $\nu$ replaced by $Q$
  provides us with the fundamental equation
  \begin{equation} \label{eq:darma}
    L_t
     =\int\E\eckB{1-
        \exp\ruB{-\ld\phi_{\chi}(t)-\int_0^\infty
          a\rub{\chi_s}
          L_{t-s}
          \,ds}}Q(d\chi)\quad\fa t\geq0.
  \end  {equation}
  Based on~\eqref{eq:darma},
  the idea for the proof of~\eqref{eq:def:L_t} is as follows.
  The term $\ld \phi_\chi(t)$ vanishes
  as $t\to\infty$.
  If $L_t$ converges to some limit, then the limit has to be a fixed-point
  of the function
  \begin{equation}
    k(z)=\int \eckB{1-\exp\ruB{-z\int_0^\infty a\rub{\chi_s}\,ds}}Q(d\chi).
  \end  {equation}
  By Lemma~\ref{l:fixed_point}, this function is (typically strictly) concave.
  Therefore, it has exactly one
  attracting fixed-point. Furthermore, this fact forces $L_t$ to converge
  as $t\to\infty$.

  We will need finiteness of $L_\infty:=\limsup_{t\to\infty} L_t$.
  Looking for a contradiction,
  we assume $L_\infty=\infty$.
  Then there exists a sequence $(t_n)_{n\in\N}$
  with $t_n\to\infty$ such that
  $L_{t_n}\leq \sup_{t\leq t_n}L_t\leq L_{t_n}+1$.
  We estimate
  \begin{equation}  \begin{split}   \label{eq:drama}
    L_{t_n}
    &\leq \int \eckB{1-\E\exp\ruB{-\ld \phi_\chi(t_n)-
       \int_0^\infty a\rub{\chi_s}\sup_{r\leq t_n}L_r\,ds}}Q(d\chi)\\
    &\leq k\ru{L_{t_n}\!+1}+\int \exp\ruB{-
       \int_0^\infty a\rub{\chi_s} L_{t_n}\,ds}
       \ruB{1-\E e^{-\ld\phi_\chi(t_n)}}Q(d\chi)\\
    &\leq k\ru{L_{t_n}\!+1}+\int \ruB{1-\E e^{-\ld\phi_\chi(t_n)}}Q(d\chi).
  \end  {split}     \end  {equation}
  The last summand converges to zero by
  Assumption~\ref{a:convergence_to_zero} and is therefore bound\-ed
  by some constant $c$.
  Inequality~\eqref{eq:drama} leads to the contradiction
  \begin{equation}
    1\leq \limn\frac{k(L_{t_n}+1)}{L_{t_n}}+\limn\frac{c}{L_{t_n}}=0.
  \end  {equation}
  The last equation is a consequence of~\eqref{gargamel}
  and the
  assumption $L_\infty=\infty$.
  Next
  we prove $ L_\infty\leq q$
  using boundedness of $(L_t)_{t\geq0}$.
  Let $(t_n)_{n\in\N}$ be a sequence
  such that $\tlimn L_{t_n}=L_\infty<\infty$.
  Then a calculation as in~\eqref{eq:drama} results in
  \begin{equation}  \begin{split}         \label{eq:duggy}
    \limn L_{t_n}
       &\leq\limsupn \int \eckB{1-\exp\ruB{-
       \int_0^\infty a\rub{\chi_s}\sup_{t\geq t_n} L_{t-s}\,ds}} Q(d\chi)\\
       &\qquad+\limsupn\int \ruB{1-\E e^{-\ld\phi_\chi(t_n)}}Q(d\chi).
  \end  {split}     \end  {equation}
  The last summand is equal to zero by
  Assumption~\ref{a:convergence_to_zero}.
  The first summand on the right-hand side of~\eqref{eq:duggy}
  is dominated by
  \begin{equation}
    \ruB{\sup_{t>0} L_t}\int\ruB{\int_0^\infty a\rub{\chi_s}\,ds}Q(d\chi)
    <\infty
  \end  {equation}
  which is finite by boundedness of $(L_t)_{t\geq0}$ and
  by Assumption~\ref{a:general:finite_man_hours}.
  Applying dominated convergence,
  we conclude that $L_\infty$ is bounded by
  \begin{equation}
    L_\infty\leq \int\eckB{1-\exp\ruB{
       -\int_0^\infty a\rub{\chi_s}
         \limsup_{t\to\infty} L_{t-s}\,ds}}Q(d\chi)
    =k\rub{L_\infty}.
  \end  {equation}
  Thus, Lemma~\ref{l:fixed_point} implies
  $\tlimsupt L_t\leq q$.

  Assume $q>0$ and suppose that $m:=\tliminft L_t=0$.
  Let $(t_n)_{n\in\N}$ be such that
  $0<L_{t_n}\geq\inf_{1\leq t\leq t_n}L_t\geq cL_{t_n}\to0$
  as $n\to\infty$ and
  $t_n+1\leq t_{n+1}\to\infty$.
  By Lemma~\ref{l:fixed_point},
  there is an $n_0$ and a $c<1$ such that
  $c\int\int_0^{t_{n_0}} a\rub{\chi_s}\,ds Q(d\chi)>1$.
  We estimate
  \begin{equation}  \begin{split}
    L_{t_n}&\geq\int\eckB{1-\exp\ruB{
       -\int_0^{t_n-1}a\rub{\chi_s}\inf_{1\leq t\leq t_n}L_{t}\,ds}}Q(d\chi)\\
   &\geq\int\eckB{1-\exp\ruB{
       - c\int_0^{t_{n_0}}a\rub{\chi_s} L_{t_n}\,ds}}Q(d\chi)\quad\fa n>n_0.
  \end  {split}     \end  {equation}
  Using dominated convergence,
  the assumption $m=0$ results in the contradiction
  \begin{equation}  \begin{split}
    1&\geq \limn\frac{1}{L_{t_n}}
     \int\eckB{1-\exp\ruB{- c L_{t_n}\int_0^{t_{n_0}} a\rub{\chi_s}\,ds}}
         Q(d\chi)\\
     &=c\int\ruB{\int_0^{t_{n_0}} a\rub{\chi_s}\,ds }Q(d\chi)>1.
  \end  {split}     \end  {equation}
  In order to prove $m\geq q$, let 
  $(t_n)_{n\in\N}$ be such that
  $\tlimn L_{t_n}=m>0$.
  An estimate as above together with dominated convergence yields
  \begin{equation}  \begin{split}
    m&=\limn L_{t_n}
    \geq\limn\int\eckB{1-\exp\ruB{
       -\int_0^{t_n}a\rub{\chi_s}\inf_{t\geq t_n}L_{t-s}\,ds}}Q(d\chi)\\
    &= \int\eckB{1-\exp\ruB{
      -\int_{0}^\infty a\rub{\chi_s}\liminft\, L_t\,ds}}Q(d\chi)
    =k(m).
  \end  {split}     \end  {equation}
  Therefore, Lemma~\ref{l:fixed_point}
  implies $\tliminft L_t=m\geq q$,
  which yields~\eqref{eq:def:L_t}.

  Finally, we finish the proof of Theorem~\ref{thm:general:extinction}.
  Applying Lemma~\ref{l:integral_laplace_trafo}, we see that
  \begin{equation}  \begin{split}  \label{eq:weseethat}
    \lefteqn{\absbb{\I\eckB{1-e^{-\ld V_t^{\phi,\nu,Q}}}
     -\int\eckB{1-\exp\ruB{-q\int_0^\infty a(\chi_s)\,ds}}\nu(d\chi)}}\\
    &\leq\int\exp\ruB{-\int_0^\infty L_{t-s}a\ru{\chi_s}\,ds}
      \E\eckB{1- e^{-\ld\phi_\chi(t)}}\nu(d\chi)\\
    &\qquad\,+\left|\int\eckB{1-\exp\ruB{-\int_0^\infty
      L_{t-s} a\ru{\chi_s}\,ds}}\nu(d\chi)\right.\\
    &\qquad\,
      -\left.\int\eckB{1-\exp\ruB{-q\int_0^\infty a\ru{\chi_s}\,ds}}
      \nu(d\chi)\right|.
  \end  {split}     \end  {equation}
  The first summand on the right-hand side of~\eqref{eq:weseethat} converges
  to zero as $t\to\infty$ by Assumption~\ref{a:convergence_to_zero}.
  By the first step~\eqref{eq:def:L_t}, $L_t\to q$ as $t\to\infty$.
  Hence, by the dominated convergence theorem
  and Assumption~\ref{a:general:finite_man_hours},
  the left-hand side of~\eqref{eq:weseethat} converges to zero as
  $t\to\infty$.
  As  $\nu$ is a probability measure by assumption, we conclude
  \begin{equation}
    \limt \E e^{-\ld V_t^{\phi,\nu,Q}}
    =\int \exp\ruB{-q\int_0^\infty a\rub{\chi}\,ds}\nu(d\chi)
    \quad\fa\ld\geq0.
  \end  {equation}
  This implies Theorem~\ref{thm:general:extinction}
  as the Laplace transform is
  convergence determining, see e.g.\ Lemma 2.1 in~\cite{Dyn89}.
\end  {proof}

%
%
\noindent
\section{The supercritical Virgin Island Model.
Proof of Theorem~\ref{thm:general:xlogx}}%
\label{sec:convergence_supercritical}
Our proof of Theorem~\ref{thm:general:xlogx} follows
the proof of Doney (1972)~\cite{Don72} for supercritical
Crump-Mode-Jagers processes.
Some changes are necessary because
the recursive equation~\eqref{eq:darma} differs from the
respective recursive equation in~\cite{Don72}.
Parts of our proof are analogous to the proof in~\cite{Don72}
which we nevertheless include here for the reason of completeness.
Lemma~\ref{l:lemma52} and Lemma~\ref{l:convergence_if_xlogx_fails}
below contain the essential part of the proof of
Theorem~\ref{thm:general:xlogx}.
For these two lemmas, we will need auxiliary lemmas which
we now provide.

We assume throughout this section that a solution $\al\in\R$
of equation~\eqref{eq:Malthusian} exists. Note that this is
implied by~\ref{a:general:finite_man_hours} and
$Q\rub{\int_0^\infty a(\chi_s)\,ds>0}>0$.
Recall the definition of $\mu^Q$ from~\eqref{eq:def:mu_nu}.

\subsection{Preliminaries}
For $\ld\geq0$, define
\begin{equation}
  H_{\al}(\psi)(\ld):=\int\eckB{1-\exp\ruB{-
     \int_0^\infty a\rub{\chi_s}\psi(\ld e^{-\al s})\,ds}}Q(d\chi)
\end  {equation}
for $\psi\in\Dpos$.
\begin{lemma} \label{l:H_contracting}
  The operator $H_{\al}$ is contracting in the sense that
  \begin{equation}
    \absb{H_{\al}(\psi_1)(\ld)-H_{\al}(\psi_2)(\ld)}
    \leq\int_0^\infty
      \absb{\psi_1(\ld e^{-\al s})-\psi_2(\ld e^{-\al s})}\mu^Q(ds)
  \end  {equation}
  for all $\psi_1,\psi_2\in\Dpos$.
\end  {lemma}
\begin{proof}
  The lemma follows immediately from $\abs{e^{-x}-e^{-y}}\leq\abs{x-y}$
  and from the definition~\eqref{eq:def:mu_nu} of $\mu^Q$.
\end  {proof}
%
%
%
\begin{lemma} \label{l:H_nondecreasing}
  The operator $H_{\al}$ is nondecreasing in the sense that
  \begin{equation}
    H_{\al}(\psi_1)(\ld)\leq H_{\al}(\psi_2)(\ld)
  \end  {equation}
  for all $\ld\geq0$ if $\psi_1(\ld)\leq\psi_2(\ld)$ for all $\ld\geq0$.
\end  {lemma}
\begin{proof}
  The lemma follows from
  $1-e^{-cx}$ being increasing in $x$ for every $c>0$.
\end  {proof}

%
%
%
For every measurable function $\psi\colon\R\times[0,\infty)\to[0,\infty)$,
define
\begin{equation}  \begin{split}  \label{eq:Hbar_reads_as}
    \bar{H}_{\al}(\psi)(t,\ld)
     :=\int\eckbb{f\ruB{\int_0^\infty a\rub{\chi_s}\psi(t-s,\ld e^{-\al s})ds}}
          Q(d\chi).
\end  {split}     \end  {equation}
for $\ld\geq0$ and $t\in\R$
where $f(x):=x-1+e^{-x}\geq0$, $x\geq0$.
If $\tilde{\psi}\colon[0,\infty)\to[0,\infty)$ is a function
of one variable, then we set
$\bar{H}_\al(\tilde{\psi})(\ld):= \bar{H}_\al(\psi)(1,\ld)$
where $\psi(t,\ld):=\tilde{\psi}(\ld)$ for $\ld\geq0$, $t\in\R$.
\begin{lemma} \label{l:barH_nondecreasing}
  The operator $\bar{H}_{\al}$ is nondecreasing in the sense that
  \begin{equation}
    \bar{H}_{\al}(\psi_1)(t,\ld)\leq \bar{H}_{\al}(\psi_2)(t,\ld)
  \end  {equation}
  for all $\ld\geq0$ and $t\in\R$
  if $\psi_1(t,\ld)\leq\psi_2(t,\ld)$ for all $\ld\geq0$, $t\in\R$.
\end  {lemma}
\begin{proof}
  The assertion follows from the basic fact that $f$ is nondecreasing.
\end  {proof}

%
%
%
\begin{lemma}  \label{l:eta_nondecreasing}
  Assume~\ref{a:general:finite_man_hours}.
  Let $id:\ld\mapsto\ld$ be the identity map.
  The function
  \begin{equation}  \label{eq:def:eta}
    \eta(\ld):=1-\frac{1}{\ld}H_{\al}(id)(\ld)
    =\frac{1}{\ld}\bar{H}_{\al}(id)(\ld),\quad \ld>0,
  \end  {equation}
  is nonnegative and nondecreasing.
  Furthermore, $\eta(0+)=0$.
\end  {lemma}
\begin{proof}
  Recall the definition of $A_\al(\chi)$ from~\eqref{eq:def:A_al}.
  By equation~\eqref{eq:Hbar_reads_as},
  we have $\ld \eta(\ld)=\int f(\ld A_\al)\,dQ$.
  Thus, $\eta$ is nonnegative.
  Furthermore, $\eta(0+)=0$ follows from the dominated convergence theorem
  and Assumption~\ref{a:general:finite_man_hours}.
  Let $x,y>0$. Then
  \begin{equation}
    \eta(x+y)-\eta(x)
      =\int\frac{x A_\al f\rub{(x+y)A_\al}-(x+y)A_\al f\rub{x A_\al}}
         {x(x+y)A_\al}dQ\geq0.
  \end  {equation}
  The inequality follows from
  $\xt f(\xt+\yt)-(\xt+\yt)f(\xt)\geq0$ for all $\xt,\yt\geq0$.
\end  {proof}

%
%
%
The following lemma, due to Athreya~\cite{Ath69},
translates the $(x \log x)$-condition~\eqref{eq:xlogx}
into an integrability condition on $\eta$ . 
For completeness, we include its proof.
\begin{lemma}  \label{l:properties_of_eta}
  Assume~\ref{a:general:finite_man_hours}.
  Let $\eta$ be the function defined in~\eqref{eq:def:eta}.
  Then
  \begin{equation}  \label{eq:two_quantities}
    \int_{0+}\frac{1}{\ld}\eta(\ld)\,d\ld<\infty\text{ and }
    \sum_{n=1}^\infty\eta(c r^n)<\infty
  \end  {equation}
  for some and then all $c>0$, $r<1$
  if and only if the $(x \log x)$-condition~\eqref{eq:xlogx} holds.
\end  {lemma}
\begin{proof}
  By monotonicity of $\eta$ (see
  Lemma~\ref{l:eta_nondecreasing}), the two quantities
  in~\eqref{eq:two_quantities} are finite or infinite at the same time.
  Fix $c>0$.
  Using Fubini's theorem and the substitution $v:=\ld A_\al$, we obtain
  \begin{equation}  \begin{split}
    \int_0^c\frac{1}{\ld}\eta(\ld)\,d\ld
    &=\int\eckbb{\int_0^c\eckB{
                   \frac{\ld A_\al-1+e^{-\ld A_\al}}{(\ld A_\al)^2}}\rub{A_\al}^2\,d\ld}dQ\\
    &=\int\eckbb{A_\al\int_0^{cA_\al}\frac{v-1+e^{-v}}{v^2}
                   \,dv}dQ.
  \end  {split}     \end  {equation}
  It is a basic fact that $\int_0^T\tfrac{1}{v^2}(v-1+e^{-v})dv\sim\log T$
  as $T\to\infty$.
\end  {proof}

\subsection{The limiting equation}
%
%
%
In the following two lemmas, we consider uniqueness and
existence of a function $\Psi$ which satisfies:
\begin{equation}  \begin{split}  \label{eq:conditions_on_Psi}
    &\text{(a)}\quad\absb{\Psi(\ld_1)-\Psi(\ld_2)}\leq\abs{\ld_1-\ld_2}
         \text{ for }\ld_1,\ld_2\geq0,\ \Psi(0)=0\\
    &\text{(b)}\quad\Psi(\ld)=H_{\al}(\Psi)(\ld)\\
    &\text{(c)}\quad\frac{\Psi(\ld)}{\ld}\to 1\qqasldO\\
    &\text{(d)}\quad0\leq\Psi(\ld_1)\leq\Psi(\ld_2)\leq\ld_2\ \fa0\leq\ld_1\leq\ld_2
             \text{ and }\limld\Psi(\ld)=q
\end  {split}     \end  {equation}
where $q\geq0$ is as in Lemma~\ref{l:fixed_point}.
Notice that the zero function does not satisfy~\eqref{eq:conditions_on_Psi}(c).
First, we prove uniqueness.
\begin{lemma}  \label{l:uniqueness_Psi}
  Assume~\ref{a:general:finite_man_hours} and $\al>0$.
  If $\Psi_1$ and $\Psi_2$
  satisfy \eqref{eq:conditions_on_Psi}, then $\Psi_1=\Psi_2$.
\end  {lemma}
\begin{proof}
  Notice that $\Psi_1(0)=\Psi_2(0)$.
  Define $\Lambda(\ld):=\tfrac{1}{\ld}\abs{\Psi_1(\ld)-\Psi_2(\ld)}$
  for $\ld>0$ and note that  $\Lambda(0+)=0$
  by~\eqref{eq:conditions_on_Psi}(c).
  From Lemma~\ref{l:H_contracting}, we obtain for $\ld>0$
  \begin{equation}  \begin{split}  \label{eq:inequality_for_Lambda}
    \Lambda(\ld)\leq\frac{1}{\ld}\int_0^\infty
      \absb{\Psi_1(\ld e^{-\al s})-\Psi_2(\ld e^{-\al s})}\mu^Q(ds)
      =\int_0^\infty\Lambda(\ld e^{-\al s})\mu^Q_\al(ds)
  \end  {split}     \end  {equation}
  where $\mu_\al^Q(ds):=e^{-\al s}\mu^Q(ds)$
  is a probability measure because
  $\al$ solves equation~\eqref{eq:general:Malthusian}.
  Let $R_i$, $i\geq 1$, be independent variables with distribution
  $\mu_\al^Q$ and note that $\E R_1<\infty$.
  We may assume that $\E R_1>0$ because $\mu^Q\rub{[0,\infty)}=0$
  implies $\Psi_i=H_\al(\Psi_i)=0$ for $i=1,2$.
  Iterating inequality~\eqref{eq:inequality_for_Lambda},
  we arrive at
  \begin{equation}  \begin{split} \label{eq:iterating_slln}
    \Lambda(\ld)\leq\E\Lambda\rub{\ld e^{-\al R_1}}
    \leq\E\Lambda\rub{\ld e^{-\al (R_1+\ldots +R_n)}}
    \lra\Lambda(0+)=0\qasn.
  \end  {split}     \end  {equation}
  The convergence in~\eqref{eq:iterating_slln} follows from
  the weak law of large numbers.
\end  {proof}
%
%
%
\begin{lemma}  \label{l:existence_Psi}
  Assume~\ref{a:general:finite_man_hours} and $\al>0$.
  There exists a solution $\Psi$
  of~\eqref{eq:conditions_on_Psi} if and only
  if the $(x \log x)$-condition~\eqref{eq:xlogx} holds.
\end  {lemma}
\begin{proof}
  Assume that~\eqref{eq:xlogx} holds. 
  Define $\Psi_0(\ld):=\ld$, $\Psi_{n+1}(\ld):=H_{\al}(\Psi_n)(\ld)$
  for $\ld\geq0$
  and $\Lambda_{n+1}(\ld):=\tfrac{1}{\ld}\absb{\Psi_{n+1}(\ld)-\Psi_n(\ld)}$
  for $\ld>0$ and $n\geq 0$.
  Recall $\mu_\al^Q$ and $(R_i)_{i\in\N}$ from the proof of
  Lemma~\ref{l:uniqueness_Psi}.
  Note that $\E R_1>0$ because of $\al>0$.
  Arguments as in the proof of Lemma~\ref{l:uniqueness_Psi} imply
  \begin{equation} \label{eq:estimate_for_Lambda_n}
    \Lambda_{n+1}(\ld)\leq \E\Lambda_n\rub{\ld e^{-\al R_1}}
    \leq\E \Lambda_1\rub{\ld e^{-\al S_n}}
  \end  {equation}
  where $S_n:=R_1+\ldots +R_n$ for $n\geq0$.
  Since $\eta\geq0$ by
  Lemma~\ref{l:eta_nondecreasing} and
  \begin{equation}
    \Psi_0(\ld)-\Psi_1(\ld)=\ld -H_{\al}(id)(\ld)=\ld\eta(\ld),
  \end  {equation}
  we see that $\eta=\Lambda_1$.
  In addition, we conclude from $\eta\geq0$ that
  $\Psi_1(\ld)\leq\Psi_0(\ld)=\ld$.
  By Lemma~\ref{l:H_nondecreasing},
  this implies inductively $\Psi_n(\ld)\leq \ld$
  for $n\geq0$, $\ld\geq0$.
  Let $\Lambda(\ld):=\sum_{n\geq1}\Lambda_n(\ld)$.
  We need to prove that $\Lambda(\ld)<\infty$. Clearly $0<\E e^{-R_1}<1$,
  so we can choose $\eps>0$ with $e^\eps\E e^{-R_1}<1$. Then
  \begin{equation}  \label{eq:summe_ueber_Sn}
    \sum_{n=0}^\infty\P\rub{S_n\leq n\eps}
    \leq
    \sum_{n=0}^\infty e^{n\eps}\E e^{-S_n}
    =
    \sum_{n=0}^\infty \ruB{e^{\eps}\E e^{-R_1}}^n<\infty.
  \end  {equation}
  Define $\bar{\eta}(\ld):=\sup_{0<u\leq \ld}\eta(u)$.
  It follows
  from \eqref{eq:estimate_for_Lambda_n}, \eqref{eq:summe_ueber_Sn},
  Lemma~\ref{l:eta_nondecreasing}
  and Lemma \ref{l:properties_of_eta} that for all $\ld>0$
  \begin{equation}  \begin{split}  \label{eq:for_Lambda}
    \Lambda(\ld)&\leq\sum_{n=0}^\infty\E\eta\rub{\ld e^{-\al S_n}}
    \leq \bar{\eta}(\ld)
     \sum_{n=0}^\infty\P(S_n\leq n\eps)
       +\sum_{n=0}^\infty \eta(\ld e^{-n\al\eps})<\infty.
  \end  {split}     \end  {equation}
  Thus, $(\Psi_n(\ld))_{n\geq0}$ is a Cauchy sequence in $[0,\ld]$.
  Hence, we conclude the existence of a function
  $\Psi$ such that $\Psi(\ld)=\limn\Psi_n(\ld)$ for every $\ld\geq0$.
  By the dominated convergence theorem, $\Psi$
  satisfies~\eqref{eq:conditions_on_Psi}(b).
  To check~\eqref{eq:conditions_on_Psi}(a), we prove that $\Psi_n$ is
  Lipschitz continuous with constant one. The induction step follows
  from Lemma~\ref{l:H_contracting}
  \begin{equation}  \begin{split}
    \absb{\Psi_{n+1}(\ld_1)-\Psi_{n+1}(\ld_2)}
    &\leq \int_0^\infty\absb{\Psi_n(\ld_1 e^{-\al s})-\Psi_n(\ld_2 e^{-\al s})}
      \mu^Q(ds)\\
    &\leq \absb{\ld_1-\ld_2}\int_0^\infty e^{-\al s}\mu^Q(ds)
    =\absb{\ld_1-\ld_2}.
  \end  {split}     \end  {equation}
  In order to
  check~\eqref{eq:conditions_on_Psi}(c), note that since 
  $\eta(0+)=0$, it follows from~\eqref{eq:for_Lambda} that $\Lambda(0+)=0$.
  Thus,
  \begin{equation}
    \absb{\frac{\Psi(\ld)}{\ld}-1}\leq
    \limsupn\frac{1}{\ld}\absb{\Psi_n(\ld)-\ld}\leq\Lambda(\ld)
    \lra0\qqasldO,
  \end  {equation}
  as required.
  Finally, monotonicity of $\Psi_n$ and $\Psi_n(\ld)\leq\ld$ for all $n\in\N$
  imply monotonicity of $\Psi$ and $\Psi(\ld)\leq\ld$, respectively.
  The last claim of~\eqref{eq:conditions_on_Psi}(d), namely $\Psi(\infty)=q$,
  follows from letting $\ld\to\infty$ in
  $\Psi(\ld)=H_\al(\Psi)(\ld)$, monotonicity of $\Psi$
  and from Lemma~\ref{l:fixed_point} together with
  $\Psi(\infty)>0$.
 
  For the ``only if''-part of the lemma,
  suppose that there exists a solution
  $\Psi$ of~\eqref{eq:conditions_on_Psi}.
  Write $\gt(\ld):=\tfrac{\Psi(\ld)}{\ld}$ for $\ld>0$.
  Since $\gt\geq0$ and $\gt(0+)=1$, there exist constants $c_1,c_2,c_3>0$
  such that $c_2\leq \gt(\ld)\leq c_3$ for all $\ld\in(0,c_1]$.
  Using~\eqref{eq:conditions_on_Psi}(b),
  $\Psi(\ld)\geq \ld c_2$ for $\ld\in(0,c_1]$
  and Lemma~\ref{l:barH_nondecreasing},
  we obtain for $\ld\in(0,c_1]$
  \begin{equation}  \begin{split}
    \gt(\ld)&=\frac{H_{\al}(\Psi)(\ld)}{\ld}
    =\frac{1}{\ld}\int_0^\infty
        \Psi(\ld e^{-\al s})\mu^Q(ds)-\frac{1}{\ld}\bar{H}_{\al}(\Psi)(\ld)\\
    &\leq \int_0^\infty \gt\rub{\ld e^{-\al s}}\mu_\al^Q(ds)-
       \frac{1}{\ld}\bar{H}_{\al}\ru{c_2  \cdot}(\ld)
    =\E \gt\rub{\ld e^{-\al R_1}}-c_2\eta(c_2\ld).
  \end  {split}     \end  {equation}
  Let $t_0$ be such that $0<c_4:=\mu_\al^Q\rub{[0,t_0]}<1$
  and write $\gt^\ast(\ld):=\sup_{u\leq \ld}\gt(u)$. Then
  \begin{equation}
    \gt^\ast(\ld)\leq c_4 \gt^\ast(\ld)+(1-c_4)\gt^\ast\rub{\ld e^{-\al t_0}}
      -c_2\eta(c_2\ld)
  \end  {equation}
  which we rewrite as
  $\gt^\ast(\ld)\leq \gt^\ast(\tau\ld)-c_5\eta(c_2\ld)$
  where $\tau:=e^{-\al t_0}$ and $c_5:=\tfrac{c_2}{1-c_4}$.
  Iterating this inequality results in
  $\gt^\ast(\ld)\leq \gt^\ast(\ld \tau^{n+1})
      -c_5\sum_{k=0}^n\eta(c_2\ld\tau^k)$
  for $n\geq0$.
  Since $\gt^\ast$ is bounded on $(0,c_1]$ this
  implies that
  $\sum_{k=0}^\infty\eta(c_2\ld\tau^k)<\infty$.
  Therefore, by Lemma~\ref{l:properties_of_eta}, the $(x \log x)$-condition
  holds.
\end  {proof}

\subsection{Convergence}
Recall $\bm$, $\I$, $m^Q$ and $L_t$ from~\eqref{eq:def:barm},
\eqref{eq:def:I},
\eqref{eq:def:mu_nu}
and~\eqref{eq:def:L_t}, respectively.
As before, let $\mu_\al^Q(ds):=e^{-\al s}\mu^Q(ds)$.
Define
\begin{equation}  \label{eq:D_ld_t}
  D(\ld,t):=\frac{m^Q(t)}{e^{\al t}\bm}
            -\frac{1}{\ld}L_t\ruB{\frac{\ld}{e^{\al t}\bm}},
\end  {equation}
$D_T(\ld):=\sup_{s\leq T}\abs{D(\ld,s)}$ and
$D_\infty(\ld):=\lim_{T\to\infty}D_T(\ld)$
for $\ld>0$ and $t,T\geq0$.
The following two lemmas follow Lemma 5.1 and Lemma 5.2, respectively,
in~\cite{Don72}.
%
%
%
\begin{lemma}  \label{l:lemma51}
  Assume~\ref{a:general:finite_moments},
  \ref{a:general:finite_man_hours},
  \ref{a:for_SKDist} and
  $\al>0$.
  If the $(x \log x)$-condition \eqref{eq:xlogx} holds, then
  $D_\infty(\ld)\to0$ as $\ld\to0$.
\end  {lemma}
\begin{proof}
  Inserting the definitions~\eqref{eq:def:mu_nu} and~\eqref{eq:def:L_t}
  of $m^Q$ and $L_t$, respectively,
  into~\eqref{eq:D_ld_t}, we see that
  \begin{equation}  \label{eq:D_pos}
    D(\ld,t)=\frac{1}{\ld}
       \I\eckB{f\ruB{\frac{\ld V^{\phi,Q,Q}_t}{e^{\al t}\bm}}}\geq  0
       \qquad \ld,t>0
  \end  {equation}
  is nonnegative
  where $f(x):=x-1+e^{-x}$, $x\geq0$.
  Insert the recursive equations~\eqref{eq:recursive_m_Q}
  and~\eqref{eq:darma} for $m^Q$ and $(L_t)_{t\geq0}$, respectively,
  into~\eqref{eq:D_ld_t} to obtain
  \begin{equation}  \begin{split}  \label{eq:estimate_D_ld_t}
    \lefteqn{D(\ld,t)
    =\int\frac{\E\phi(t,\chi)}{e^{\al t}\bm}\,Q(d\chi)
     +\int_0^t \frac{m^Q(t-s)}{e^{\al(t-s)}\bm}\mu^Q_\al(ds)}\\
     &\quad -\frac{1}{\ld}\int\eckB{1-
      \E\exp\ruB{-\ld\frac{\phi_{\chi}(t)}{e^{\al t}\bm}-\int_0^t a\rub{\chi_s}
       L_{t-s}\rub{\frac{\ld}{e^{\al t}\bm}}\,ds}}Q(d\chi)\\
    &= \int_0^t\eckB{ \frac{m^Q(t-s)}{e^{\al(t-s)}\bm}
      -\frac{1}{\ld e^{-\al s}}
       L_{t-s}\ruB{\frac{\ld e^{-\al s}}{e^{\al(t-s)}\bm}}} \mu^Q_\al(ds)\\
    &\quad +\frac{1}{\ld }\int_0^t 
       L_{t-s}\ruB{\frac{\ld e^{-\al s}}{e^{\al (t-s)}\bm}} \mu^Q(ds)\\
    &\qquad\quad-\frac{1}{\ld}\int\eckB{1-
      \exp\ruB{-\int_0^t a\rub{\chi_s}
       L_{t-s}\rub{\frac{\ld e^{-\al s}}{e^{\al (t-s)}\bm}}\,ds}}Q(d\chi)\\
    &\quad + \int
      \E\eckB{\frac{1-\exp\rub{-\ld\frac{\phi_\chi(t)}{e^{\al t}\bm}}}{\ld}}
    \eckB{1- \exp\ruB{-\int_0^t a\rub{\chi_s}
       L_{t-s}\rub{\frac{\ld}{e^{\al t}\bm}}\,ds}}Q(d\chi)\\
    &\quad + \frac{1}{\ld}\int
    \E\eckbb{\frac{\ld\phi(t,\chi)}{e^{\al t}\bm}-1
    +\exp\ruB{-\frac{\ld\phi(t,\chi)}{e^{\al t}\bm}}}  Q(d\chi)\\
    &=:\int_0^t D(\ld e^{-\al s},t-s)\mu^Q_\al(ds)
     +\frac{1}{\ld}\bar{H}_\al\ruB{(t,\ld)\mapsto
              L_t\rub{\frac{\ld}{e^{\al t}\bm}}}
     +T_1+T_2
  \end  {split}     \end  {equation}
  where $T_1$ and $T_2$ are suitably defined.
  Inequality~\eqref{eq:D_pos} implies
  \begin{equation} \label{eq:L_bdd_by_id}
    L_t\ruB{\frac{\ld}{e^{\al t}\bm}}\leq \ld \frac{m^Q(t)}{e^{\al t}\bm}
    \leq\ld c_1\qquad t,\ld\geq0
  \end  {equation}
  where $c_1$ is a finite constant. The last inequality is a
  consequence of Theorem~\ref{thm:general:expected_man_hours},
  equation~\eqref{eq:general:rescaled_VIM_supercritical},
  with $\nu$ replaced by $Q$.
  Lemma~\ref{l:barH_nondecreasing} together
  with~\eqref{eq:L_bdd_by_id} implies
  \begin{equation}  \label{eq:estimate_barH}
    \frac{1}{\ld}\bar{H}_\al\ruB{ (t,\ld)\mapsto
    L_t\rub{\frac{\ld}{e^{\al t}\bm}}}
    \leq \frac{1}{\ld}\bar{H}_\al\ruB{c_1\mal id}
    =c_1\eta(\ld c_1).
  \end  {equation}
  Using $1-e^{-x}\leq x$, inequality \eqref{eq:L_bdd_by_id}
  and $x-1+e^{-x}\leq\tfrac{1}{2}x^2$, $x\geq0$,
  we see that the expressions $T_1$ and $T_2$ are bounded above by
  \begin{equation}  \begin{split}  \label{eq:inequality_T}
    T_1&\leq\int\eckbb{\frac{\E\phi_\chi(t)}{e^{\al t}\bm}
     \int_0^t a\rub{\chi_s}\ld e^{-\al s} c_1 ds}Q(d\chi)
     \leq c_2\ld\\
    T_2&\leq\frac{\ld}{2}\int\E\ruB{\frac{\phi_\chi(t)}{e^{\al t}\bm}}^2 Q(d\chi)
    \leq c_3 \ld
  \end  {split}     \end  {equation}
  for all $\ld,t>0$ where $c_2, c_3$ are finite constants which
  are independent of $t>0$ and $\ld>0$.
  Such constants exist by Assumption~\ref{a:for_SKDist}.
  Taking supremum over $t\in[0,T]$ 
  in~\eqref{eq:estimate_D_ld_t}
  and inserting~\eqref{eq:estimate_barH}
  and~\eqref{eq:inequality_T} results in
  \begin{equation}
    D_T(\ld)\leq\int_0^T D_T(\ld e^{-\al s})\mu_\al^Q(ds)
     + c_1\eta(\ld c_1)+c_2\ld+c_3\ld\qquad \fa\ld>0.
  \end  {equation}
  Choose $t_0>0$ such that $c_4:=\mu_\al^Q([0,t_0])\in(0,1)$.
  Then, by monotonicity of $D_T$,
  \begin{equation}
    D_T(\ld)\leq c_4 D_T(\ld)+(1-c_4)D_T(\ld e^{-\al t_0})
     + c_1\eta(\ld c_1)+c_2\ld+c_3\ld
  \end  {equation}
  for all $\ld>0$.
  Hence, $D_T$ is bounded by
  \begin{equation}
    D_T(\ld)\leq D_T(\ld e^{-\al t_0})+c_5\eta(c_1\ld)+c_6\ld
    \qquad \fa\ld>0
  \end  {equation}
  where $c_5:=\tfrac{c_1}{1-c_4}$ and $c_6:=\tfrac{c_2+c_3}{1-c_4}$.
  Iterate this inequality to obtain
  \begin{equation}  \begin{split} \label{eq:for_D_T}
    D_T(\ld)&\leq D_T(\ld e^{-\al t_0 n})
      +\sum_{k=0}^{n-1}\ruB{c_5\eta(c_1\ld e^{-\al t_0 k})
           +c_6\ld e^{-\al t_0 k}}\\
     &\lran D_T(0+)
      +\sum_{k=0}^{\infty}\ruB{c_5\eta(c_1\ld e^{-\al t_0 k})
           +c_6\ld e^{-\al t_0 k}}.
  \end  {split}     \end  {equation}
  Now we need to prove $D_T(0+)=0$.
  Looking at~\eqref{eq:D_pos} and using $f(x)\leq x^2/2$, we see that
  $D_T(\ld)$ is bounded by
  $\tfrac{\ld}{2\mb^2}\sup_{t\leq T}\I\eckb{\rub{V_t^{\phi,Q,Q}}^2}$.
  This is finite because of inequality~\eqref{eq:second_moment_V} with $\nu=Q$
  together with~\ref{a:for_SKDist}.
  Therefore $D_T(0+)=0$.
  Letting $T\to\infty$ in~\eqref{eq:for_D_T}, we obtain
  \begin{equation}  \begin{split}  \label{eq:D_infty}
    D_\infty(\ld)
      \leq c_5\sum_{k=0}^\infty \eta\rub{\ld c_1 e^{-\al t_0 k}}
         + \ld c_6\sum_{k=0}^\infty e^{-\al t_0 k}
         \quad \fa\ld>0.
  \end  {split}     \end  {equation}
  The right-hand side is finite by Lemma~\ref{l:properties_of_eta}.
  By Lemma~\ref{l:eta_nondecreasing}, we know that $\eta(0+)=0$
  and that $\eta$ is nondecreasing.
  Letting $\ld\to0$ in~\eqref{eq:D_infty} and using
  the dominated convergence theorem implies
  $D_\infty(\ld)\to 0$ as $\ld\to0$.
\end  {proof}
%
%
%
\begin{lemma}  \label{l:lemma52}
  Assume~\ref{a:general:finite_moments},
  \ref{a:general:finite_man_hours},
  \ref{a:convergence_to_zero},
  \ref{a:for_SKDist} and $\al>0$
  If  the  $(x \log x)$-condition \eqref{eq:xlogx} holds, then 
  \begin{equation}   \label{eq:convergence_xlogx}
    L_t\ruB{\frac{\ld}{\bm e^{\al t}}}\to\Psi(\ld)\qqast
  \end  {equation}
  for every $\ld\geq0$
  where $\Psi$ is the unique solution
  of~\eqref{eq:conditions_on_Psi}.
\end  {lemma}
\begin{proof}
  The case $\ld=0$ is trivial.
  For $\ld>0,t\geq0$, define
  \begin{equation}
    J(\ld,t):=
      \frac{1}{\ld}\rubb{L_t\rub{\frac{\ld}{\bm e^{\al t}}}-\Psi(\ld)}.
  \end  {equation}
  Furthermore, let $J_T(\ld):=\sup_{t\geq T}\abs{J(\ld,t)}$ and
  $J_\infty(\ld):=\lim_{T\to\infty}J_T(\ld)$ for $\ld>0$.
  We will prove $J_\infty(\ld)=0$ for $\ld>0$.
  By Theorem~\ref{thm:general:expected_man_hours}
  and~\eqref{eq:conditions_on_Psi}(c),
  \begin{equation}
    \absb{J(\ld,t)+D(\ld,t)}\leq\absb{\frac{m^Q(t)}{\bm e^{\al t}}-1}
      +\absb{\frac{\Psi(\ld)}{\ld}-1}
    \lrat \absb{\frac{\Psi(\ld)}{\ld}-1}
    \lraldO 0.
  \end  {equation}
  Hence, $J_\infty(0+)=0$ by Lemma~\ref{l:lemma51}.
  Using~\eqref{eq:darma}
  and~\eqref{eq:conditions_on_Psi}(b),
  we estimate
  \begin{equation}  \begin{split} \label{eq:lange_rechnung_1}
    \lefteqn{\Bigl| \ld J(\ld,2t)
     -\int 1-\E\exp\ruB{-\tfrac{\ld}{\bm e^{\al 2t}}\phi_\chi(2t)
       -\int_0^t a\rub{\chi_s} L_{2t-s}
           \rub{\tfrac{\ld}{\bm e^{\al 2 t}}} ds}Q(d\chi)}\\
   &\quad+\int 1-\exp\ruB{-\int_0^t a\rub{\chi_s}
     \Psi\rub{\ld e^{-\al s}}ds}Q(d\chi)\Bigr|\\
   &=\Bigl| 
     \int \E\exp\ruB{-\tfrac{\ld}{\bm e^{\al 2t}}\phi_\chi(2t)}
       \Bigl\{\exp\ruB{-\int_0^t a\rub{\chi_s}
          L_{2t-s}\rub{\tfrac{\ld}{\bm e^{\al 2t}}}ds}\\
   &\qquad\qquad\qquad\qquad\qquad\qquad
            -\exp\ruB{-\int_0^{2t} a\rub{\chi_s}
            L_{2t-s}\rub{\tfrac{\ld}{\bm e^{\al 2t}}} ds}\Bigr\}
       Q(d\chi)\\
   &\quad-\int \exp\ruB{-\int_0^t a\rub{\chi_s} \Psi\rub{\ld e^{-\al s}}ds}
   +\exp\ruB{-\int_0^\infty a\rub{\chi_s} \Psi\rub{\ld e^{-\al s}}ds}
    Q(d\chi)\Bigr|\\
   &\leq
    \int \int_t^\infty a\rub{\chi_s}
      \curlB{\sup_{u\geq0}L_u\rub{\tfrac{\ld}{\bm}}
      +\Psi\rub{\ld e^{-\al s}}}ds\,Q(d\chi)\\
   &\leq c\int\int_t^\infty  a\rub{\chi_s} ds\,Q(d\chi)
  \end  {split}     \end  {equation}
  for a suitable constant $c$.
  The last inequality uses boundedness of $(L_t)_{t\geq0}$, 
  see the proof of Theorem~\ref{thm:general:extinction},
  and of $\Psi$.
  By Assumption~\ref{a:general:finite_man_hours},
  the right-hand side of~\eqref{eq:lange_rechnung_1} converges
  to zero as $t\to\infty$.
  Fix $\ld>0$ and let $(t_n)_{n\geq1}$ be
  such that $\limn \abs{J(\ld,2t_n)}=J_\infty(\ld)$.
  With this, we get
  \begin{equation}  \begin{split} \label{eq:lange_rechnung_2}
    \Bigl|&\int 1-
       \E\exp\ruB{-\tfrac{\ld}{\bm e^{\al 2t_n}}\phi_\chi(2t_n)
       -\int_0^{t_n} a\rub{\chi_s} L_{2t_n-s}
         \rub{\tfrac{\ld}{\bm e^{\al 2{t_n}}}} ds}Q(d\chi)\\
   &\qquad -\int 1-\exp\ruB{-\int_0^{t_n} a\rub{\chi_s}
     \Psi\rub{\ld e^{-\al s}}ds}Q(d\chi)\Bigr|\\
   &\leq
     \int 1-\E\exp\ruB{-\tfrac{\ld}{\mb}\phi_\chi(2t_n)}Q(d\chi)\\
     &\qquad+ \int\int_0^{t_n}a\rub{\chi_s}
       \absB{L_{2t_n-s}\rub{\tfrac{\ld e^{-\al s}}{\bm e^{\al (2t_n-s)}}}
          -\Psi(\ld e^{-\al s})} ds \,Q(d\chi)\\
   &\leq
     \int 1-\E\exp\rub{-\tfrac{\ld}{\mb}\phi_\chi(2t_n)}Q(d\chi)
     + \ld\int_0^{t_n} J_{t_n}(\ld e^{-\al s}) \mu_\al^Q(ds)\\
   &\lran
     \ld\int_0^{\infty} J_{\infty}(\ld e^{-\al s}) \mu_\al^Q(ds).
  \end  {split}     \end  {equation}
  The convergence in~\eqref{eq:lange_rechnung_2} follows
  from~\ref{a:convergence_to_zero}
  and from the dominated convergence theorem
  together with Assumption~\ref{a:general:finite_man_hours}.
  Recall $(R_i)_{i\geq 1}$ from the proof of Lemma \ref{l:uniqueness_Psi}.
  Putting \eqref{eq:lange_rechnung_1}
  and~\eqref{eq:lange_rechnung_2} together, we arrive at
  \begin{equation}  \begin{split}
    \lefteqn{J_\infty(\ld)=\limn \abs{J(\ld,2t_n)}
    \leq \int_0^{\infty}
       J_{\infty}\rub{\ld e^{-\al s}}
       \mu_\al^Q(ds)}\\
    &=\E J_\infty\rub{\ld e^{-\al R_1}}
    \leq\ldots\leq
    \E J_\infty\rub{\ld e^{-\al (R_1+\ldots+R_m)}}\lram J_\infty(0+)=0.
  \end  {split}     \end  {equation}
  This finishes the proof.
\end  {proof}
%
%
%
If the $(x \log x)$-condition fails to hold, then the rescaled
supercritical Virgin Island Model converges to zero. The proof of this
assertion follows Kaplan~\cite{Kap75}.
\begin{lemma}  \label{l:convergence_if_xlogx_fails}
  Assume~\ref{a:general:finite_moments},
  \ref{a:general:finite_man_hours},
  \ref{a:convergence_to_zero} and $\al>0$.
  If the $(x \log x)$-condition \eqref{eq:xlogx} fails to hold, then
  \begin{equation}  \label{eq:convergence_if_xlogx_fails}
     L_t\ruB{\frac{\ld}{\bm e^{\al t}}}\lra 0\qqast
  \end  {equation}
  for every $\ld\geq0$
\end  {lemma}
\begin{proof}
  Define $K(\ld,t):=\tfrac{1}{\ld}L_{t}(\ld e^{-\al t})$ for $\ld>0$
  and $K(0,t):=\I\rub{V_t^{\phi,Q,Q}}e^{-\al t}$.
  It suffices to prove
  \begin{equation}
    K_\infty(\ld):=\limT K_T(\ld):=\limT \sup_{t\geq T}K(\ld,t)=0.
  \end  {equation}
  Assume that $K_\infty(\ld_0)=:\dl>0$ for some $\ld_0>0$.
  We will prove that the $(x \log x)$-condition~\eqref{eq:xlogx} holds.
  An elementary calculation shows that $\ld\mapsto \tfrac{1}{\ld}(1-e^{-\ld})$
  is decreasing.
  Thus, both $K(\ld,t)$ and $K_\infty(\ld)$ are decreasing in $\ld$.
  Furthermore, by Theorem~\ref{thm:general:expected_man_hours},
  \begin{equation} \label{eq:K_bdd}
    \dl\leq K_\infty(\ld)\leq \sup_{t\geq0} K(\ld,t)
    \leq \sup_{t\geq0}\frac{\E V^{\phi,Q,Q}_t}{e^{\al t}}=:C
    <\infty\quad\fa\ld\leq\ld_0.
  \end  {equation}
  Fix $t_0>0$, $\ld\leq\ld_0$
  and let $t\geq 2t_0$.
  Inserting the recursive equation~\eqref{eq:darma},
  \begin{equation}  \begin{split}  \label{eq:estimate_for_K}
    \lefteqn{\ld K(\ld,t)=L_t(\ld e^{-\al t})}\\
    &=\int\eckB{1-\E\exp\rub{-\frac{\ld}{e^{\al t}}\phi_{\chi}(t)}
      \exp\ruB{-\int_0^t a\rub{\chi_s}L_{t-s}(\ld e^{-\al t})ds}}Q(d\chi)\\
    &\leq \sup_{u\geq t_0}\int \ruB{1-\E e^{-\ld \phi_\chi(u)}} Q(d\chi)\\
     &\qquad\qquad\,+\int\ruB{1-\exp\ruB{-\int_0^{t-t_0}
         a\rub{\chi_s}L_{t-s}(\ld e^{-\al s} e^{-\al (t-s)})ds}}Q(d\chi)\\
     &\qquad\qquad\,+\int\ruB{1-\exp\ruB{-\int_{t-t_0}^{t}
         a\rub{\chi_s}L_{t-s}(\ld e^{-\al s} e^{-\al (t-s)})ds}}Q(d\chi)\\
    &=:T_1+T_2+T_3.
  \end  {split}     \end  {equation}
  By Assumption~\ref{a:convergence_to_zero}, the first term converges
  to zero uniformly in $t\geq2t_0$
  as $t_0\to\infty$.
  For the third term, we use inequality~\eqref{eq:K_bdd} to obtain
  \begin{equation}  \begin{split}
    T_3&=\int\ruB{1-\exp\rub{-\int_{t-t_0}^t a\rub{\chi_s}
      K(\ld e^{-\al s},t-s)\ld e^{-\al s}ds}}Q(d\chi)\\
    &\leq\int \int_{t_0}^\infty a\rub{\chi_s}C \ld e^{-\al s} \,ds Q(d\chi).
  \end  {split}     \end  {equation}
  The right-hand side converges to zero uniformly in $t\geq2t_0$
  as $t_0\to\infty$ by Assumption~\ref{a:general:finite_man_hours}.
  The second term is bounded above by
  \begin{equation}  \begin{split}
     T_2
     \leq \int\ruB{1-\exp\rub{-\int_0^{\infty}a\rub{\chi_s}
       K_{t_0}(\ld e^{-\al s})\ld e^{-\al s}ds}}Q(d\chi).
  \end  {split}     \end  {equation}
  Recall $(R_i)_{i\geq 1}$ from the proof of Lemma~\ref{l:uniqueness_Psi}.
  Define $S_0=0$ and
  $S_n:=R_1+\ldots+R_n$, $n\geq1$.
  Taking supremum over $t\geq 2 t_0$ in~\eqref{eq:estimate_for_K}
  and letting $t_0\to\infty$, we arrive at
  \begin{equation}  \begin{split}  \label{eq:inequality_K}
    \lefteqn{K_\infty(\ld)\leq\frac{1}{\ld}\int\ruB{1-\exp\rub{-
      \int_0^\infty a\rub{\chi_s}K_\infty(\ld e^{-\al s})
       \ld e^{-\al s}ds}}Q(d\chi)}\\
    &=\int_0^\infty K_\infty\rub{\ld e^{-\al s}}e^{-\al s}\mu^Q(ds)
       -\frac{1}{\ld}\bar{H}_\al\ruB{\ldt\mapsto K_\infty(\ldt)\ldt}(\ld)\\
    &\leq \E\eckB{ K_\infty \rub{\ld e^{-\al R_1}}}
       -\frac{1}{\ld}\bar{H}_\al\rub{\ldt\mapsto\dl\ldt}(\ld)
    = \E\eckB{ K_\infty \rub{\ld e^{-\al R_1}}}
      -\dl\eta(\dl \ld)\\
    &\leq\cdots\leq\E\eckB{ K_\infty \rub{\ld e^{-\al S_n}}}
      -\dl\sum_{k=0}^{n-1}\E\eta\rub{\dl\ld e^{-\al S_k}}
  \end  {split}     \end  {equation}
  for all $n\geq0$.
  The second inequality follows from $\dl\leq K_\infty(\ldt)$ for
  $\ldt\leq\ld_0$  and Lemma~\ref{l:barH_nondecreasing}.
  Boundedness of $K_\infty$ on $(0,\ld_0]$, see \eqref{eq:K_bdd},
  implies
  \begin{equation}
   \sum_{k=0}^\infty \eta\rub{\dl\ld e^{-\al S_k}}<\infty\qquad\text{a.s.}.
  \end  {equation}
  By the law of large numbers, we know that $S_k\leq k(\E R_1+\eps)$ for
  large $k$ a.s.
  Hence,
  \begin{equation}
    \sum_{k=0}^\infty\eta\rub{\dl\ld r^k}<\infty
  \end  {equation}
  where $r=e^{-\al (\E R_1+\eps)}\in(0,1)$.
  Therefore, 
  the $(x \log x)$-condition~\eqref{eq:xlogx} holds
  by Lemma~\ref{l:properties_of_eta}.
  This finishes the proof.
\end  {proof}
%
\begin{proof}[\textbf{\upshape Proof of Theorem~\ref{thm:general:xlogx}}]
  Assume that the $(x \log x)$-condition~\eqref{eq:xlogx} holds.
  Insert \eqref{eq:convergence_xlogx} into~\eqref{eq:integral_laplace_trafo}
  and use
  Assumption~\ref{a:convergence_to_zero} to obtain
  \begin{equation}  \label{eq:def:Psit}
    \E\eckbb{\exp\ruB{-\frac{\ld V^{\phi,\nu,Q}_t}{\mb e^{\al t}}}}
    \lrat\int\eckbb{\exp\ruB{-\int_0^\infty \Psi(\ld e^{-\al s}) a(\chi_s)\,ds}}
       \nu(d\chi)
  \end  {equation}
  for $\ld\geq0$. For this, we applied the dominated convergence theorem
  together with Assumption~\ref{a:general:finite_man_hours}.
  Denote the right-hand side of~\eqref{eq:def:Psit} by $\Psit(\ld)$
  and note that $\Psit$ is continuous and satisfies $\Psit(0+)=1$.
  A standard result, e.g.\ Lemma 2.1 in~\cite{Dyn89},
  provides us with the existence of a random variable $W\geq0$ such that
  $\E e^{-\ld W}=\Psit(\ld)$ for all $\ld\geq0$.
  This proves the weak
  convergence~\eqref{eq:weak_convergence_rescaled_VIM}
  as the Laplace transform is convergence determining.
  Note that
  \begin{equation}
    \P(W=0)=\Psit(\infty)
    =\int \eckbb{\exp\ruB{-\int_0^\infty \Psi(\infty)a(\chi_s)\,ds}}\nu(d\chi)
  \end  {equation}
  by the dominated convergence theorem. Furthermore,
  \begin{equation}
    \E W=\limldO \frac{1-\Psit(\ld)}{\ld}
    =\int \eckB{\int_0^\infty e^{-\al s}a(\chi_s)ds}\,\nu(d\chi).
  \end  {equation}
  If the $(x \log x)$-condition fails to hold,
  then $\E\eckb{1-\exp\rub{-\ld V_t^{\phi,\nu,Q}/e^{\al t}}}\to 0$
  as $t\to\infty$ follows by inserting~\eqref{eq:convergence_if_xlogx_fails}
  into~\eqref{eq:integral_laplace_trafo} together
  with~\ref{a:convergence_to_zero}.
\end  {proof}

%
\label{sec:vor_excursions}
\noindent
\section{Excursions from a trap of one-dimensional diffusions.\\ Proof of Theorem~\ref{thm:existence_excursion_measure}}%
\label{sec:excursions_from_a_trap_of_one_dimensional_diffusions}
\newcommand{\Ypf}{\ensuremath{Y^{\uparrow}}}%
Recall the Assumptions~\ref{a:A1}, \ref{a:A2}, \ref{a:S_bar},
\ref{a:finite_man_hours} and~\ref{a:finite_man_hours_squared}
from Section~\ref{sec:main_results}.
The process $(Y_t)_{t\geq0}$,
the excursion set $U$
and
the scale function $\bS$
have been defined in~\eqref{eq:def:Y},
in~\eqref{eq:def:U}
and
in~\eqref{eq:def:bar_s}, respectively.
The stopping time $T_\eps$ has been introduced shortly after~\eqref{eq:def:U}.

In this section, we define the excursion measure $\barQY$ 
and prove the convergence result of
Theorem~\ref{thm:existence_excursion_measure}.
We follow Pitman and Yor~\cite{PY82} in the construction of the
excursion measure.
Under Assumptions~\ref{a:A1} and~\ref{a:A2},
zero is an absorbing point for $Y$.
Thus, we cannot simply start
in zero and wait until the process returns to zero.
Informally speaking, we instead condition the process to converge
to infinity.
One way to achieve this is by Doob's h-transformation.
Note that $\rub{\bS(Y_{t\wedge T_\eps})}_{t\geq0}$
is a bounded martingale
for every $\eps>0$, 
see Section V.28 in~\cite{RW2}.
In particular,
\begin{equation}
  \E^y \eckb{\bS(Y_{t\wedge T_\eps})}=\bS(y)
\end  {equation}
for every $y< \eps$.
For $\eps>0$, consider the diffusion 
$\ru{Y^{\upa,\eps}_t}_{t\geq0}$ on $[0,\infty)$
-- to be called the $\upa$-diffusion stopped at time $T_\eps$ --
defined by the
semigroup $(T_t^\eps)_{t\geq0}$ where
\begin{equation}  \label{eq:def:semigroup_upa}
  T_t^\eps f(y):=\frac{1}{\bS(y)}\E^y\eckb{\bS(Y_{t\wedge T_\eps})
     f(Y_{t\wedge T_\eps})},\quad y>0,t\geq0, f\in\C_b\rub{[0,\infty),\R}.
\end  {equation}
The sequence of processes $\rub{(Y_t^{\upa,\eps})_{t\geq0},\eps>0}$
is consistent in the sense that
\begin{equation}
  \Law[y]{Y^{\upa,\eps+\dl}_{\centerdot\wedge T_{\eps}}}
  =\Law[y]{Y^{\upa,\eps    }_{\centerdot}}
\end  {equation}
for all $0\leq y\leq\eps$ and $\dl>0$. Therefore, we may define a process
$Y^\upa=(Y_t^\upa)_{0\leq t\leq T_\infty}$
which coincides with $(Y_t^{\upa,\eps})_{t\geq0}$ until time $T_\eps$
for every $\eps>0$. Note that the $\upa$-diffusion possibly
explodes in finite time.

The following important observation of Williams
has been quoted by Pitman and Yor~\cite{PY82}.
Because we assume that zero is an exit boundary for $(Y_t)_{t\geq0}$,
zero is an entrance boundary but not an exit boundary for the
$\upa$-diffusion.
More precisely, the $\upa$-diffusion started at its entrance boundary zero
and run up to the last time it hits a level $y>0$ is described by
Theorem 2.5 of Williams~\cite{Wil74}
as the time reversal back from 
$T_0$ of the $\downa$-diffusion started at $y$, where the
$\downa$-diffusion is the process $(Y_t)_{t\geq0}$ conditioned on $T_0<\infty$.
Hence, the process $\rub{\Ypf_t}_{t\geq0}$ may be started in zero but
takes strictly positive values at positive times.

Pitman and Yor~\cite{PY82} define the excursion measure $\QY$ as
follows. Under 
\begin{equation}
  \QY(\centerdot|T_\eps<T_0),
\end  {equation}
that is, conditional on
``excursions reach level $\eps$'', an excursion follows the
$\upa$-diffusion until time $T_\eps$ and then follows the 
dynamics of $(Y_t)_{t\geq0}$.
In addition, $\barQY\rub{T_\eps<T_0}=\tfrac{1}{\Sb(\eps)}$.
With this in mind, define a process
$\Yh^{\eps}:=\rub{\Yh^\eps_t}_{t\geq0}$
which satisfies
\begin{eqnarray}
  \Lawb[y]{\ru{\Yh^{\eps}_{t\wedge T_\eps}}_{t\geq0}}
  \!\!&=&\!\! \Lawb[y]{\ru{Y^{\upa,\eps}_{t}}_{t\geq0}}\label{eq:BED1}\\
  \Lawb[y]{\ru{\Yh^{\eps}_{T_\eps+t}}_{t\geq0}}
  \!\!&=&\!\! \Lawb[\eps]{\ru{Y_{t}}_{t\geq0}}\label{eq:BED2}
\end  {eqnarray}
for $y\geq0$.
In addition,  $(\Yh^{\eps}_t,t\leq T_\eps)$
and $(\Yh^{\eps}_t,t\geq T_\eps)$ are independent.
Define the excursion measure $\QY$ on $U$ by
\begin{equation}
  \1_{T_\eps<T_0}\QY(d\chi):=
  \frac{1}{\bS(\eps)} \P^0\rub{\Yh^\eps\in d\chi},\quad \eps>0.
\end  {equation}
This is well-defined if
\begin{equation}  \label{eq:for_well_defined}
  \1_{T_{\eps+\dl}<T_0}
  \frac{1}{\bS(\eps)} \P^0\rub{\Yh^\eps\in d\chi}
  =
  \frac{1}{\bS(\eps+\dl)} \P^0\rub{\Yh^{\eps+\dl}\in d\chi}
\end  {equation}
holds for all $\eps,\dl>0$.
The critical part here is the path between $T_\eps$ and $T_{\eps+\dl}$.
Therefore, \eqref{eq:for_well_defined} follows from
\begin{equation}  \begin{split}
  \lefteqn{\frac{1}{\bS(\eps)}
  \E^\eps\eckb{F(Y)\1_{T_{\eps+\dl}<T_0}}
  =\frac{1}{\bS(\eps+\dl)}
  \E^\eps\eckb{F(Y)|{T_{\eps+\dl}<T_0}}}\\
  &=\frac{1}{\bS(\eps+\dl)}
  \E^\eps\eckb{F(\Yh^{\eps+\dl})}
  =\frac{1}{\bS(\eps+\dl)}
  \E^0\eckb{F(\Yh^{\eps+\dl}_{T_\eps+\centerdot})}.
\end  {split}     \end  {equation}
The first equality is equation~\eqref{eq:def:scale_function}
with $c=0$, $y=\eps$ and $b=\eps+\dl$.
The last equality is the strong Markov property of $Y^{\upa,\eps+\dl}$.
The last but one equality is the following lemma.

%
%
%
\begin{lemma}  \label{l:conditioned_on_T_eps}
  Assume~\ref{a:A1} and~\ref{a:A2}.
  Let $0<y<\eps$. Then
  \begin{equation} \label{eq:conditioned_on_T_eps}
    \Lawb[y]{\,Y\,|\,T_\eps<T_0}=\Lawb[y]{\Yh^{\eps}}.
  \end  {equation}
\end  {lemma}
\begin{proof}
  We begin with the proof of independence of
  $(\Yh^{\eps}_t,t\leq T_\eps)$ and of
  $(\Yh^{\eps}_t,t\geq T_\eps)$.
  Let $F$ and $G$
  be two bounded continuous functions on the path space.
  Denote by $\MCF_{T_\eps}$ the $\sigma$-algebra generated
  by $(Y_t)_{t\leq T_\eps}$.
  Then
  \begin{equation}  \begin{split} \label{eq:indie}
    \lefteqn{\E^y\eckb{F\rub{Y_{T_\eps\wedge\centerdot}}
               G\rub{Y_{T_\eps+\centerdot}}|T_\eps<T_0}}\\
    &= \E^y\eckB{F\rub{Y_{T_\eps\wedge\centerdot}}
      \E^y\eckb{G\rub{Y_{T_\eps+\centerdot}}|\MCF_{T_\eps}}|T_\eps<T_0}\\
    &= \E^y\eckb{F\rub{Y_{T_\eps\wedge\centerdot}} |T_\eps<T_0}
      \E^\eps\eckb{G\rub{Y_{\centerdot}}}.
  \end  {split}     \end  {equation}
  The last equality is the strong Markov property of $Y$.
  Choosing $F\equiv 1$ in~\eqref{eq:indie}
  proves that the left-hand side
  of~\eqref{eq:conditioned_on_T_eps} satisfies~\eqref{eq:BED2}.
  In addition, equation~\eqref{eq:indie} proves the desired
  independence.
  For the proof of
  \begin{equation} \label{eq:for_theProof_of}
    \P^y\rub{\ru{Y^{\upa,\eps}_{t}}_{t\geq0}}=
    \P^y\rub{(Y_{t\wedge T_\eps})_{t\geq0}|T_\eps<T_0},
  \end  {equation}
  we repeatedly apply the semigroup~\eqref{eq:def:semigroup_upa}
  of $(Y_t^{\upa,\eps})_{t\geq0}$ to obtain
  \begin{equation}  \label{eq:prodigod}
    \E^y\eckB{\prod_{i=1}^nf_i\rub{Y^{\upa,\eps}_{t_i}}}
    =\frac{1}{\bS(y)}\E^y\eckB{
     \bS(Y_{t_n\wedge T_\eps})
     \prod_{i=1}^nf_i\rub{Y_{t_i\wedge T_\eps}}}
  \end  {equation}
  for bounded, continuous functions $f_1,...,f_n$ and time points
  $0\leq t_1<...<t_n$.
  By equation~\eqref{eq:def:scale_function} with $c=0$,
  \begin{equation}
    \bS(Y_{t_n\wedge T_\eps})=
    \bS(\eps)\P^{Y_{t_n\wedge T_\eps}}\eckb{T_\eps<T_0}
    =\bS(\eps)
    \E^{y}\eckb{\1_{T_\eps<T_0}|\MCF_{t_n\wedge T_\eps}}
  \end  {equation}
  $\P^y$--almost surely where $\MCF_{t_n\wedge T_\eps}$ is the $\sigma$-algebra generated
  by $(Y_s)_{s\leq t_n\wedge T_\eps}$.
  Insert this identity in the right-hand side
  of~\eqref{eq:prodigod} to obtain
  \begin{equation}
    \E^y\eckB{\prod_{i=1}^nf_i\rub{Y^{\upa,\eps}_{t_i}}}
    =\frac{1}{\P^y\rub{T_\eps<T_0}}\E^y\eckB{
     \1_{T_\eps<T_0}
     \prod_{i=1}^nf_i\rub{Y_{t_i\wedge T_\eps}}}.
  \end  {equation}
  This proves~\eqref{eq:for_theProof_of}
  because finite-dimensional distributions determine the law of a process.
\end  {proof}

Now we prove convergence to the excursion measure $\barQY$.
\begin{proof}[\textbf{\upshape Proof of Theorem~\ref{thm:existence_excursion_measure}}]
  Let $F\colon\mathbf{C}\rub{[0,\infty),[0,\infty)}\to\R$
  be a bounded continuous function for which
  there exists an $\eps>0$ such
  that $F(\chi)\1_{T_0<T_\eps}=0$ for every path $\chi$.
  Let $0<y<\eps$.
  By Lemma~\ref{l:conditioned_on_T_eps}, we obtain
  \begin{equation}  \begin{split} \label{eq:selters1}
    \frac{1}{\bS(y)}\E^yF(Y)
    &=\frac{1}{\bS(\eps)\P^y(T_\eps<T_0)}
      \E^y \eckb{F(Y)\1_{T_\eps<T_0}}\\
    &=\frac{1}{\bS(\eps)} \E^y F(\Yh^{\eps})
    =\frac{1}{\bS(\eps)}\E^0 F(\Yh^{\eps}_{T_y+\centerdot}).
  \end  {split}     \end  {equation}
  The last equality
  is the strong Markov property of the $\upa$-diffusion.
  The random time $T_y$ converges to zero almost surely as $y\to0$.
  Another observation we need is that every continuous
  path $(\chi_t)_{t\geq0}$ is
  uniformly continuous on any compact set $[0,T]$.
  Hence, the sequence of paths $\rub{\ru{\chi_{T_y+t}}_{t\geq 0},y>0}$
  converges locally uniformly
  to the path $\rub{\chi_t}_{t\geq0}$ almost surely as $y\to0$.
  Therefore, the dominated convergence theorem implies
  \begin{equation}  \begin{split} \label{eq:selters2}
    \lim_{y\to0}
    \E^0 F(\Yh^{\eps}_{T_y+\centerdot})
    =\E^0 \lim_{y\to0} F(\Yh^{\eps}_{T_y+\centerdot})
    =\E^0 F(\Yh^{\eps}_\centerdot).
  \end  {split}     \end  {equation}
  Putting~\eqref{eq:selters1} and~\eqref{eq:selters2} together,
  we arrive at
  \begin{equation}  \begin{split}
    \lim_{y\to0} \frac{1}{\bS(y)}\E^yF(Y)
    &=\frac{1}{\bS(\eps)} \E^0 F(\Yh^{\eps})
    =\int F(\chi)\QY(d\chi),
  \end  {split}     \end  {equation}
  which proves the theorem.
\end  {proof}

We will employ Lemma~\ref{l:conditioned_on_T_eps} to calculate
explicit expressions for some functionals of $\barQY$.
For example, we will prove in Lemma~\ref{l:Q_explicit_man_hours}
together with Lemma~\ref{l:finite_man_hours} that
\begin{equation} \label{eq:zwischentext}
  \int\ruB{\int_0^\infty a\rub{\chi_s}\,ds}\QY(d\chi)
  =\int_0^\infty \frac{a(z)}{g(z)\bs(z)}\,dz
\end  {equation}
provided that Assumptions~\ref{a:A1}, \ref{a:A2}
and~\ref{a:finite_man_hours} hold.
Equation~\eqref{eq:zwischentext} shows that
condition~\eqref{eq:general:bed_fuer_extinction} and
condition~\eqref{eq:bed_fuer_extinction} are equivalent.
The following lemmas prepare for the proof of~\eqref{eq:zwischentext}.

%
%
%
\begin{lemma}  \label{l:mit_T_b}
  Assume~\ref{a:A1} and~\ref{a:A2}.
  Let $f\in\C\rub{[0,\infty),[0,\infty)}$
  have compact support in $(0,\infty)$.
  Furthermore, let the continuous function $\psi\colon[0,\infty)\to\R$
  be nonnegative and nondecreasing.
  Then
  \begin{equation} \label{eq:kalkhofe}
    \frac{1}{\bS(y)}\E^y\eckbb{\ruB{\int_0^{T_b}\psi(s)f(Y_s)\,ds}^m}
    \lrayO\int\eckbb{\ruB{
       \int_0^{T_b}\psi(s)f(\chi_s)\,ds}^m}\QY(d\chi)
  \end  {equation}
  for every $b\leq\infty$ and $m\in\N_{\geq0}$.
\end  {lemma}
\begin{proof}
  W.l.o.g.\ assume $m\geq1$.
  Let $\eps>0$ be such that $\eps<\inf\supp f$ and let $y<\eps$.
  Using Lemma~\ref{l:conditioned_on_T_eps}, we see that the
  left-hand side of~\eqref{eq:kalkhofe} is equal to
  \begin{equation*}  \begin{split}
    &\frac{1}{\bS(y)}\E^y\eckbb{\ruB{\int_0^{T_b}
        \psi(s)f(Y_s)\,ds}^m\1_{T_\eps<T_0}}
    =\frac{1}{\bS(\eps)}\E^y\eckbb{\ruB{\int_0^{T_b}\psi(s)
      f(\Yh^{\eps}_s)\,ds}^m}\\
    &=\frac{1}{\bS(\eps)}\E^0\eckbb{\ruB{\int_{T_y}^{T_b}\psi(s-T_y)
        f(\Yh^{\eps}_{s})\,ds}^m}
    \lrayO \int\rubb{\int_0^{T_b}\psi(s)
        f(\chi_s)\,ds}^m\QY(d\chi).
  \end  {split}     \end  {equation*}
  The second equality is the strong Markov property of $Y^{\upa,\eps}$
  and the change of variable $s\mapsto s-T_y$.
  For the convergence, we applied the monotone convergence theorem.
\end  {proof}

The explicit formula on the right-hand side of~\eqref{eq:zwischentext}
originates in the explicit formula~\eqref{eq:explicit_man_hours} below,
which we recall from the literature.
%
%
\begin{lemma}  \label{l:explicit_man_hours}
  Assume~\ref{a:A1} and~\ref{a:A2}.
  If $f\in \mathbf{C}_b[0,\infty)$ or
  $f\in\mathbf{C}\rub{[0,\infty),[0,\infty)}$,
  then
  \begin{equation}   \label{eq:explicit_man_hours}
    \E^y\ruB{\int_0^{T_0\wedge T_b}f(Y_s)\,ds}
    =\int_0^{b}\ruB{f(z)\frac{\bS(b)-\bS(y\vee z)}{\bS(b)}
      \frac{\bS(y\wedge z)}{g(z)\bs(z)}} \,dz
  \end  {equation}
  for all $0\leq y\leq b<\infty$.
\end  {lemma}
\begin{proof}
  See e.g.\ Section 15.3 of Karlin and Taylor~\cite{KT2}.
\end  {proof}

Let $(\Yt_t)_{t\geq0}$ be a Markov process with \cadlag\ sample paths
and state space $E$ which is a Polish space.
For an open set $O\subset E$, denote by $\tau$ the first exit time
of $(\Yt_t)_{t\geq0}$
from the set $O$. Notice that $\tau$ is a stopping time.
For $m\in\N_0$, define
\begin{equation}  \label{eq:def:wmn}
  w_m(y):=\E^y\eckbb{\ruB{{\int_0^{\tau}f(\Yt_s)\,ds}}^m},\ y\in E,m\in\N_0,
\end  {equation}
for a given function $f\in\mathbf{C}\rub{O,[0,\infty)}$.
In the following lemma, we derive an expression for $w_2$
for which Lemma~\ref{l:explicit_man_hours} is applicable.
%
\begin{lemma}  \label{l:conversion}
  Let $(\Yt_t)_{t\geq0}$ be a time-homogeneous
  Markov process with \cadlag\ sample paths
  and state space $E$ which is a Polish space.
  Let $w_m$ be as in~\eqref{eq:def:wmn} with an open set $O\subset E$
  and with a function $f\in\mathbf{C}\rub{O,[0,\infty)}$.
  Then
  \begin{eqnarray}
    \E^y\ruB{{\int_0^{\tau}sf(\Yt_s)\,ds}}
      \!\!\!\!&=&\!\!\!\!\E^y\ruB{\int_0^{\tau} w_1(\Yt_s)\,ds}
        \label{eq:wEE}\\
    \E^y\eckbb{\ruB{{\int_0^{\tau}f(\Yt_s)\,ds}}^2}
      \!\!\!\!&=&\!\!\!\!\E^y\ruB{\int_0^{\tau}2 f(\Yt_s)w_1(\Yt_s)\,ds}
        \label{eq:wZZ}
  \end  {eqnarray}
  for all $y\in E$.
\end  {lemma}
\begin{proof}
  Let $y\in E$ be fixed.
  For the proof of~\eqref{eq:wEE},
  we apply Fubini to obtain
  \begin{equation}  \begin{split}  \label{eq:schorsch}
    &\E^y\ruB{\int_0^{\tau}\int_0^sdr f(\Yt_s)\,ds}
    =\E^y\ruB{\int_0^{\tau}\int_r^{\tau}f(\Yt_s)\,ds\,dr}\\
    &=\int_0^{\infty}\E^y\ruB{\1_{r<\tau}
       \int_0^{\infty} \1_{s+r<\tau}f(\Yt_{s+r})ds}\,dr.
  \end  {split}     \end  {equation}
  The last equality follows from Fubini and a change of variables.
  The stopping time $\tau$ can be expressed
  as $\tau=F\rub{\ru{\Yt_u}_{u\geq0}}$
  with a suitable path functional $F$.
  Furthermore, $\tau$ satisfies
  \begin{equation}
    \{r<\tau\}\cap\{s+r<\tau\}
    =\{r<\tau\}\cap\{s< F\rub{\ru{\Yt_{u+r}}_{u\geq0}}\}
  \end  {equation}
  for $r,s\geq0$.
  Therefore, the right-hand side of~\eqref{eq:schorsch} is equal to
  \begin{equation}  \begin{split}
    \lefteqn{\int_0^{\infty}\E^y\ruB{\1_{r<\tau}
     \int_0^{\infty}\1_{s<F\rub{\ru{\Yt_{u+r}}_{u\geq0}}}f(\Yt_{s+r})\,ds}\,dr}\\
    &=\int_0^{\infty}\E^y\ruB{\1_{r<\tau}\E^{\Yt_r}\eckb{
     \int_0^{\infty}\1_{s<\tau}f(\Yt_{s})\,ds}}\,dr
    =\E^y\ruB{\int_0^{\tau}w_1(\Yt_r)\,dr}.
  \end  {split}     \end  {equation}
  The last but one equality is the Markov property of $(\Yt_t)_{t\geq0}$.
  This proves~\eqref{eq:wEE}.
  For the proof of~\eqref{eq:wZZ}, break the symmetry in the
  square of $w_2(y)$ to see that $w_2(y)$ is equal to
  \begin{equation}  \begin{split}
    \lefteqn{\E^y\ruB{2\int_0^{\tau}\eckB{f(\Yt_r)\int_r^{\tau}f(\Yt_s)\,ds}dr}}\\
    &=2\int_0^{\infty}\E^y\ruB{\1_{r<\tau}f(\Yt_r)
       \E^{\Yt_r}\eckb{\int_0^{\tau}f(\Yt_{s})\,ds}}\,dr
    =\E^y\ruB{\int_0^{\tau}2f(\Yt_s)w_1(\Yt_s)\,ds}.
  \end  {split}     \end  {equation}
  This finishes the proof.
\end  {proof}

We will need that $(Y_t)_{t\geq0}$ dies out in finite time. The following
lemma gives a condition for this. Recall $\bS(\infty):=\lim_{y\to\infty}\bS(y)$.
%
%
\begin{lemma}   \label{l:finite_time_extinction}
  Assume~\ref{a:A1} and \ref{a:A2}.
  Let $y>0$.
  Then
  the solution $(Y_t)_{t\geq0}$ of
  equation~\eqref{eq:def:Y} hits zero in finite time almost surely
  if and only if
  $\bS(\infty)=\infty$.
  If $\bS(\infty)<\infty$, then $(Y_t)_{t\geq0}$ converges to infinity
  as $t\to\infty$ on the event $\{T_0=\infty\}$ almost surely.
\end  {lemma}
\begin{proof}
  On the event $\{Y_t\leq K\}$, we have that
  \begin{equation}
    \P^{Y_t}\rub{\exists s\colon Y_s=0}\geq\P^K\rub{T_0<\infty}>0
  \end  {equation}
  almost surely.
  The last inequality follows from
  Lemma 15.6.2 of~\cite{KT2}
  and
  Assumption~\ref{a:A2}.
  Therefore, Theorem 2 of Jagers~\cite{Ja92} implies that,
  with probability one, either $(Y_t)_{t\geq0}$ hits
  zero in finite time or converges to infinity as $t\to\infty$.
  With equation~\eqref{eq:def:scale_function},
  we obtain
  \begin{equation}  \label{eqdickidom}
    \P^y\rub{\lim_{t\to\infty} Y_t=\infty}
    =\lim_{b\to\infty}\P^y\rub{Y\text{ hits }b\text{ before }0}
    = \lim_{b\to\infty}\frac{\bS(y)}{\bS(b)}=\frac{\bS(y)}{\bS(\infty)}.
  \end  {equation}
  This proves the assertion.
\end  {proof}

The following lemma 
makes Assumption~\ref{a:finite_man_hours} more transparent.
It proves that \ref{a:finite_man_hours} holds if
and only if the expected area under $\rub{a(Y_t)}_{t\geq0}$ is finite.
%
%
%
\begin{lemma}   \label{l:finite_man_hours}
  Assume~\ref{a:A1} and \ref{a:A2}.
  Assumption \ref{a:finite_man_hours} holds if and only
  if
  \begin{equation}
    \E^y\ruB{\int_0^\infty a(Y_s)\,ds}<\infty\qquad\fa y>0.
  \end  {equation}
  If Assumption~\ref{a:finite_man_hours} holds, then
  $\bS(\infty)=\infty$  and
  \begin{equation}  \label{eq:finite_man_hours_Y}
    \E^y\ruB{\int_0^\infty f\rub{Y_s}\,ds}
    =\int_0^\infty \bS\rub{y\wedge z}\frac{f\ru{z}}{g(z)\bs(z)}\,dz<\infty
  \end  {equation}
  for all $y\geq0$ and $f\in\C\rub{[0,\infty),[0,\infty)}$
  with $c_f:=\sup_{z>0}f(z)/z<\infty$.
\end  {lemma}
\begin{proof}
  Let $c_1,c_2$ be the constants from~\ref{a:A1}.
  In equation~\eqref{eq:explicit_man_hours},
  let $b\to\infty$ and apply monotone convergence to obtain
  \begin{equation}  \begin{split} \label{eq:bleistift}
    \E^y\ruB{\int_0^{\infty}f(Y_s)\,ds}
    &=\int_0^{\infty}\ruB{f(z)\eckB{1-\frac{\bS(y\vee z)}{\bS(\infty)}}
      \frac{\bS(y\wedge z)}{g(z)\bs(z)}} \,dz.
  \end  {split}     \end  {equation}
  Hence, if Assumption~\ref{a:finite_man_hours} holds,
  then Assumption~\ref{a:A2}
  implies that the  right-hand side
  of~\eqref{eq:bleistift} is finite because
  $f(z)\leq c_f z\leq\tfrac{c_f}{c_1}a(z)$, $z>0$.
  Therefore, the left-hand side of~\eqref{eq:bleistift} with $f(\cdot)$
  replaced by $a(\cdot)$ is finite.
  Together with $\lim_{x\to\infty}a(x)=\infty$, this implies that
  $(Y_t)_{t\geq0}$ does not converge to infinity with
  positive probability as $t\to\infty$.
  Thus Lemma~\ref{l:finite_time_extinction} implies $\bS(\infty)=\infty$
  and equation~\eqref{eq:bleistift}
  implies~\eqref{eq:finite_man_hours_Y}.

  Now we prove that Assumption~\ref{a:finite_man_hours} holds
  if the left-hand side of~\eqref{eq:bleistift} with $f(\cdot)$
  replaced by $a(\cdot)$ is finite.
  Again, $\lim_{x\to\infty}a(x)=\infty$
  and Lemma~\ref{l:finite_time_extinction} imply $\bS(\infty)=\infty$.
  Using monotonicity of $S$, we obtain for $x>0$
  \begin{equation}  \begin{split}
    \int_x^\infty \frac{a(z)}{g(z)\bs(z)}\,dz
    \leq \frac{1}{\bS(x)}\int_0^\infty a(z)\frac{\bS(x\wedge z)}{g(z)\bs(z)}\,dz.
  \end  {split}     \end  {equation}
  The right-hand side is finite because~\eqref{eq:bleistift} 
  with $f(\cdot)$ replaced by $a(\cdot)$ is finite.
  Therefore, Assumption~\ref{a:finite_man_hours} holds.
\end  {proof}

%
%
\begin{lemma}  \label{l:bS_durch_y_bounded}
  Assume~\ref{a:A1}, \ref{a:S_bar} and let $n\in\N_{\geq1}$.
  If $\int_1^\infty \tfrac{y^n}{g(y)\sb(y)}\,dy<\infty$, then
  \begin{equation}
    \sup_{y\in(0,\infty)}\frac{y^n}{\bS(y)}<\infty.
  \end  {equation}
\end  {lemma}
\begin{proof}
  It suffices to prove $\liminf_{y\to\infty}\tfrac{\bS(y)}{y^n}>0$
  because $\tfrac{y^n}{\bS(y)}$ is locally bounded in $(0,\infty)$
  and $\bS^{'}(0)\in(0,\infty)$ by Assumption~\ref{a:S_bar}.
  By Assumption~\ref{a:A1}, $g(y)\leq c_g y^2$ for all $y\geq 1$
  and a constant $c_g<\infty$.
  Let $0\leq x\mapsto \psi(x):=1-(1-x)^+\wedge 1$.
  Thus,
  \begin{equation} \label{eq:estimate_c_g}
    \infty>\int_1^\infty \frac{y^n}{g(y)\bs(y)}\,dy
    \geq\frac{1}{c_g}\int_1^\infty\frac{y^{n-1}}{y\bs(y)}\,dy
    \geq\frac{1}{c_g}\int_1^\infty\frac{1}{y}\mal
          \ruB{1-\psi\rub{\frac{\bs(y)}{y^{n-1}}}}\,dy.
  \end  {equation}
  The last inequality follows from
  $\tfrac{1}{z}\geq \1_{z\leq 1}\geq 1-\psi(z)$, $z>0$.
  Consequently,
  \begin{equation} \label{eq:one_equal_to_psi}
    1=\limz\frac{\int_1^z\frac{1}{y}\psi\rub{\frac{\sb(y)}{y^{n-1}}}\,dy}{\log(z)}
    =\limz\frac{\frac1z\psi\rub{\frac{\bs(z)}{z^{n-1}}}}{\frac1z}
    =\limz\psi\rub{\frac{\bs(z)}{z^{n-1}}}.
  \end  {equation}
  The proof of the second equation in~\eqref{eq:one_equal_to_psi}
  is similar to the proof of the lemma
  of L'Hospital. From~\eqref{eq:one_equal_to_psi}, we conclude
  $\liminf_{y\to\infty}\tfrac{\bs(y)}{y^{n-1}}\geq1$ which implies
  \begin{equation}
       \liminf_{z\to\infty}\frac{\int_0^z \bs(y) \,dy}{z^n}
   \geq\liminf_{z\to\infty}\frac{\int_0^z y^{n-1}\,dy}{z^n}
    =\frac{1}{n}.
  \end  {equation}
  This finishes the proof.
\end  {proof}

Now we prove
equation~\eqref{eq:zwischentext}.
Recall $\bS(\infty):=\lim_{y\to\infty}\bS(y)$.
Define $\wb_0\equiv1$ and
\begin{equation}
  \wb_1(z):=\int_0^\infty f(u)\frac{\bS(z\wedge u)}{g(u)\bs(u)}\,du,\quad z\geq0
\end  {equation}
for $f\in\mathbf{C}\rub{[0,\infty),[0,\infty)}$.
If $\bS(\infty)=\infty$, then $\wb_1(z)$ is the monotone
limit of the right-hand side
of~\eqref{eq:explicit_man_hours} as $b\to\infty$.
%
\begin{lemma}  \label{l:Q_explicit_man_hours}
  Assume \ref{a:A1}, \ref{a:A2}
  and $\bS(\infty)=\infty$.
  Let $f\in\mathbf{C}\rub{[0,\infty),[0,\infty)}$.
  Then
  \begin{eqnarray}
    \int\ruB{\int_0^\infty f(\chi_s)\,ds}^m\QY(d\chi)
      \!\!\!\!&=&\!\!\!\!
        \int_0^\infty f(z)\frac{m \wb_{m-1}(z)}{g(z)\bs(z)}\,dz
         \label{eq:Qfunctional}\\
    \int\ruB{\int_0^\infty sf(\chi_s)\,ds}\QY(d\chi)
      \!\!\!\!&=&\!\!\!\!
        \int_0^\infty \wb_1(z)\frac{1}{g(z)\bs(z)}\,dz\label{eq:QwOI}
  \end  {eqnarray}
  for $m=1,2$.
  If~\ref{a:finite_man_hours} holds and if $f(z)/z$ is bounded,
  then~\eqref{eq:Qfunctional} is finite for $m=1$.
  If~\ref{a:finite_man_hours_squared} holds and if $f(z)/z$ is
  bounded, then~\eqref{eq:Qfunctional} is finite for $m=2$.
\end  {lemma}
\begin{proof}
  Choose $f_\eps\in\mathbf{C}\rub{[0,\infty),[0,\infty)}$
  with compact support in $(0,\infty)$ for every $\eps>0$ such that
  $f_\eps\upa f$ as $\eps\to0$.
  Fix $\eps>0$ and $b\in(0,\infty)$.
  Lemma~\ref{l:mit_T_b} proves that
  \begin{equation}  \begin{split}  \label{eq:loc_limit}
    \limyO \frac{1}{\bS(y)}\E^y\eckbb{\ruB{ \int_0^{T_b}f_{\eps}(Y_s)\,ds}^m}
    =\int\ruB{\int_0^{T_b} f_\eps\ru{\chi_s}\,ds}^m \QY\ru{d\chi}.
  \end  {split}     \end  {equation}
  Let $w_m^b(y)$ be as in~\eqref{eq:def:wmn} with $\tau$ replaced by $T_b$
  and $f$ replaced by $f_\eps$.
  Fix $m\in\{1,2\}$.
  Lemma~\ref{l:conversion} and Lemma~\ref{l:explicit_man_hours} provide
  us with an expression for
  the left-hand side of equation \eqref{eq:loc_limit}. Hence,
  \begin{equation}  \begin{split} \label{eq:fast_fertig}
    \lefteqn{\int\ruB{\int_0^{T_b} f_{\eps}(\chi_s)\,ds}^m\QY(d\chi)}\\
      &=\limyO\frac{1}{\bS(y)} \int_0^{b}f_{\eps}(z)m\, w_{m-1}^b\,(z)
         \frac{\bS(b)-\bS(y\vee z)}{\bS(b)}
      \frac{\bS(y\wedge z)}{g(z)\bs(z)} \,dz\\
      &=\int_0^{b}f_{\eps}(z)\,m\, w_{m-1}^b(z)\ruB{1-\frac{\bS( z)}{\bS(b)}}
      \frac{1}{g(z)\bs(z)} \,dz.
  \end  {split}     \end  {equation}
  The last equation follows from dominated convergence and
  Assumption~\ref{a:A2}.
  Note that
  the hitting time $T_b\rub{\ru{\chi_t}_{t\geq0}}\to\infty$ as $b\to\infty$
  for every continuous path $\ru{\chi_t}_{t\geq0}$.
  By Lemma~\ref{l:explicit_man_hours} and the monotone convergence theorem,
  $w_{m-1}^b(y)\nearrow\wb_{m-1}(y)$ as $b\nearrow\infty$.
  Let $b\to\infty$, $\eps\to0$ and
  apply monotone convergence to arrive at equation~\eqref{eq:Qfunctional}.

  Similar arguments prove~\eqref{eq:QwOI}.
  Instead of~\eqref{eq:loc_limit}, consider
  \begin{equation}  \begin{split}  \label{eq:loc_limit_mit_s}
    \limyO \frac{1}{\bS(y)}\E^y\ruB{ \int_0^{T_b}sf_{\eps}(Y_s)\,ds}
    =\int\ruB{\int_0^{T_b}s f_\eps\ru{\chi_s}\,ds} \QY\ru{d\chi}
  \end  {split}     \end  {equation}
  which is implied by Lemma~\ref{l:mit_T_b}.
  Furthermore, instead of applying
  Lemma~\ref{l:explicit_man_hours} to equation~\eqref{eq:loc_limit},
  apply equation~\eqref{eq:wEE} together with
  equation~\eqref{eq:explicit_man_hours}.

  For the rest of the proof, assume that $f(z)/z$ is bounded
  by $c_f$.
  Let $c_1,c_2$ be the constants from~\ref{a:A1}.
  Note that $f(z)\leq c_f z\leq\tfrac{c_f}{c_1}a(z)$.
  Consider the right-hand side of~\eqref{eq:Qfunctional}.
  If $m=1$, then
  the integral over $[1,\infty)$ is finite by
  Assumption~\ref{a:finite_man_hours}.
  If $m=2$, then
  the integral over $[1,\infty)$ is finite by
  Assumption~\ref{a:finite_man_hours_squared}.
  The integral over $[0,1)$ is finite because of~\ref{a:A2}
  and
  \begin{equation}  \label{eq:a_leq_Sq}
    a(z)\leq c_2 z\leq \cb \Sb(z)\quad z\in[0,1],
  \end  {equation}
  where $\cb$ is a finite constant.
  The last inequality
  in~\eqref{eq:a_leq_Sq} follows from
  Lemma~\ref{l:bS_durch_y_bounded}.
\end  {proof}

%
%
%
The convergence~\eqref{eq:def:wlim} of
Theorem~\ref{thm:existence_excursion_measure}
also holds for $(\chi_s)_{s\geq0}\mapsto f(\chi_t)$,
$t$ fixed, if $f(y)/y$ is a bounded function.
For this, we first estimate the moments of $(Y_t)_{t\geq0}$.
\begin{lemma} \label{l:finite_moments_Y}
  Assume~\ref{a:A1}.
  Let $(Y_t)_{t\geq0}$ be a solution of
  equation~\eqref{eq:def:Y} and
  let $T$ be finite.
  Then, for every $n\in\N_{\geq2}$,
  there exists a constant $c_T$ such that
  \begin{equation}  \begin{split}  \label{eq:moments}
    \sup_{t\leq T}\E^y\eckb{{Y_{\tau\wedge t}}}\leq c_T{y},\quad
    \E^y\eckb{\sup_{t\leq T} {Y_t}^n}\leq c_T({y}+{y}^n)
  \end  {split}     \end  {equation}
  for all $y\geq0$ and every stopping time $\tau$.
\end  {lemma}
\begin{proof}
  The proof is fairly standard and uses
  It\^o's formula and Doob's $L_p$-inequality.
\end  {proof}

%
%
%
\begin{lemma}  \label{l:functional_fixed_t}
  Assume~\ref{a:A1}, \ref{a:A2} and~\ref{a:S_bar}.
  Let $f\colon[0,\infty)\to[0,\infty)$ be a continuous
  function
  such that $f(y)\leq c_f y\vee y^n$
  for some $n\in\N_{\geq 1}$,
  some constant $c_f<\infty$
  and
  for all $y\geq0$.
  If $\int_1^\infty \tfrac{y^n}{g(y)\sb(y)}\,dy<\infty$, then
  \begin{equation} \label{eq:fixed_time}
    \int f(\chi_t)\QY(d\chi)
    =\limyO\frac{1}{\bS(y)}\E^yf(Y_t)
    =\E^0\eckb{\frac{1}{\bS(Y^\upa_t)}f(Y^\upa_t)\1_{t<T_\infty}}
  \end  {equation}
  is bounded in $t>0$.
\end  {lemma}
\begin{proof}
  Choose $f_\eps\in\mathbf{C}\rub{[0,\infty),[0,\infty)}$
  with compact support in $(0,\infty)$ for every $\eps>0$ such that
  $f_\eps\upa f$ pointwise as $\eps\to0$.
  By Theorem~\ref{thm:existence_excursion_measure},
  \begin{equation}  \label{eq:fixed_time_eps}
    \int f_\eps(\chi_t)\QY(d\chi)
    =\limyO\frac{1}{\bS(y)}\E^yf_{\eps}(Y_t).
  \end  {equation}
  The left-hand side of~\eqref{eq:fixed_time_eps} converges 
  to the left-hand side of~\eqref{eq:fixed_time}
  as $\eps\to0$
  by the monotone convergence theorem. Hence, the first equality
  in~\eqref{eq:fixed_time} follows from~\eqref{eq:fixed_time_eps}
  if the limits $\limeps$ and $\limyO$ can be interchanged.
  For this, we prove the second equality
  in~\eqref{eq:fixed_time}.
  
  Let $b\in(0,\infty)$.
  The $\upa$-diffusion is a strong Markov process. Thus,
  by~\eqref{eq:def:semigroup_upa},
  \begin{equation}  \begin{split}  \label{eq:strong_MP}
    \lefteqn{\limyO\frac{1}{\bS(y)}\E^y \eckB{f(Y_t)\1_{t<T_b}}
    =\limyO\E^y\eckB{\frac{f(Y_t^\upa)}{\bS(Y_t^\upa)}\1_{t<T_b}}}\\
    &=\E^0\eckB{\limyO\frac{f(Y_{t+T_y}^\upa)}{\bS(Y_{t+T_y}^\upa)}
            \1_{{t+T_y}<T_b}}
    =\E^0\eckB{\frac{f(Y_t^\upa)}{\bS(Y_t^\upa)}\1_{t<T_b}}.
  \end  {split}     \end  {equation}
  The second equality follows from the dominated convergence theorem
  because of
  \begin{equation}
    \sup_{0<y\leq b}\frac{f(y)}{\bS(y)}
    \leq c_f\sup_{0<y\leq b}\frac{y\vee y^n}{\bS(y)}
    <\infty.
  \end  {equation}
  Right-continuity of the function
  $t\mapsto \tfrac{f(Y_t^\upa)}{\bS(Y_t^\upa)}\1_{t<T_b}$
  implies the last equality in~\eqref{eq:strong_MP}.
  Now we let $b\to\infty$ in~\eqref{eq:strong_MP}
  and apply monotone convergence to obtain
  \begin{equation}  \label{eq:after_b_to_infty}
    \limb\,\limyO\frac{1}{\bS(y)}\E^y \eckB{f(Y_t)\1_{t<T_b}}
    =\E^0\eckB{\frac{f(Y_t^\upa)}{\bS(Y_t^\upa)}\1_{t<T_\infty}}.
  \end  {equation}
  The following estimate justifies the interchange of the
  limits $\limb$ and $\limyO$
  \begin{equation}  \begin{split} \label{eq:my_together}
    \lefteqn{\absB{\limyO\frac{1}{\bS(y)}\E^y f(Y_t)
      -\limb\,\limyO\frac{1}{\bS(y)}\E^y\eckb{ f(Y_t)\1_{t<T_b}}}}\\
    &\leq c_f\limb\sup_{y\leq 1}\frac{1}{\bS(y)}
       \E^y \eckB{Y_t\vee Y_t^n \1_{\sup_{s\leq t}Y_s \geq b}}\\
    &\leq c_f\limb\frac{1}{b}\sup_{y\leq 1}\frac{y}{\bS(y)}
       \sup_{y\leq 1}\frac{1}{y}
         \E^y\eckB{ \sup_{s\leq t}\rub{Y_s^2+Y_s^{n+1}}}=0.
  \end  {split}     \end  {equation}
  The last equality follows from $\bS^{'}(0)\in(0,\infty)$
  and from Lemma~\ref{l:finite_moments_Y}.
  Putting \eqref{eq:my_together} and~\eqref{eq:after_b_to_infty} together,
  we get
  \begin{equation}  \begin{split}  \label{eq:after_interchange}
    \limyO\frac{1}{\bS(y)}\E^y \eckb{f(Y_t)}
    =\limb\,\limyO\frac{1}{\bS(y)}\E^y \eckb{f(Y_t)\1_{t<T_b}}
    =\E^0\eckB{\frac{f(Y_t^\upa)}{\bS(Y_t^\upa)}\1_{t<T_\infty}}.
  \end  {split}     \end  {equation}
  Note that~\eqref{eq:after_interchange} is bounded in $t>0$ 
  because of $f(y)\leq c_f y\vee y^n$ and
  Lemma~\ref{l:bS_durch_y_bounded}.

  We finish the proof of the first equality in~\eqref{eq:fixed_time}
  by proving that the limits $\limeps$ and $\limyO$ on the
  right-hand side of~\eqref{eq:fixed_time_eps} interchange.
  \begin{equation}  \begin{split}
    \lefteqn{\absB{\limeps\,\limyO\frac{1}{\bS(y)}\E^y f_\eps (Y_t)
       -\limyO\frac{1}{\bS(y)}\E^y f(Y_t)}}\\
    &\leq\limeps\,\limyO\frac{1}{\bS(y)}\E^y\eckb{f(Y_t)-f_{\eps}(Y_t)}
    =\limeps\,\E^0\eckB{\frac{f(Y_t^\upa)-f_{\eps}(Y_t^\upa)}{\bS(Y_t^\upa)}
       \1_{t<T_\infty}}
    =0.
  \end  {split}     \end  {equation}
  The first equality is~\eqref{eq:after_interchange}
  with $f$ replaced by $f-f_\eps$.
  The last equality follows from the dominated convergence theorem.
  The function $f_\eps/\bS$ converges to $f/\bS$ for every $y>0$ as
  $\eps\to0$. Note that $Y_t^\upa>0$ almost surely for $t>0$.
  Integrability of $\tfrac{f(Y_t^\upa)}{\bS(Y_t^\upa)}\1_{t<T_\infty}$
  follows from finiteness of~\eqref{eq:after_interchange}.
\end  {proof}

We have settled equation~\eqref{eq:zwischentext}
in Lemma~\ref{l:Q_explicit_man_hours}
(together with
Lemma~\ref{l:finite_man_hours}).
A consequence of the finiteness of this equation is that
$\liminf_{t\to\infty}\int \chi_t \,d\barQY=0$. In the proof of
the extinction result for the Virgin Island Model, we will need
that $\int \chi_t \,d\barQY$ converges to zero as $t\to\infty$.
This convergence will follow from equation~\eqref{eq:zwischentext}
if $[0,\infty)\ni t\mapsto \int \chi_t \,d\barQY$ is globally upward
Lipschitz continuous. 
We already know that this function is bound\-ed in $t$ by
Lemma~\ref{l:functional_fixed_t}.

%
%
%
\begin{lemma}  \label{l:integrand_ist_nullfolge}
  Assume~\ref{a:A1}, \ref{a:A2} and \ref{a:S_bar}.
  Let $n\in\N_{\geq1}$.
  If $\int_1^\infty \tfrac{y^n}{g(y)\sb(y)}\,dy<\infty$,
  then
  \begin{equation}
    \limt \int \chi^n_t\,\barQY(d\chi)=0.
  \end  {equation}
\end  {lemma}
\begin{proof}
  We will prove that the function
  $[0,\infty)\ni t\mapsto \int \chi_t^n d\barQY$ is globally upward
  Lipschitz continuous. The assertion then follows from the finiteness
  of~\eqref{eq:Qfunctional} with $f(z)$ replaced by $z^n$ and with $m=1$.
  Recall $\tau_K$, $c_h$ and $c_g$ from the proof of
  Lemma~\ref{l:finite_moments_Y}.
  From~\eqref{eq:def:Y} and It\^o's lemma,
  we obtain for $y>0$ and $0\leq s\leq t$
  \begin{equation}
    \frac{1}{\bS(y)}\E^y\rub{Y^n_{t\wedge\tau_K}}
    -\frac{1}{\bS(y)}\E^y\ruB{Y^n_{s\wedge\tau_K}}
    \leq \ct\int_s^t\frac{1}{\bS(y)}
      \E^y\rub{Y_{r\wedge\tau_K}^n+Y_{r\wedge\tau_K}^{n-1}}\,dr
  \end  {equation}
  where $\ct:= n\rub{c_h+(n-1)c_g}$.
  Letting $K\to\infty$ and then $y\to0$, we conclude from
  the dominated convergence theorem, Lemma~\ref{l:finite_moments_Y}
  and Lemma~\ref{l:functional_fixed_t} that
  \begin{equation}  \label{eq:upward_Lippe}
    \int \chi_t^n - \chi_s^n\,\barQY(d\chi)
    \leq \ct \int_s^t
       \E^0\eckB{\frac{\rub{Y_r^\upa}^n+\rub{Y_r^\upa}^{n-1}}
                      {\bS(Y_r^\upa)}\1_{r<T_\infty}}\,dr
    \leq \ct c_S\abs{t-s}
  \end  {equation}
  for some constant $c_S$.
  The last inequality follows from Lemma~\ref{l:bS_durch_y_bounded}.
  Inequality~\eqref{eq:upward_Lippe} implies upward Lipschitz continuity
  which finishes the proof.
\end  {proof}

%
%
\noindent
\section{Proof of Theorem~\ref{thm:extinction_VIM},
Theorem~\ref{thm:expected_man_hours} and of Theorem~\ref{thm:xlogx}}%
\label{sec:proof_main_theorems}

We will derive Theorem~\ref{thm:extinction_VIM} from
Theorem~\ref{thm:general:extinction} and
Theorem~\ref{thm:expected_man_hours} from
Theorem~\ref{thm:general:expected_man_hours}.
Thus, we need to check that
Assumptions~\ref{a:general:finite_moments}, 
\ref{a:general:finite_man_hours}
and \ref{a:convergence_to_zero}
with $\phi(t,\chi):=\chi_t$,
$\nu:=\Law[x]{Y}$ and $Q:=Q_Y$
hold under~\ref{a:A1}, \ref{a:A2}, \ref{a:S_bar}
and \ref{a:finite_man_hours}.
Recall that $Q_Y=\Sb^{'}(0)\barQY$
and $\bs(0)=\bS^{'}(0)s(0)$.
Assumption~\ref{a:general:finite_moments} follows from
Lemma~\ref{l:finite_moments_Y}
and Lemma~\ref{l:functional_fixed_t}.
Lemma~\ref{l:finite_man_hours} and
Lemma~\ref{l:Q_explicit_man_hours}
imply~\ref{a:general:finite_man_hours}.
Lemma~\ref{l:finite_time_extinction} together with
Lemma~\ref{l:finite_man_hours} implies that $(Y_t)_{t\geq0}$ hits zero in
finite time almost surely.
The second assumption in~\ref{a:convergence_to_zero} is
implied by Lemma~\ref{l:integrand_ist_nullfolge} with $n=1$ and
Assumption~\ref{a:finite_man_hours}.
By Theorem~\ref{thm:general:extinction}, we now know that the total mass
process $(V_t)_{t\geq0}$ dies out if and only if
condition~\eqref{eq:general:bed_fuer_extinction} is satisfied.
However, by Lemma~\ref{l:Q_explicit_man_hours} with $m=1$ and
$f(\cdot)=a(\cdot)$, 
condition~\eqref{eq:general:bed_fuer_extinction} is equivalent to
condition~\eqref{eq:bed_fuer_extinction}.
This proves Theorem~\ref{thm:extinction_VIM}

For an application of Theorem~\ref{thm:general:expected_man_hours}, note
that $f^\nu$ and $f^Q$ are integrable by Lemma \ref{l:finite_man_hours}
and Lemma~\ref{l:Q_explicit_man_hours}, respectively.
In addition, Lemma~\ref{l:finite_man_hours}
and Lemma~\ref{l:Q_explicit_man_hours} show that
\begin{equation*}
  w_{id}(x)=\E^x\int_0^\infty Y_t\,dt=\int_0^\infty f^\nu(t)\,dt
  \text{ and }
  w_a^{'}(0)=\int \rubb{\int_0^\infty a(\chi_s)\,ds} Q_Y(d\chi).
\end  {equation*}
Similar equations hold for $w_{id}^{'}(0)$ and $w_a(x)$.
Moreover, the denominators in~\eqref{eq:man_hours_critical}
and~\eqref{eq:general:man_hours_critical} are equal by
Lemma~\ref{l:Q_explicit_man_hours}, equation~\eqref{eq:QwOI},
together with Lemma~\ref{l:finite_man_hours}.
Therefore, equations~\eqref{eq:man_hours_subcritical}
and~\eqref{eq:man_hours_critical} follow from
equations~\eqref{eq:general:man_hours_subcritical}
and~\eqref{eq:general:man_hours_critical}, respectively.
In the supercritical case, \eqref{eq:supercritical_renewal_theory}
holds because of
\begin{equation}
  \sum_{k=0}^\infty \sup_{k\leq t\leq k+1}e^{-\al t}\int \chi_t Q(d\chi)
  \leq \sup_{t\geq0}\int \chi_t Q(d\chi)\sum_{k=0}^\infty e^{-\al(k+1)}
\end  {equation}
and Lemma~\ref{l:integrand_ist_nullfolge} with $n=1$ together with
Assumption~\ref{a:finite_man_hours}.
Furthermore, Lemma \ref{l:functional_fixed_t}
together with Lemma~\ref{l:bS_durch_y_bounded} and
the dominated convergence theorem
implies continuity of $f^Q$.
Therefore, Theorem~\ref{thm:general:expected_man_hours}
implies~\eqref{eq:general:man_hours_supercritical}
which together with~\eqref{eq:general:rescaled_VIM_supercritical}
reads as~\eqref{eq:man_hours_supercritical}.

Theorem~\ref{thm:xlogx} is a corollary of Theorem~\ref{thm:general:xlogx}.
For this, we need to check \ref{a:for_SKDist}.
The expression in~\eqref{eq:for_SKDist_1} is finite because of
Lemma~\ref{l:functional_fixed_t} with $f(\cdot)=(a(\cdot))^2$ and
Assumption~\ref{a:finite_man_hours_squared}.
Assumption~\ref{a:A1} provides us with $c_1 y\leq a(y)$ for all $y\geq0$
and some $c_1>0$. Thus,
\begin{equation}
  \int_1^\infty \frac{y^2}{g(y)\sb(y)}\,dy
  \leq\frac{1}{c_1}\int_1^\infty a(y)\frac{y+ w_a(y)}{g(y)\sb(y)}\,dy
\end  {equation}
which is finite by~\ref{a:finite_man_hours_squared}.
Lemma~\ref{l:integrand_ist_nullfolge} with $n=2$ and
Lemma~\ref{l:finite_moments_Y} show that $\int \chi_t^2 \barQY(d\chi)$
is bounded in $t\geq0$.
Furthermore, H\"older's inequality implies
\begin{equation}  \label{eq:hoelders_inequality}
  \ruB{\int \eckB{\chi_t\int_0^t a(\chi_s)\,ds}\barQY(d\chi)}^2
  \leq\int \chi_t^2\,\barQY(d\chi)
     \int\ruB{\int_0^\infty a(\chi_s)\,ds}^2\barQY(d\chi)
\end  {equation}
which is bounded in $t\geq0$ because
of Lemma~\ref{l:Q_explicit_man_hours} with $m=2$, $f(\cdot)=a(\cdot)$
and because of Assumption~\ref{a:finite_man_hours_squared}.
Therefore, we may apply Theorem~\ref{thm:general:xlogx}. Note that the
limiting variable is not identically zero because of
\begin{equation}
  \int\ruB{ A_\al \log^+(A_\al)}dQ\leq\int\rub{A_\al}^2 dQ\leq
  \int\ruB{\int_0^\infty a(\chi_s)\,ds}^2 \barQY(d\chi)<\infty.
\end  {equation}
The right-hand side is finite
because of Lemma~\ref{l:Q_explicit_man_hours} with $m=2$, $f(\cdot)=a(\cdot)$
and because of Assumption~\ref{a:finite_man_hours_squared}.

\bigskip
\noindent
\textbf{\large{Acknowledgments:}}
We thank Anton Wakolbinger and Don Dawson
for inspiring discussions and valuable comments.
Also we thank two anonymous referees for their helpful comments and suggestions.

\hyphenation{Sprin-ger}

\end{document}